\documentclass{amsart}
\usepackage{graphicx}
\newtheorem{theorem}{Theorem}[section]
\newtheorem{corollary}[theorem]{Corollary}
\newtheorem{lemma}[theorem]{Lemma}
\newtheorem{proposition}[theorem]{Proposition}
\theoremstyle{definition}
\newtheorem{definition}[theorem]{Definition}
\theoremstyle{example}
\newtheorem{example}[theorem]{Example}
\theoremstyle{remark}
\newtheorem{remark}{Remark}[theorem]
\numberwithin{equation}{section}

\begin{document}
\title{Dynamically Consistent Nonlinear Evaluations and Expectations}%
\author{Shige Peng}%
\address{Institute of Mathematics, Shandong University,\\
Jinan, 250100, China \\
}%
\email{Email: peng@sdu.edu.cn}%

\thanks{This research is supported in part by The National Natural Science
Foundation of China No. 10131040. This reversion is made after the
author's visit, during November 2003, to Institute of Mathematics
and System Science, Academica Sinica, where he gives a series of
lectures on this paper. He thanks Zhiming Ma and Jia-an Yan, for
their fruitful suggestions, critics and warm encouragements. He
also thanks to Claude Dellacherie for his
suggestions and critics. }%

\date{First version: 16 August 2003,
This version 31 March 2004}%
\begin{abstract}
{\small How an agent (or a firm, an investor, a financial market)
evaluates
a contingent claim, say a European type of derivatives $X$, with maturity $t$%
? In this paper we study a dynamic evaluation  of this problem. We
denote by $\left\{ \mathcal{F}_t\right\} _{t\geq 0}$, the
information acquired  by this agent. The value $X$ is known at the
maturity $t$ means
that $X$ is an $\mathcal{F}_t$--measurable random variable. We denote by $%
\mathcal{E}_{s,t}[X]$ the evaluated value of $X$ at the time $s\leq t$. $%
\mathcal{E}_{s,t}[X]$ is $\mathcal{F}_s$--measurable since his
evaluation is based on his information at the time $s$. Thus
$\mathcal{E}_{s,t}[\cdot ]$ is an operator that maps an
$\mathcal{F}_t$--measurable random variable to
an $\mathcal{F}_s$--measurable one. A system of operators $\left\{ \mathcal{E%
}_{s,t}[\cdot ]\right\} _{0\leq s\leq t<\infty }$ is called $\mathcal{F}_t$%
--consistent evaluations if it satisfies the following conditions: (A1) $%
\mathcal{E}_{s,t}[X]\geq \mathcal{E}_{s,t}[Y]$, if $X\geq Y$; (A2) $\mathcal{%
E}_{t,t}[X]=X$; (A3) $\mathcal{E}_{r,s}\mathcal{E}_{s,t}[X]=\mathcal{E}%
_{r,t}[X]$, for $r\leq s\leq t$; (A4) $1_A\mathcal{E}_{s,t}[X1_A]=1_A%
\mathcal{E}_{s,t}[X]$, if $A\in \mathcal{F}_s$. }

{\small In the situation where $\mathcal{F}_t$ is generated by a
Brownian
motion, we propose the so-called $g$--evaluation defined by $\mathcal{E}%
_{s,t}^g[X]:=y_s $, where $y$ is the solution of the backward
stochastic
differential equation with generator $g$ and with the terminal condition $%
y_t=X$. This $g$--evaluation satisfies (A1)--(A4). We also provide
examples to determine the function $g=g(y,z)$ by testing. }

{\small The main result of this paper is as follows: if a given $\mathcal{F}%
_t$--consistent evaluation is $\mathcal{E}^{g_\mu }$--dominated, i.e., (A5) $%
\mathcal{E}_{s,t}[X]-\mathcal{E}_{s,t}[X^{\prime }]\leq
\mathcal{E}^{g_\mu
}[X-X^{\prime }]$, for a large enough $\mu >0$, where $g_\mu =\mu (|y|+|z|)$%
, then $\mathcal{E}_{s.t}[\cdot ]$ is a $g$--evaluation}
\end{abstract}
\maketitle

\tableofcontents

\section{Introduction\label{ss1}}

We are interested in the following dynamically consistent evaluation of
risky assets: Let $\eta =(\eta _t)_{t\geq 0}$ be a $d$--dimensional process,
it may be the prices of stocks in a financial market, the rates of
exchanges, the rates of local and global inflations etc. We assume that at
each time $t\geq 0$, the information for of an agent (a firm, a group of
people or even a financial market) is the history of $\eta $ during the time
interval $[0,t]$. Namely, his actual filtration is
\[
\mathcal{F}_t=\sigma \{\eta _s;s\leq t\}.
\]
We denote the set of all real valued $\mathcal{F}_t $--measurable random
variables by $m\mathcal{F}_t $. Under this notation an $\eta $--underlying
derivative $X$, with maturity $T\in [0,\infty )$, is an $\mathcal{F}_T $%
--measurable random variable, i.e., $X\in m\mathcal{F}_T$. We denote this
evaluated value at the time $t$ by $\mathcal{E}_{t,T}[X]$. It is reasonable
to assume that $\mathcal{E}_{t,T}[X]$ is $\mathcal{F}_t $--measurable. In
other words,
\[
\mathcal{E}_{t,T}[X]:m\mathcal{F}_T\longrightarrow m\mathcal{F}_t.
\]
In particular
\[
\mathcal{E}_{0,T}[X]:m\mathcal{F}_T\longrightarrow \mathbf{R}.
\]

We will make the following \textbf{Axiomatic Conditions}\ for $(\mathcal{E}%
_{t,T}[\cdot ])_{0\leq t\leq T<\infty }$:

\smallskip\

\noindent (A1) Monotonicity: $\mathcal{E}_{t,T}[X]\geq \mathcal{E}%
_{t,T}[X^{\prime }],\;$if $X\geq X^{\prime };$

\noindent (A2)\ \ \ $\mathcal{E}_{T,T}[X]=X,$\ $\forall X\in m\mathcal{F}_T $%
. Particularly $\mathcal{E}_{0,0}[c]=c$;\

\noindent (A3) Dynamical consistency:\textbf{\ }$\mathcal{E}_{s,t}[\mathcal{E%
}_{t,T}[X]]=\mathcal{E}_{s,T}[X]$, if $s\leq t\leq T$;

\noindent (A4) ``Zero--one law'':\textbf{\ }for each $t\leq T$, $1_A\mathcal{%
E}_{t,T}[X]=1_A\mathcal{E}_{t,T}[1_AX],\;\;\forall A\in \mathcal{F}_t.$%
\textbf{\ }

or, more specially,

\noindent (A4') \textbf{\ }for each $t\leq T$, $1_A\mathcal{E}_{t,T}[X]=%
\mathcal{E}_{t,T}[1_AX],\;\;\forall A\in \mathcal{F}_t.$\textbf{\ }

\smallskip\

\begin{remark}
The meaning of (A1) and (A2) are obvious. Condition (A3) means that the
evaluated value $\mathcal{E}_{t,T}[X]$ can be also treated as a derivative
with the maturity $t$. At a time $s\leq t$, the ``price'' of this derivative
evaluated by $\mathcal{E}_{s,t}[\mathcal{E}_{t,T}[X]]$ is the same as the
``price'' of the original derivative $X$ with maturity $T$, i.e., $\mathcal{E%
}_{s,T}[X]$.
\end{remark}

\begin{remark}
The meaning of condition (A4) is: at time $t$, the agent knows whether $\eta
_{\cdot \wedge t}$ is in $A$. If $\eta _{\cdot \wedge t}$ is in $A$, then
the value $\mathcal{E}_{t,T}[X]$ is the same as $\mathcal{E}_{t,T}[1_AX]$
since the two outcomes $X$ and $1_AX$ are exactly the same.
\end{remark}

\smallskip

It is clear that, to investigate this abstract evaluation problem, we need
to introduce some regulation condition of $\mathcal{E}$. In this paper the
information $\mathcal{F}_t$ will be limited to the $\sigma $--filtration of
some $d$--dimensional Brownian motion, and $X$ will be assumed to be
square--integrable, i.e., $X\in L^2(\mathcal{F}_T)$.

A condition stronger than (A2) is:

(A2')\ \ \ $\mathcal{E}_{s,t}[X]=X,$\ $\forall 0\leq s\leq t$, $\forall X\in
m\mathcal{F}_s $.

The meaning is that the market has zero interest rate for a non--risky asset
$X$. In this case we can define $\mathcal{E}[X|\mathcal{F}_t]:=\mathcal{E}%
_{t,T}[X]$, for a sufficiently large $T$, and $\mathcal{E}[X]:=\mathcal{E}[X|%
\mathcal{F}_0]$. It is easy to check that
\[
\mathcal{E}[1_A\mathcal{E}[X|\mathcal{F}_t]]=\mathcal{E}[1_AX].
\]
$\mathcal{E}[X|\mathcal{F}_t]$ is called the $\mathcal{E}$--conditional
expectation of $X$ under $\mathcal{F}_t$. It satisfies all properties of a
classical expectation, with one exception that it can be a nonlinear
operator. $\{\mathcal{E}[X|\mathcal{F}_t]\}_{0\leq t\leq T}$ is called an $%
\mathcal{F}_t$--consistent nonlinear expectation.

A typical filtration-consistent nonlinear expectation, called $g$%
-expectation and denoted by $\{\mathcal{E}_g[X|\mathcal{F}_t]\}_{0\leq t\leq
T}$, was introduced in \cite[Peng1997]{Peng1997}. A significant feature of
this $g$--expectation is that the value of $\mathcal{E}_g[X|\mathcal{F}_t]$
is uniquely determined by a simple function $g(t,y,z)$ with $g(t,y,0)\equiv
0 $. In fact $(\mathcal{E}_g[X|\mathcal{F}_t])_{0\leq t\leq T}$ is the
solution of the backward stochastic differential equation (BSDE in short)
with the function $g$ as its generator and with $X$ as its terminal
condition at the terminal time $T$. It is then not surprising that the
behavior of $\mathcal{E}_g[\cdot |\mathcal{F}_t]$ is entirely characterized
by this concrete function $g$. For example, $\mathcal{E}_g[\cdot |\mathcal{F}%
_t]$ is a linear (conventional) expectation if and only if $g$ is
independent of $y$ and is a linear function of $z$, i.e., $g$ has a form $%
g=b_t\cdot z$; $\mathcal{E}_g[X|\mathcal{F}_t]$ is concave (resp. convex) in
$X$ if and only if $g$ is concave (resp. convex) in $(y,z)$, etc. For an
interesting application of $g$--expectations to the utility in stochastic
continuous--time setting with ambiguity (or ``model uncertainty'' referred
by Hansen and Sargent and Anderson, Hansen and \cite[Sargent]{AHS}, see
\cite[Chen and Epstein, 2002]{Chen-Epstein2002}.

$g$--expectations also have very interesting mathematical properties. A
nonlinear Doob--Meyer's decomposition theorem for $g$--supermartingales was
obtained by \cite[Peng, 1999]{Peng1999}, for the case of Brownian
filtration, and then by Chen and Peng 1998 \cite[Chen \& Peng, 1998]
{Chen-Peng1998} for a general filtration. In the case where the assumption $%
g(t,y,0)\equiv 0$ does not hold, we have to denote the solution of the
related BSDE by $\mathcal{E}_{s,t}^g[X]$ instead of $\mathcal{E}_g[X|%
\mathcal{F}_s]$. $\{\mathcal{E}_{s,t}^g[\cdot ]\}_{0\leq s\leq t\leq T}$
satisfies (A1)--(A4). $\mathcal{E}_{s,t}^g[\cdot ]$ can be applied to a
wider situation in economics and finance.

The application of BSDE to the pricing of contingent claims in a financial
market was studied in \cite[El Karoui et al., 1997]{EPQ1997}. Most of the
results in \cite{EPQ1997} can be interpreted in the language of $\mathcal{E}%
_{s,t}^g[\cdot ]$. Other recent results in $g$--expectations are in \cite
{BCHMP00}, \cite{chen98}, \cite{CHMP2002}, \cite{Chen-Peng1998}, \cite
{Chen-Peng2001}, \cite{Chen-Peng2000}, \cite{Peng1997}, \cite{Peng1999},
\cite{Peng2002}, \cite{Peng2003a}, \cite{Peng2003b} where some cases are
studied in depth. For nonlinear evaluations, see \cite[Peng, 2002]{Peng2002}%
, \cite[Peng, 2003]{Peng2003a} and \cite[Peng, 2003]{Peng2003b}).

An interesting problem is: are the notions of $g$--expectations and $g$%
-evaluations general enough to represent all ``enough regular''
filtration-consistent nonlinear expectations and evaluations? In this case
we can then concentrate ourselves to find the corresponding function $g$
which determine entirely the evaluation.

For the case of filtration--consistent expectations, we have partially
solved the problem in \cite{CHMP2002}: If the assumptions (A1)--(A4) plus
(A2') hold and if for a large enough $\mu >0$, the nonlinear expectation $%
\mathcal{E}[\cdot ]$ is dominated by the `$g^\mu $--expectation' $\mathcal{E}%
^{g^\mu }[\cdot ]$ with $g=\mu |z|$, and furthermore, if $\mathcal{E}[X+\eta
|\mathcal{F}_t]=\mathcal{E}[X|\mathcal{F}_t]+\eta $ for all $\mathcal{F}_t$%
-measurable $\eta $, then, there exists a unique $g$, independent of $y$,
such that $\mathcal{E}[\cdot ]=\mathcal{E}_g[\cdot ]$.

The main objective of this paper is to prove this problem for the general
case of filtration--consistent evaluation: (see Theorem \ref{m7.1}) if a
filtration consistent evaluation $\{\mathcal{E}_{s,t}[\cdot ]\}_{0\leq s\leq
t\leq T}$ satisfies (A1)--(A4) plus the corresponding $\mathcal{E}^{g_\mu }$%
--dominated conditions (see (A5) in Section \ref{ss3}), then the mechanism
of this seemingly very abstract evaluation $\mathcal{E}_{s,t}[\cdot ]$ can
be entirely determined by a simple function $g(t,y,z)$. This means that,
there exists a unique function $g$ such that, for each $0\leq s\leq t\leq T$
and for each $X\in L^2(\mathcal{F}_t)$, the value $\mathcal{E}_{s,t}[X]$ is
the solution of the BSDE with generator $g$ and with terminal condition $X$.
The result of this paper have non trivially generalized our previous result
of \cite{CHMP2002}: condition (A2') and ``$\mathcal{E}_{t,T}[X+\theta
]=\theta $,\ $\forall X\in m\mathcal{F}_T$ and $\eta \in m\mathcal{F}_t$ ''
are not at all required.

It is worth to point out that the well--known Black--Scholes option pricing
formula is a case where $g$ is a linear function. But in our axiomatic
condition (A1)--(A4) as well as the regularity condition (A5), neither the
arbitrage free condition, which is a principle argument in Black--Scholes
theory, nor utility maximization has been involved. Another point is that
the model of the price $\eta $ of the underlying stocks is not specified.
This gives us a large freedom to determine the function $g$ in each specific
situation. We also explain how the function $g$ can be determined by simply
testing the agent's evaluation. This testing method is very useful to
determine an agent's behavior under risk.

The paper is organized as follows: in section \ref{ss2}, we give a rigorous
setting of the notion of $\mathcal{F}_t$--consistent nonlinear evaluation
and its special case: $\mathcal{F}_t$--expectations in subsection \ref{ss2.1}%
. We then give a concrete $\mathcal{F}_t$--evaluation: $\mathcal{E}^g$%
--evaluations in subsection \ref{ss2.2}. The main result, Theorem \ref{m7.1}%
, will be presented in section \ref{ss3}. We also provide some examples and
explain how to find the function $g$ through by testing the input--output
data. This main theorem will be proved in Section \ref{ss8}, with several
propositions served as lemmas for the proof given in Sections \ref{ass3}--%
\ref{ss9}. Although the whole paper is focused to prove Theorem
\ref{m7.1}, many preparative results of this paper have their own
interests, e.g., the existence and uniqueness of BSDE under
$\mathcal{E}$ (Theorem \ref{m5.1}); the new nonlinear
supermartingale decomposition theorem of Doob--Meyer's type
(Theorem \ref{m6.1}), using a new and intrinsic formulation (see
Remark \ref{m6.1Rem2} after Theorem \ref{m6.1}). This
decomposition theorem has also the extension of
$\mathcal{E}_{s,t}[\cdot ]$ to $\mathcal{E}_{\sigma ,\tau }[\cdot
]$ with stopping times $\sigma $ and $\tau $ and the related
optional stopping theorem (Theorem \ref{m8.13} and Theorem
\ref{m8.14}). Mathematically, some of them are more fundamental
than Theorem \ref{m7.1}. In particular, the nonlinear
decomposition theorems of Doob--Meyer's type, i.e., Proposition
\ref{p2.3} and Theorem \ref{m8.1} play crucial roles in the proof
of Theorem \ref{m7.1}. Theorem \ref{m8.1} has also an intersting
interpretation in finance (see Remark \ref{m6.1Rem1}).

Another application of the dynamical expectations and evaluations
is to risk measures. Axiomatic conditions for a (one step)
coherent risk measure was introduced by  Artzner, Delbaen, Eber
and Heath 1999 \cite{ADEH1999} and, for a convex risk measure, by
F\"ollmer and Schied (2002) cite{Fo-Sc}. Rosazza Gianin (2003)
studied dynamical risk measures using the notion of
$g$--expectations in \cite{rosazza} (see also \cite{Peng2003b}) in
which  (B1)--(B4) are satisfied. In fact conditions (A1)-(A4), as
well as their special situation (B1)--(B4) (see Proposition
\ref{m2a2}) provides
 an ideal characterization  of the dynamical behaviors of a the
a risk measure. But in this paper we emphasis the study of the
mechanism of the evaluation to a further payoff, for which is, in
general,  the translation property in risk measure is not
satisfied.

\section{Basic setting and $\mathcal{E}^g$--evaluations by BSDE\label{ss2}}

\subsection{Basic setting\label{ss2.1}}

Let $(\Omega ,\mathcal{F},P)$ be a probability space and let $(B_t)_{t\geq
0} $ be a $d$--dimensional Brownian motion defined in this space. We denote
by $(\mathcal{F}_t)_{t\geq 0}$ the natural filtration generated by $B$,
i.e.,
\[
\mathcal{F}_t:=\sigma \{\sigma \{B_s,\;s\leq t\}\cap \mathcal{N}\}.
\]
Here $\mathcal{N}$ is the collection of all $P$--null subsets. For each $%
t\in [0,\infty )$, we denote by

\begin{itemize}
\item  $L^2(\mathcal{F}_t):=$\{the space of all real valued $\mathcal{F}_t$%
--measurable random variables such that $E[|\xi |^p]<\infty $\}.
\end{itemize}

\begin{definition}
\label{d2.1}A system of operators:
\[
\mathcal{E}_{s,t}[X]:X\in L^2(\mathcal{F}_t)\rightarrow L^2(\mathcal{F}%
_s),\;T_0\leq s\leq t\leq T_1
\]
is called an $\mathcal{F}_t$--\textrm{consistent nonlinear evaluation}
defined on $[T_0,T_1]$ if it satisfies the following properties: for each $%
T_0\leq r\leq s\leq t$ and for each $X$, $X^{\prime }\in L^2(\mathcal{F}_t)$%
, \\ \textbf{(A1)} $\mathcal{E}_{s,t}[X]\geq \mathcal{E}_{s,t}[X^{\prime }]$%
, a.s., if $X\geq X^{\prime }$, a.s.; \\ \textbf{(A2)} $\mathcal{E}%
_{t,t}[X]=X$, a.s.; \\ \textbf{(A3)} $\mathcal{E}_{r,s}[\mathcal{E}%
_{s,t}[X]]=\mathcal{E}_{r,t}[X]$, a.s.; \\ \textbf{(A4)} $1_A\mathcal{E}%
_{s,t}[X]=1_A\mathcal{E}_{s,t}[1_AX]$, a.s.\ $\forall A\in \mathcal{F}_s$.
\end{definition}

We will often consider (A1)--(A4) plus an additional condition: \\\textbf{(A4%
}$_0$\textbf{)} $\mathcal{E}_{s,t}[0]=0$, a.s.\ $\forall 0\leq s\leq t\leq T$%
.

\begin{proposition}
\label{m2a4}(A4) plus (A4$_0$) is equivalent to \\ \textbf{(A4')} $1_A%
\mathcal{E}_{s,t}[X]=\mathcal{E}_{s,t}[1_AX]$, a.s.\ $\forall A\in \mathcal{F%
}_s$.
\end{proposition}
\smallskip\noindent\textbf{Proof. }
It is clear that (A4') implies (A4). $\mathcal{E}%
_{s,t}[0]\equiv 0$ can be derived by putting $A=\emptyset $ in (A4'). On the
other hand, (A4) plus the additional condition implies
\[
1_{A^C}\mathcal{E}_{s,t}[1_AX]=1_{A^C}\mathcal{E}_{s,t}[1_{A^c}1_AX]=0.
\]
We thus have
\begin{eqnarray*}
\mathcal{E}_{s,t}[1_AX] &=&1_{A^C}1_A\mathcal{E}_{s,t}[X]+1_A1_A\mathcal{E}%
_{s,t}[X] \\
&=&1_A\mathcal{E}_{s,t}[X].
\end{eqnarray*}
$\Box $\medskip\

\begin{proposition}
\label{mA4} (A4) is equivalent to, for each $0\leq s\leq t$ and $X,X^{\prime
}\in L^2(\mathcal{F}_t)$,
\begin{equation}
\mathcal{E}_{s,t}[1_AX+1_{A^C}X^{\prime }]=1_A\mathcal{E}_{s,t}[X]+1_{A^C}%
\mathcal{E}_{s,t}[X^{\prime }],\;\mathrm{a.s.}\;\forall A\in \mathcal{F}_s.
\label{eA4}
\end{equation}
\end{proposition}

\smallskip\noindent\textbf{Proof. }(A4) $\Rightarrow $ (\ref{eA4}): We let $%
Y=1_AX+1_{A^C}X^{\prime }$. Then, by (A4)
\[
1_A\mathcal{E}_{s,t}[Y]=1_A\mathcal{E}_{s,t}[1_AY]=1_A\mathcal{E}%
_{s,t}[1_AX]=1_A\mathcal{E}_{s,t}[X].
\]
Similarly
\[
1_{A^C}\mathcal{E}_{s,t}[Y]=1_{A^C}\mathcal{E}_{s,t}[1_{A^C}Y]=1_{A^C}%
\mathcal{E}_{s,t}[1_{A^C}X^{\prime }]=1_{A^C}\mathcal{E}_{s,t}[X^{\prime }].
\]
Thus (\ref{eA4}) from $1_A\mathcal{E}_{s,t}[Y]+1_{A^C}\mathcal{E}%
_{s,t}[Y]=1_A\mathcal{E}_{s,t}[X]+1_{A^C}\mathcal{E}_{s,t}[X^{\prime }]$.

(\ref{eA4}) $\Rightarrow $ (A4): It is simply because of
\begin{eqnarray*}
1_A\mathcal{E}_{s,t}[1_AX] &=&1_A\mathcal{E}_{s,t}[1_AX+1_{A^C}(1_AX)] \\
&=&1_A(1_A\mathcal{E}_{s,t}[X]+1_{A^C}\mathcal{E}_{s,t}[1_{A^C}X]) \\
&=&1_A\mathcal{E}_{s,t}[X].
\end{eqnarray*}
$\Box $\medskip\

\begin{remark}
At time $t$, the agent knows the value of $1_A$. (A4) means that,
if $1_A=1$ then the evaluated value $\mathcal{E}_{s,t}[1_AX]$
should be the same
as $\mathcal{E}_{s,t}[X]$ since the two outcomes $X(\omega )$ and $%
(1_AX)(\omega )$ are exactly the same. (A4) is applied to the
evaluation of a final outcome $X$ plus some ``dividend''
$(D_s)_{s\geq 0}$.
\end{remark}

If, instead of (A2), we set

\medskip
\noindent\textbf{(A2')} $\mathcal{E}_{s,t}[X]=X$, a.s., for each $T_0\leq
s\leq t\leq T$, and $X\in L^2(\mathcal{F}_s)$. \medskip

\noindent Then we define
\[
\mathcal{E}[X|\mathcal{F}_t]:=\mathcal{E}_{t,T}[X],\;X\in L^2(\mathcal{F}%
_T).
\]
We observe that this notion describes all $\mathcal{E}_{s,t}[X]$ since, when
$X\in L^2(\mathcal{F}_t)$, $\mathcal{E}_{s,t}[X]=\mathcal{E}[X|\mathcal{F}%
_t] $.

\begin{proposition}
\label{m2a2}With (A2'), the system of operators
\[
\mathcal{E}[\cdot |\mathcal{F}_t]:L^2(\mathcal{F}_T)\rightarrow L^2(\mathcal{%
F}_t)
\]
is a $\mathcal{F}_t$--\textrm{consistent nonlinear expectation}, i.e., it
satisfies, for each $T_0\leq s\leq t\leq T$, and $X\in L^2(\mathcal{F}_T)$\\
\textbf{(B1)} $\mathcal{E}[X|\mathcal{F}_t]\geq \mathcal{E}[X|\mathcal{F}_t]$%
, a.s., if $X\geq X^{\prime }$, a.s. \\ \textbf{(B2)} $\mathcal{E}[X|%
\mathcal{F}_t]=X$, a.s., if $X\in L^2(\mathcal{F}_t)$; \\ \textbf{(B3)} $%
\mathcal{E}[\mathcal{E}[X|\mathcal{F}_t]|\mathcal{F}_s]=\mathcal{E}[X|%
\mathcal{F}_s]$, a.s.; \\ \textbf{(B4)} $\mathcal{E}[1_AX|\mathcal{F}_t]=1_A%
\mathcal{E}[X|\mathcal{F}_t]$, a.s.\ $\forall A\in \mathcal{F}_t$.
\end{proposition}

\smallskip\noindent\textbf{Proof. }(B1)--(B3) are easy. Since (A2') implies $\mathcal{E%
}_{s,t}[0]=0$, thus, by Proposition \ref{m2a4}, (A4') and then (B4) holds. $%
\Box $\medskip\

We have the following immediate result

\begin{proposition}
\label{m2.3}Let $T_0<T_1<T_2<\cdots <T_N$ be given and, for $i=0,1,2,\cdots
,N-1$, let
\[
\mathcal{E}_{s,t}^i[X]:X\in L^2(\mathcal{F}_t)\rightarrow L^2(\mathcal{F}%
_s),\;T_i\leq s\leq t\leq T_{i+1}
\]
be an $\mathcal{F}_t$--consistent evaluation defined on $[T_i,T_{i+1}]$ in
the sense of Definition \ref{d2.1}. Then there exists a unique $\mathcal{F}_t
$--consistent evaluation $\mathcal{E}[\cdot ]$ defined on $[T_0,T_N]$
\[
\mathcal{E}_{s,t}[X]:X\in L^2(\mathcal{F}_t)\rightarrow L^2(\mathcal{F}%
_s),\;T_0\leq s\leq t\leq T_N
\]
such that, for each $i=0,1,\cdots ,N-1$, and for each $T_i\leq s\leq t\leq
T_{i+1}$,
\begin{equation}
\mathcal{E}_{s,t}[X]=\mathcal{E}_{s,t}^i[X],\;\forall X\in L^2(\mathcal{F}%
_t).  \label{e2.01}
\end{equation}
\end{proposition}

\smallskip\noindent\textbf{Proof. }It suffices to prove the case $N=2$. Because after
we then can apply this result to prove the cases $[T_0,T_3]=[T_0,T_2]\cup
[T_2,T_3]$, $\cdots $, $[T_0,T_N]=[T_0,T_{N-1}]\cup [T_{N-1},T_N]$. We
define
\begin{equation}
\mathcal{E}_{s,t}[X]=\left\{
\begin{array}{cll}
\hbox{(i)} & \mathcal{E}_{s,t}^1[X]\hbox{,} & T_0\leq s\leq t\leq T_1; \\
\hbox{(ii)} & \mathcal{E}_{s,t}^2[X], & T_1\leq s\leq t\leq T_2; \\
\hbox{(iii)} & \mathcal{E}_{s,T_1}^1[\mathcal{E}_{T_1,t}^2[X]] &
T_1\leq s<T_1<t\leq T_2.
\end{array}
\right.  \label{e2.02}
\end{equation}
It is clear that, on $[T_0,T_2]$, $\mathcal{E}_{s,t}[\cdot ]$ satisfies (A1)
and (A2). To prove (A3) we only need to check the relation
\[
\mathcal{E}_{r,s}[\mathcal{E}_{s,t}[X]]=\mathcal{E}_{r,t}[X],\;T_0\leq r\leq
s\leq t\leq T_1
\]
for two cases: $T_0\leq r\leq s\leq T_1\leq t\leq T_2$ and $T_0\leq r\leq
T_1\leq s\leq t\leq T_2$. For the first case
\begin{eqnarray*}
\mathcal{E}_{r,s}[\mathcal{E}_{s,t}[X]] &=&\mathcal{E}_{r,s}^1[\mathcal{E}%
_{s,T_1}^1[\mathcal{E}_{T_1,t}^2[X]]] \\
&=&\mathcal{E}_{r,T_1}^1[\mathcal{E}_{T_1,t}^2[X]] \\
&=&\mathcal{E}_{r,t}[X].
\end{eqnarray*}
For the second case
\begin{eqnarray*}
\mathcal{E}_{r,s}[\mathcal{E}_{s,t}[X]] &=&\mathcal{E}_{r,T_1}^1[\mathcal{E}%
_{T_1,s}^2[\mathcal{E}_{s,t}^2[X]]] \\
\ &=&\mathcal{E}_{r,T_1}^1[\mathcal{E}_{T_1,t}^2[X]] \\
\ &=&\mathcal{E}_{r,t}[X].
\end{eqnarray*}

We now prove (A4). Again it suffices to check the case $T_0\leq s\leq
T_1\leq t\leq T_2$. In this case, for each $A\in \mathcal{F}_s\subset
\mathcal{F}_{T_1}$, (A4) is derived from
\begin{eqnarray*}
1_A\mathcal{E}_{s,t}[X] &=&1_A\mathcal{E}_{s,T_1}^1[\mathcal{E}_{T_1,t}^2[X]]
\\
&=&1_A\mathcal{E}_{s,T_1}^1[1_A\mathcal{E}_{T_1,t}^2[X]] \\
&=&1_A\mathcal{E}_{s,T_1}^1[\mathcal{E}_{T_1,t}^2[1_AX]] \\
&=&1_A\mathcal{E}_{s,t}[1_AX].
\end{eqnarray*}

It remains to prove the uniqueness of $\mathcal{E}[\cdot ]$. Let $\mathcal{E}%
^a[\cdot ]$ be an $\mathcal{F}_t$--consistent evaluation such that,\
\[
\mathcal{E}_{s,t}^a[X]=\mathcal{E}_{s,t}^i[X],\;\forall X\in L^2(\mathcal{F}%
_t),\;i=1,2.
\]
We then have, when $T_0\leq s\leq t\leq T_1$ and $T_1\leq s\leq t\leq T_2$, $%
\mathcal{E}_{s,t}^a[X]\equiv \mathcal{E}_{s,t}[X],\;\forall X\in L^2(%
\mathcal{F}_t)$. For the remaining case, i.e., $T_0\leq s<T_1<t\leq T_1$,
since $\mathcal{E}^a$ satisfies (A3),
\begin{eqnarray*}
\mathcal{E}_{s,t}^a[X] &=&\mathcal{E}_{s,T_1}^a[\mathcal{E}_{T_1,t}^a[X]] \\
&=&\mathcal{E}_{s,T_1}^1[\mathcal{E}_{T_1,t}^2[X]] \\
&=&\mathcal{E}_{s,t}[X],\;\forall X\in L^2(\mathcal{F}_t).
\end{eqnarray*}
Thus $\mathcal{E}_{s,t}^a[\cdot ]=\mathcal{E}_{s,t}[\cdot ]$. This completes
the proof. $\Box $\medskip\

\begin{remark}
\label{m2.3r}(i) In the remaining of this paper, we mainly consider the
situation $t\in [0,T]$ for a fixed $T$. The conclusions can be extended to $%
[0,\infty )$, using the above Proposition. (ii) The argument of the above
Proposition \ref{m2.3} can be also applied to a filtration different from $\{%
\mathcal{F}_t\}_{t\geq 0}$, e.g., $\{\mathcal{F}_{t\wedge \tau }\}_{t\geq 0}$%
, where $\tau $ is an $\mathcal{F}_t$--stopping time.
\end{remark}

\subsection{$\mathcal{E}^g$--evaluations induced by BSDE\label{ss2.2}}

In the remaining of this paper, we limited ourselves within the time
interval $[0,T]$ for some fixed $T>0$. The results of this paper can be
extended to the whole interval $[0,\infty )$, using Proposition \ref{m2.3}.
We need the following notations. Let $p\geq 1$ and $\tau \leq T$ be a given $%
\mathcal{F}_t$--stopping time.

\begin{itemize}
\item  $L^p(\mathcal{F}_\tau ;R^m):=$\{the space of all $R^m$--valued $%
\mathcal{F}_\tau $--measurable random variables such that $E[|\xi
|^p]<\infty $\};

\item  $L_{\mathcal{F}}^p(0,\tau ;R^m):=$\{$R^m$--valued and $\mathcal{F}_t$%
--predictable stochastic processes such that $E\int_0^\tau |\phi
_t|^pdt<\infty $\};

\item  $D_{\mathcal{F}}^p(0,\tau ;R^m):=$\{all RCLL processes in $L_{%
\mathcal{F}}^p(0,\tau ;R^m)$ such that $E[\sup_{0\leq t\leq \tau }|\phi
_t|^p]<\infty $\};

\item  $S_{\mathcal{F}}^p(0,\tau ;R^m):=$\{all continuous processes in $D_{%
\mathcal{F}}^p(0,\tau ;R^m)$ \};

\item  $\mathcal{S}_T:=$\{the collection of all $\mathcal{F}_t$--stopping
times bounded by $T$\};

\item  $\mathcal{S}_T^0:=$\{$\tau \in \mathcal{S}_T$ and $\cup
_{i=1}^n\{\tau =t_i\}=\Omega $, with some deterministic $0\leq t_1<\cdots
<t_N$\}.
\end{itemize}

In the case $m=1$, we denote them by $L^p(\mathcal{F}_\tau )$, $L_{\mathcal{F%
}}^p(0,\tau )$, $D_{\mathcal{F}}^p(0,\tau )$ and $S_{\mathcal{F}}^p(0,\tau )$%
. We recall that all elements in $D_{\mathcal{F}}^2(0,T)$ are $\mathcal{F}_t$%
--predictable.

For each given $t\in [0,T]$ and $X\in L^2(\mathcal{F}_t)$, we solve the
following BSDE

\begin{equation}
Y_s=X+\int_s^tg(r,Y_r,Z_r)dr-\int_s^tZ_rdB_r,\;s\in [0,t],  \label{tsBSDE}
\end{equation}
where the unknown is the pair of the adapted processes $(Y,Z)$. Here the
function
\[
g:(\omega ,t,y,z)\in \Omega \times [0,T]\times R\times R^d\rightarrow R
\]
satisfies the following basic assumptions for each $\forall y,y^{\prime }\in
R,\;z,z^{\prime }\in R^d$
\begin{equation}
\left\{
\begin{array}{rrl}
\hbox{(i)} &  & g(\cdot ,y,z)\in L_{\mathcal{F}}^2(0,T)\hbox{;} \\
\hbox{(ii)} &  & |g(t,y,z)-g(t,y^{\prime },z^{\prime })|\leq \mu
(|y-y^{\prime }|+|z-z^{\prime }|)\;.
\end{array}
\right.  \label{h2.1}
\end{equation}
In some cases it is interesting to consider the following situation:
\begin{equation}
\left\{
\begin{array}{rrl}
\hbox{(a)} & g(\cdot ,0,0) & \equiv 0, \\
\hbox{(b)} & g(\cdot ,y,0) & \equiv 0,\;\forall y\in R.
\end{array}
\ \right.  \label{h2.2}
\end{equation}
Obviously (b) implies (a). This kind of BSDE was intrduced by Bismut \cite
{Bismut73}, \cite{Bismut78} for the case where $g$ is a linear function of $%
(y,z)$. Pardoux and Peng \cite{Pardoux-Peng1990} obtained the following
result (see Theorem \ref{th2.1} for a more general situation): for each $%
X\in L^2(\mathcal{F}_t)$, there exists a unique solution $(Y,Z)\in S_{%
\mathcal{F}}^2(0,t)\times L_{\mathcal{F}}^2(0,t;R^d)$ of the BSDE (\ref
{tsBSDE}).

\begin{definition}
\label{EgstX}We denote by $\mathcal{E}_{s,t}^g[X]:=Y_s$, $0\leq s\leq t$.
\end{definition}

We thus define a system of operators
\[
\mathcal{E}_{s,t}^g[\cdot ]:L^2(\mathcal{F}_t)\rightarrow L^2(\mathcal{F}_s)%
\hbox{, }0\leq s\leq t\leq T.
\]
We will prove that $(\mathcal{E}_{s,t}^g[\cdot ])_{0\leq s\leq t\leq T}$
forms an $\mathcal{F}_t$--consistent evaluation on $[0,T]$. This evaluation
is entirely determined by the simple function $g$.

\section{Main result: $\mathcal{E}_{s,t}[\cdot ]$ is determined by a
function $g$\label{ss3}}

From now on the system $\mathcal{E}_{s,t}[\cdot ]:L^2(\mathcal{F}%
_t)\rightarrow L^2(\mathcal{F}_s)$, $0\leq s\leq t\leq T$, is always a fixed
$\mathcal{F}_t$--consistent nonlinear evaluation, i.e., satisfying
(A1)--(A4), with additional assumptions (A4$_0$) and the following ${%
\mathcal{E}}^{g_\mu}$--domination assumption:

\smallskip\

\noindent \textbf{(A5)} there exists a sufficiently large number $\mu >0$
such that, for each $0\leq s\leq t\leq T$,
\begin{equation}
\mathcal{E}_{s,t}[X]-\mathcal{E}_{s,t}[X^{\prime }]\leq \mathcal{E}%
_{s,t}^{g_\mu }[X-X^{\prime }],\;\;\forall X,X^{\prime }\in L^2(\mathcal{F}%
_t),  \label{e3.1}
\end{equation}
where the function $g_\mu $ is
\begin{equation}
g_\mu (y,z):=\mu |y|+\mu |z|,\;(y,z)\in R\times R^d.  \label{e2.9}
\end{equation}

The \textbf{main theorem} of this paper is:

\begin{theorem}
\label{m7.1}Let $\mathcal{E}_{s,t}[\cdot ]:L^2(\mathcal{F}_t)\rightarrow L^2(%
\mathcal{F}_s)$, $0\leq s\leq t\leq T$, satisfy (A1)--(A4), (A4$_0$) and
(A5). Then there exists a function $g(\omega ,t,y,z)$ satisfying (\ref{h2.1}%
) with $g(s,0,0)\equiv 0$, such that, for each $0\leq s\leq t\leq T$,
\begin{equation}
\mathcal{E}_{s,t}[X]=\mathcal{E}_{s,t}^g[X],\;\;\forall X\in L^2(\mathcal{F}%
_t).  \label{e7.1}
\end{equation}
\end{theorem}

\begin{remark}
The case where $\mathcal{E}_{s,t}[\cdot ]$ satisfy (A1)--(A5), without (A4$_0
$), can be obtained as corollaries of the this main theorem. This will be
given in Corollaries \ref{m7.1a} and \ref{m7.1b}. In this more general
situation the condition $g(s,0,0)\equiv 0$ is not imposed.
\end{remark}

We consider some special situations of our theorem.

\begin{example}
If moreover, $g(s,y,0)\equiv 0$. Then, by \cite{Peng1997}, (A2') holds.
Thus, according to Proposition \ref{m2a2}, $\mathcal{E}_{s,t}^g[\cdot ]$
becomes an $\mathcal{F}_t$--consistent nonlinear expectation:
\[
\mathcal{E}[X|\mathcal{F}_t]=\mathcal{E}_g[X|\mathcal{F}_t]:=\mathcal{E}%
_{s,t}^g[X]=\mathcal{E}_{s,T}^g[X].
\]
This is the so called $g$--expectation introduced in \cite{Peng1997}.
\end{example}

This extends non trivially the result obtained in \cite{CHMP2002}, (see also
\cite{Peng2003b} for a more systematical presentation and explanations in
finance), where we needed a more strict domination condition plus the
following assumption
\[
\mathcal{E}[X+\eta
|\mathcal{F}_t]=\mathcal{E}[X|\mathcal{F}_t]+\eta ,\;\forall \eta
\in \mathcal{F}_t\hbox{.}
\]
Under these assumptions we have proved in \cite{CHMP2002} that there exists
a unique function $g=g(s,z)$, with $g(s,0)\equiv 0$, such that $\mathcal{E}%
_g[X]\equiv \mathcal{E}[X]=\mathcal{E}[X|\mathcal{F}_0]$.

\begin{example}
Consider a financial market consisting of $d+1$ assets: one bond and $d$
stocks. We denote by $P_0(t)$ the price of the bond and by $P_i(t)$ the
price of the $i$-th stock at time $t$. We assume that $P_0$ is the solution
of the ordinary differential equation: $dP_0(t)=r(t)P_0(t)dt,$ and $%
\{P_i\}_{i=1}^d$ is the solution of the following SDE
\begin{eqnarray*}
dP_i(t) &=&P_i(t)[b_i(t)dt+{\sum }_{j=1}^d\sigma _{ij}(t)dB_t^j], \\
P_i(0) &=&p_i,\quad i=1,\cdots ,d.
\end{eqnarray*}
Here $r$ is the interest rate of the bond; $\{b_i\}_{i=1}^d$ is the rate of
the expected return, $\{\sigma _{ij}\}_{i,j=1}^d$ the volatility of the
stocks. We assume that $r$, $b$, $\sigma $ and $\sigma ^{-1}$ are all $%
\mathcal{F}_t$--adapted and uniformly bounded processes on $[0,\infty )$.
Black and Scholes have solved the problem of the market evaluation of an
European type of derivative $X\in L^2(\mathcal{F}_T)$ with maturity $T$. In
the point of view of BSDE, the problem can be treated as follows: consider
an investor who has, at a time $t\leq T$, $n_0(t)$ bonds and $n_i(t)$ $i$%
-stocks, $i=1,\cdots ,d$, i.e., he invests $n_0(t)P_0(t)$ in bond and $\pi
_i(t)=n_i(t)P_i(t)$ in the $i$-th stock. $\pi (t)=(\pi _1(t),\cdots ,\pi
_d(t))$, $0\le t\le T$ is an $R^d$ valued, square-integrable and adapted
process. We define by $y(t)$ the investor's wealth invested in the market at
time $t$:
\[
y(t)=n_0(t)P_0(t)+{\sum }_{i=1}^d\pi _i(t).
\]
We make the so called self--financing assumption: in the period $[0,T]$, the
investor does not withdraw his money from, or put his money in his account $%
y_t$. Under this condition, his wealth $y(t)$ evolves according to
\[
dy(t)=n_0(t)dP_0(t)+{\sum }_{i=1}^dn_i(t)dP_i(t).
\]
or
\[
dy(t)=[r(t)y(t)+{\sum }_{i=1}^d(b_i(t)-r(t))\pi _i(t)]dt+{\sum }%
_{i,j=1}^d\sigma _{ij}(t)\pi _i(t)dB_t^j.
\]
We denote $g(t,y,z):=-r(t)y-{\sum }_{i,j=1}^d(b_i(t)-r(t))\sigma
_{ij}^{-1}(t)z_j$. Then, by the variable change $z_j(t)={\sum }%
_{i=1}^d\sigma _{ij}(t)\pi _i(t)$, the above equation is
\[
-dy(t)=g(t,y(t),z(t))dt-z(t)dB_t.
\]
We observe that the function $g$ satisfies (\ref{h2.1}). It follows from the
existence and uniqueness theorem of BSDE (Theorem \ref{th2.1}) that for each
derivative $X\in L^2(\mathcal{F}_T)$, there exists a unique solution $%
(y(\cdot ),z(\cdot ))\in L_{\mathcal{F}}^2(0,T;R^{1+d})$ with the terminal
condition $y_T=X$. This meaning is significant: in order to replicate the
derivative $X$, the investor needs and only needs to invest $y(t)$ at the
present time $t$ and then, during the time interval $[t,T]$ and then to
perform the portfolio strategy $\pi _i(s)=\sigma _{ij}^{-1}(s)z_j(s)$.
Furthermore, by Comparison Theorem of BSDE, if he wants to replicate a $%
X^{\prime }$ which is bigger than $X$, (i.e., $X^{\prime }\geq X$, a.s., $%
P(X^{\prime }\geq X)>0$), then he must pay more, i.e., there this no
arbitrage opportunity. This $y(t)$ is called the Black--Scholes price, or
Black--Scholes evaluation, of $X$ at the time $t$. We define, as in (\ref
{e2.4}), $\mathcal{E}_{t,T}^g[X]=y_t$. We observe that the function $g$
satisfies (b) of condition (\ref{h2.2}). It follows from Proposition \ref
{p2.1} that $\mathcal{E}_{t,T}^g[\cdot ]$ satisfies properties (A1)--(A4)
for $\mathcal{F}_t$--consistent evaluation.
\end{example}

\begin{example}
An very important problem is: if we know that the evaluation of an
investigated agent is a $g$--evaluation $\mathcal{E}^g$, how to find this
function $g$. We now consider a case where $g$ depends only on $z$, i.e., $%
g=g(z):\mathbf{R}^d\rightarrow \mathbf{R}$. In this case we can find such $g$
by the following testing method. Let $\bar z\in \mathbf{R}^d$ be given. We
denote $Y_s:=\mathcal{E}_{s,T}^g[\bar z(B_T-B_t)]$, $s\in [t,T]$, where $t$
is the present time. It is the solution of the following BSDE
\[
Y_s=\bar z(B_T-B_t)+\int_s^Tg(Z_u)du-\int_s^TZ_udB_u,\;s\in [t,T].
\]
It is seen that the solution is $Y_s=\bar z(B_s-B_t)+\int_s^Tg(\bar z)ds$, $%
Z_s\equiv \bar z$. Thus
\[
\mathcal{E}_{t,T}^g[\bar z(B_T-B_t)]=Y_t=g(\bar z)(T-t),
\]
or
\begin{equation}
g(\bar z)=(T-t)^{-1}\mathcal{E}_{t,T}^g[\bar z(B_T-B_t)].  \label{e2.exm1}
\end{equation}
Thus the function $g$ can be tested as follows: at the present time $t$, we
ask the investigated agent to evaluate $\bar z(B_T-B_t)$. We thus get $%
\mathcal{E}_{t,T}^g[\bar z(B_T-B_t)]$. Then $g(\bar z)$ is obtained by (\ref
{e2.exm1}).
\end{example}

\begin{remark}
The above test works also for the case $g:[0,\infty )\times \mathbf{R}%
^d\rightarrow \mathbf{R}$, or for a more general situation $g=\gamma
y+g_0(t,z)$.
\end{remark}

An interesting problem is, in general, how to find the function $g$ through
a testing of the input--output behaviour of $\mathcal{E}^g[\cdot ]$? Let $%
b:R^n\longmapsto R^n$, $\bar \sigma :R^n\longmapsto R^{n\times d}$ be two
Lipschitz functions.
\[
X_s^{t,x}=x+\int_t^sb(X_r^{t,x})dr+\int_t^s\sigma (X_r^{t,x})dB_r,\;\;s\geq
t.
\]
The following result was obtained in Proposition 2.3 of \cite{BCHMP00}.

\begin{proposition}
We assume that the generator $g$ satisfies (\ref{h2.1})\textbf{. }We also
assume that, for each fixed $(y,z)$, $g(\cdot ,y,z)\in D_{\mathcal{F}}^2(0,T)
$. Then for each $(t,x,p,y)\in [0,\infty )\times R^n\times R^n\times R$, we
have
\[
L^2\hbox{--}\lim_{\epsilon \rightarrow 0}\frac 1\epsilon [\mathcal{E}%
_{t,t+\epsilon }^g[y+p\cdot (X_{t+\epsilon }^{t,x}-x)]-y]=g(t,y,\sigma
^T(x)p)+p\cdot b(x).\label{limE}
\]
\end{proposition}

\section{A more general formulation: $\mathcal{E}_{s,t}^g[\cdot ;K]$%
--evaluation\label{ass3}}

To prove our main result, we need to introduce a more general type of $%
\mathcal{F}_t$--consistent evaluations $\mathcal{E}_{s,t}[\cdot ;K]$ induced
by $\mathcal{E}_{s,t}[\cdot ]$, for each given process $K\in D_{\mathcal{F}%
}^2(0,T)$. In finance, $K$ often represents a dividend and/or a consumption
process. We will firstly consider $\mathcal{E}_{s,t}^g[\cdot ;K]$. For
technical convenience, we will directly consider stopping times $%
\sigma\leq\tau\leq T$ in the place of deterministic times $s$ and $t$.

Let $\tau \in \mathcal{S}_T$ be a given stopping time. We consider the
following backward stochastic differential equation:
\begin{equation}
Y_s=X+K_\tau -K_{s\wedge \tau }+\int_{s\wedge \tau }^\tau
g(r,Y_r,Z_r)dr-\int_{s\wedge \tau }^\tau Z_rdB_r,\;s\in [0,T].  \label{tBSDE}
\end{equation}
Here the pair $(Y,Z)$ is the unknown process to be solved. $X\in L^2(%
\mathcal{F}_\tau )$, $K\in D_{\mathcal{F}}^2(0,T)$ are given.

We recall the following basic results of BSDE.

\begin{theorem}
\label{th2.1}(\cite{Pardoux-Peng1990}, \cite{Peng1997a}) We assume (\ref
{h2.1}). Then there exists a unique solution $(Y,Z)\in L_{\mathcal{F}%
}^2(0,\tau ;R\times R^d)$ of BSDE (\ref{tBSDE}). We denote it by
\begin{equation}
(Y_s^{\tau ,X,K},Z_s^{\tau ,X,K})=(Y_s,Z_s),\;s\in [0,\tau ].  \label{e2.1a}
\end{equation}
We have
\[
Y^{\tau ,X,K}+K\in S_{\mathcal{F}}^2(0,\tau ).
\]
and the estimate
\begin{eqnarray}
\ \ \ \ \ \ \ \ E\int_0^\tau |Z_s^{\tau ,X,K}|^2ds+E[\sup_{s\in [0,\tau
]}|Y_s^{\tau ,X,K}-K_s|^2]  \label{e2.1} \\
\ \leq CE[(X+K_\tau )^2]+CE\int_0^\tau (K_s^2+|g(s,0,0)|^2)ds,  \nonumber
\end{eqnarray}
where the constant $C$ depends only on $\mu $ and $T$. Furthermore, let $%
X^{\prime }\in L^2(\mathcal{F}_\tau )$, $K^{\prime }\in D_{\mathcal{F}%
}^2(0,T)$ be also given. Then we have
\begin{eqnarray}
\ \ \ \ \ \ \ \ \ \ E[\sup_{s\in [0,\tau ]}|Y_s^{\tau ,X,K}-Y_s^{\tau
,X^{\prime },K^{\prime }}+K_s-K_s^{\prime }|^2]+E\int_0^\tau |Z_s^{\tau
,X,K}-Z_s^{\tau ,X^{\prime },K^{\prime }}|^2ds  \label{e2.2} \\
\ \leq CE[(X-X^{\prime }+K_\tau -K_\tau ^{\prime })^2]+CE\int_0^\tau
(K_s-K_s^{\prime })^2ds,  \nonumber
\end{eqnarray}
where the constant $C$ depends only on $\mu $ and $T$.
\end{theorem}

\smallskip\noindent\textbf{Proof. } In \cite{Pardoux-Peng1990} (see also \cite{EPQ1997}%
), the result of BSDE is for $\tau =T$ and $K_t=\int_0^t\phi _sds$ for some $%
\phi \in L_{\mathcal{F}}^2(0,T)$. The present situation can be treated by
defining (see \cite{Peng1997})
\begin{equation}
\begin{array}{rl}
\bar Y_s & :=Y_s+K_s, \\
\bar g(s,y,z) & :=g(s,y-K_s,z)1_{[0,\tau ]}(s)
\end{array}
\label{e-gbar}
\end{equation}
and considering the following equivalent BSDE
\begin{equation}
\bar Y_s=X+K_\tau +\int_s^T\bar g(r,\bar Y_r,Z_r)dr-\int_s^TZ_rdB_r,\;s\in
[0,T].  \label{TBSDE}
\end{equation}
It is clear that $\bar Y_s\equiv X+K_s$, $Z_s\equiv 0$ on $[\tau ,T]$. Since
$\bar g$ is a Lipschitz function with the same Lipschitz constant $\mu $ and
\[
\bar g(\cdot ,0,0)=g(\cdot ,-K_{\cdot },0)1_{[0,\tau ]}(\cdot )\in L_{%
\mathcal{F}}^2(0,T),
\]
thus, by \cite{Pardoux-Peng1990}, \cite{Peng1997a}, the BSDE (\ref{TBSDE})
has a unique solution $(\bar Y,Z)$. We also have
\[
E\int_0^\tau |Z|^2ds+E[\sup_{s\in [0\tau ]}|\bar Y_s|^2]\leq CE[(X+K_\tau
)^2]+CE\int_0^\tau (K_s^2+|g(s,0,0)|^2)ds,
\]
where the constant $C$ depends only on $\mu $ and $T$. We thus have estimate
(\ref{e2.1}). Moreover, let $(\bar Y^{\prime },Z^{\prime })$ denotes the
solution of the (\ref{TBSDE}) with $X^{\prime }$ and $K^{\prime }$ in the
place of $X$ and $K$. We have the following classical estimate:
\begin{eqnarray*}
&&\ \ \ \ \ E[\sup_{s\in [0,\tau ]}|\bar Y_s-\bar Y_s^{\prime
}|^2]+E\int_0^\tau |Z_s-Z_s^{\prime }|^2ds \\
\ &\leq &CE[(X-X^{\prime }+K_\tau -K_\tau ^{\prime })^2]+CE\int_0^\tau
(K_s-K_s^{\prime })^2ds.
\end{eqnarray*}
where $C$ is the same as in (\ref{e2.1}). We then have
(\ref{e2.2}). The proof is complete. $\Box $\medskip\

We introduce a new notation.
\begin{definition}
\label{d2.2}We denote, for $\sigma ,\tau \in S_T$, $\sigma \leq
\tau $,
\begin{eqnarray}
\mathcal{E}_{\sigma ,\tau }^g[X;K_{\cdot }]:=Y_\sigma ^{\tau ,X,K}  \label{e2.3} \\
\mathcal{E}_{\sigma ,\tau }^g[X]:=\mathcal{E}_{\sigma ,\tau }^g[X;0].
\label{e2.4}
\end{eqnarray}
This notion generalizes that of $\mathcal{E}_{s,t}^g[\cdot ]$ in Definition
\ref{EgstX}. Clearly when (\ref{h2.2})--(a) is satisfied, we have $\mathcal{E%
}_{\sigma ,\tau }^g[0]=\mathcal{E}_{\sigma ,\tau }^g[0;0]=0$. In particular $%
\mathcal{E}_{s,t}^g[0]\equiv 0$, $0\leq s\leq t\leq T$.
\end{definition}

\begin{remark}
About the notations $\mathcal{E}^g[\cdot ]$. This notation was firstly
introduced in \cite{Peng1997} in the case where $g$ satisfies (\ref{h2.2}%
)--(b). In this situation it is easy to check that
\[
\mathcal{E}_{s,t}^g[X]\equiv \mathcal{E}_{s,T}^g[X],\;\forall 0\leq s\leq
t\leq T.
\]
In other words, $\mathcal{E}^g$--is a nonlinear expectation, called $g$%
--expectation. The general situation, i.e., without (\ref{h2.2}) was
introduced in \cite{Peng1997a} and \cite{Chen-Peng2001}.
\end{remark}

By the above existence and uniqueness theorem, we have for each stopping
times $0\leq \rho \leq \sigma \leq \tau \leq T$ and for each $X\in L^2(%
\mathcal{F}_\tau )$ and $K\in D_{\mathcal{F}}^2(0,T)$,
\begin{equation}
\mathcal{E}_{\rho ,\sigma }^g[\mathcal{E}_{\sigma ,\tau
}^g[X;K_{\cdot }];K_{\cdot }]=\mathcal{E}_{\rho ,\tau
}^g[X;K_{\cdot }],\;\hbox{a.s.} \label{e2.5}
\end{equation}
It is also easy to check that, with the notation $g_{-}(t,y,z):=-g(t,-y,-z)$
\begin{equation}
-\mathcal{E}_{\sigma ,\tau }^g[X;K_{\cdot }]=\mathcal{E}_{\sigma ,\tau
}^{g_{-}}[-X;-K_{\cdot }].  \label{e2.5.1}
\end{equation}

We will see that $\{\mathcal{E}_{t,T}^g\left[ X\right]\}_{0\leq t\leq T}$, $%
X\in L^2(\mathcal{F}_T)$ form an $\mathcal{F}_t$--consistent nonlinear
evaluation. The following monotonicity property is the comparison theorem of
BSDE.

\begin{theorem}
\label{th2.2}We assume (\ref{h2.1}). If we assume that the elements in the
above theorem satisfy $X\geq X^{\prime }$, a.s., and that $K-K^{\prime }$ is
an increasing process. Then, for each stopping times $0\leq \sigma \leq \tau
\leq T$, we have
\begin{equation}
\mathcal{E}_{\sigma ,\tau }^g[X;K_{\cdot }]\geq
\mathcal{E}_{\sigma ,\tau }^g[X^{\prime };K_{\cdot }^{\prime
}]\hbox{, a.s.}  \label{e2.6}
\end{equation}
In particular,
\begin{equation}
\mathcal{E}_{\sigma ,\tau }^g[X]\geq \mathcal{E}_{\sigma ,\tau
}^g[X^{\prime }]\hbox{, a.s.}  \label{e2.7}
\end{equation}
If $A\in D_{\mathcal{F}}^2(0,T)$ is an increasing process, then
\begin{equation}
\mathcal{E}_{\sigma ,\tau }^g[X;A_{\cdot }]\geq
\mathcal{E}_{\sigma ,\tau }^g[X]\hbox{.}  \label{e2.8}
\end{equation}
\end{theorem}

\smallskip\noindent\textbf{Proof. }The case $K_t\equiv K_t^{\prime }\equiv 0$ is the
classical comparison theorem of BSDE. The present general situation, see
\cite{Peng1997a} or \cite{Peng2003b}.

We recall the special function $g_\mu (y,z)$ defined in (\ref{e2.9})

\begin{corollary}
\label{c2.1}$\mathcal{E}^g$ is dominated by $\mathcal{E}^{g_\mu }$ in the
following sense: for each stopping times $0\leq \sigma \leq \tau \leq T$ and
$X$, $X^{\prime }\in L^2(\mathcal{F}_\tau )$, $K$, $K^{\prime }\in D_{%
\mathcal{F}}^2(0,T)$, we have
\begin{equation}
\mathcal{E}_{\sigma ,\tau }^g[X;K_{\cdot }]-\mathcal{E}_{\sigma
,\tau }^g[X^{\prime };K_{\cdot }^{\prime }]\leq
\mathcal{E}_{\sigma ,\tau }^{g_\mu }[X-X^{\prime };(K-K^{\prime
})_{\cdot }]\hbox{, a.s.}  \label{e2.10}
\end{equation}
where $\mu $ is the Lipschitz constant of $g$ given in (\ref{h2.1}).
\end{corollary}

\smallskip\noindent\textbf{Proof. }By the definition of $\mathcal{E}^g[\cdot ]$, The
processes defined by $Y_s=\mathcal{E}_{s\wedge \tau ,\tau }^g[X;K_{\cdot }]$
and $Y_s^{\prime }=\mathcal{E}_{s\wedge \tau ,\tau }^g[X^{\prime };K_{\cdot
}^{\prime }]$ solve respectively the following BSDEs on $[0,\tau ]$:
\begin{eqnarray*}
Y_s &=&X+K_\tau -K_{s\wedge \tau }+\int_{s\wedge \tau }^\tau
g(r,Y_r,Z_r)dr-\int_{s\wedge \tau }^\tau Z_rdB_r, \\
Y_s^{\prime } &=&X^{\prime }+K_\tau ^{\prime }-K_{s\wedge \tau }^{\prime
}+\int_{s\wedge \tau }^\tau g(r,Y_r^{\prime },Z_r^{\prime })dr-\int_{s\wedge
\tau }^\tau Z_r^{\prime }dB_r.
\end{eqnarray*}
We denote $\hat Y=Y-Y^{\prime }$, $\hat Z=Z-Z^{\prime }$ and
\[
\hat K_t=K_t-K_t^{\prime }+\int_0^t[-g_\mu (\hat Y_s,\hat
Z_s)+g(s,Y_s,Z_s)-g(s,Y_s,Z_s)]ds.
\]
$(\hat Y,\hat Z)$ solves BSDE
\[
\hat Y_s=X-X^{\prime }+\hat K_\tau -\hat K_{s\wedge \tau }+\int_{s\wedge
\tau }^\tau g_\mu (r,\hat Y_r,\hat Z_r)dr-\int_{s\wedge \tau }^\tau \hat
Z_rdB_r.
\]
We compare it to the BSDE
\[
\bar Y_s=X-X^{\prime }+(K-K^{\prime })_\tau -(K-K^{\prime })_{s\wedge \tau
}+\int_{s\wedge \tau }^\tau g_\mu (r,\bar Y_r,\bar Z_r)dr-\int_{s\wedge \tau
}^\tau \bar Z_rdB_r.
\]
Since $d(K-K^{\prime }-\hat K)_s\geq 0$, thus, by comparison theorem, i.e.,
Theorem \ref{th2.2}, $\bar Y_s\geq \hat Y_s=Y_s-Y_s^{\prime }$. We thus have
(\ref{e2.10}). $\Box $\medskip\

\begin{proposition}
\label{p2.0}We have the following uniform estimate: for each $X\in L^2(%
\mathcal{F}_T)$ and $g_0(\cdot )\in L_{\mathcal{F}}^2(0,T)$
\begin{equation}
E[(\mathcal{E}_{t,T}^{g_\mu }[X;\int_0^{\cdot }g_0(s)ds])^2]\leq
E[|X|^2]e^{\beta (T-t)}+E[\int_t^Te^{\beta (s-t)}|g_0(s)|^2ds].
\label{e2.gm}
\end{equation}
and
\begin{equation}
E[\sup_{t\in [0,T]}(\mathcal{E}_{t,T}^{g_\mu }[X;\int_0^{\cdot
}g_0(s)ds])^2]\leq CE[|X|^2+\int_0^T|g_0(s)|^2ds].  \label{e2.gm1}
\end{equation}
where $\beta =2\mu ^2+2\mu +2$. The constant $C$ depends only on $T$ and the
Lipschitz constant $\mu $ in (\ref{h2.1}).
\end{proposition}

\smallskip\noindent\textbf{Proof. }$Y_t=\mathcal{E}_{t,T}^{g_\mu }[X;\int_0^{\cdot
}g_0(s)ds]$ satisfies the following BSDE on $[0,T]$:
\begin{equation}
Y_t=X+\int_t^Tg_0(s)ds+\mu \int_t^T(|Y_s|+|Z_s|)ds-\int_t^TZ_sdB_s.
\label{e2.gm2}
\end{equation}
We apply It\^o's formula to $|Y_t|^2e^{\beta t}$:
\begin{eqnarray*}
&&E[|Y_t|^2e^{\beta t}+\int_t^Te^{\beta s}(|Z_s|^2+\beta |Y_s|^2)ds \\
&=&E[|X|^2]e^{\beta T}+E\int_t^Te^{\beta s}2Y_s(g_0(s)+\mu |Y_s|+\mu |Z_s|)ds
\\
&\leq &E[|X|^2]e^{\beta T}+E\int_t^Te^{\beta s}[2\mu ^2+2\mu
+1)|Y_s|^2+|g_0(s)|^2+\frac 12|Z_s|^2]ds.
\end{eqnarray*}
We thus have (\ref{e2.gm}) and
\[
E\int_0^Te^{\beta s}(|Z_s|^2+|Y_s|^2)ds\leq E[|X|^2]e^{\beta
T}+E\int_0^Te^{\beta s}|g_0(s)|^2ds.
\]
With (\ref{e2.gm2}), we now apply BDG's inequality to $Y_t^2$. Then (\ref
{e2.gm1}) follows.

We now can assert that

\begin{proposition}
\label{p2.1}Let $g$ satisfies (\ref{h2.1}). Then, for each fixed $K\in D_{%
\mathcal{F}}^2(0,T)$, the system of operators
\begin{equation}
\mathcal{E}_{s,t}^g[X;K_{\cdot }]:L^2(\mathcal{F}_t)\rightarrow L^2(\mathcal{%
F}_s),\;0\leq s\leq t\leq T.  \label{e2.11}
\end{equation}
defined in (\ref{e2.3}) is an $\mathcal{F}_t$--consistent nonlinear
evaluation, i.e., it satisfies (A1)--(A4) of Definition \ref{d2.1}.
\end{proposition}

\smallskip\noindent\textbf{Proof. } (A1) is given by (\ref{e2.6}). (A2) is clearly
true by the definition. (A3) is proved by (\ref{e2.5}). We now consider
(A4). In fact we can prove stronger results: for each stopping times $0\leq
\sigma \leq \tau \leq T$ and $X\in L^2(\mathcal{F}_\tau )$, we have
\begin{equation}
1_A\mathcal{E}_{\sigma ,\tau }^g[X;K_{\cdot }]=1_A\mathcal{E}_{\sigma ,\tau
}^g[1_AX;K_{\cdot }],\;\;\forall A\in \mathcal{F}_\tau .  \label{e2.12}
\end{equation}
as well as
\begin{equation}
1_A\mathcal{E}_{\sigma ,\tau }^g[X;K_{\cdot }]=\mathcal{E}_{\sigma ,\tau
}^{g_{\sigma ,A}}[1_AX;K_{\cdot }^{\sigma ,A}],\;\;\forall A\in \mathcal{F}%
_\tau ,  \label{e2.12a}
\end{equation}
where we set
\begin{eqnarray}
g_{\sigma ,A}(t,y,z) &:&=1_{[0,\sigma )}(t)g(t,y,z)+1_{[\sigma
,T]}(t)1_Ag(t,y,z),  \label{e2.12b} \\
K_t^{\sigma ,A} &:&=1_{[0,\sigma )}(t)K_t+1_{[\sigma ,T]}(t)1_A(K_t-K_\sigma
).  \label{e2.12c}
\end{eqnarray}
We will give the proof of (\ref{e2.12}). The proof of (\ref{e2.12a}) is
similar. According to BSDE (\ref{tBSDE}) for each stopping time $\rho \in
[\sigma ,\tau ]$, $Y_\rho :=\mathcal{E}_{\rho ,\tau }^g[X;K_{\cdot }]$ and $%
\bar Y_\rho :=\mathcal{E}_{\rho ,\tau }^g[1_AX;K_{\cdot }]$ solve
respectively
\[
Y_\rho =X+K_\tau -K_\rho +\int_\rho ^\tau g(r,Y_r,Z_r)dr-\int_\rho ^\tau
Z_rdB_r,
\]
and
\[
\bar Y_\rho =1_AX+K_\tau -K_\rho +\int_\rho ^\tau g(r,\bar Y_r,\bar
Z_r)dr-\int_\rho ^\tau \bar Z_rdB_r
\]
We multiply $1_A$, $A\in \mathcal{F}_\sigma $ on both sides of the above two
BSDEs. Since $1_Ag(r,Y_r,Z_r)=1_Ag(r,Y_r1_A,Z_r1_A)$, we have
\[
1_AY_\rho =1_AX+1_AK_\tau -1_AK_\rho +\int_\rho ^\tau
1_Ag(r,1_AY_r,1_AZ_r)dr-\int_\rho ^\tau 1_AZ_rdB_r,
\]
and
\[
1_A\bar Y_\rho =1_AX+1_AK_\tau -1_AK_\rho +\int_\rho ^\tau 1_Ag(r,1_A\bar
Y_r,1_A\bar Z_r)dr-\int_\rho ^\tau 1_A\bar Z_rdB_r.
\]
It is clear that $1_AY_\rho $ and $1_A\bar Y_\rho $ satisfy exactly the same
BSDE with the same terminal condition on $[\sigma ,\tau ]$. By uniqueness of
BSDE, $1_AY_\rho \equiv 1_A\bar Y_\rho $ on $[\sigma ,\tau ]$, i.e., $1_A%
\mathcal{E}_{\sigma ,\tau }^g[X;K_{\cdot }]\equiv 1_A\mathcal{E}_{\sigma
,\tau }^g[1_AX;K_{\cdot }]$. The proof is complete. $\Box $\medskip\

We now consider nonlinear martingales induced by $\mathcal{E}^g$.

\begin{definition}
\label{d2.3}Let $K\in D_{\mathcal{F}}^2(0,T)$ be given. A process $Y\in D_{%
\mathcal{F}}^2(0,T)$ is said to be an $\mathcal{E}^g[\cdot ;K]$--martingale
(resp. $\mathcal{E}^g[\cdot ;K]$--supermartingale, $\mathcal{E}^g[\cdot ;K]$%
--submartingale) if for each $0\leq s\leq t\leq T$%
\begin{equation}
\mathcal{E}_{s,t}^g[Y_t;K_{\cdot }]=Y_s\hbox{, (resp. }\leq
Y_s\hbox{, }\geq Y_s\hbox{)}.  \label{e2.13}
\end{equation}
\end{definition}

\begin{remark}
\label{r2.1}If $(y,z)\in L_{\mathcal{F}}^2(0,T;R\times R^d)$ solves the BSDE
\[
y_s=y_t+K_t-K_s+\int_s^tg(r,y_r,z_r)dr-\int_s^tz_rdB_r,\hbox{
}s\leq t.
\]
It is clear that $(-y,-z)$ solves
\begin{eqnarray*}
-y_s=-y_t+(-K_t)-(-K_s) \\
\ \ \ \ \ \ \ \ +\int_s^t[-g(r,-(-y_r),-(-z_r))dr-\int_s^t(-z_r)dB_r.
\end{eqnarray*}
Thus, if $y$ is an $\mathcal{E}^g[\cdot ;K_{\cdot }]$--martingale (resp. $%
\mathcal{E}^g[\cdot ;K]$--supermartingale, $\mathcal{E}^g[\cdot ;K]$%
--submartingale), then $-y$ is an $\mathcal{E}_{s,t}^{g_{*}}[\cdot
;-K_{\cdot }]$--martingale (resp. $\mathcal{E}^{g_{*}}[\cdot ;K]$%
--submartingale, $\mathcal{E}^{g_{*}}[\cdot ;K]$--supermartingale), where we
denote
\[
g_{*}(t,y,z):=-g(t,-y,-z).
\]
Therefor many results concerning $\mathcal{E}^g[\cdot ;K]$--supermartingales
can be also applied to situations of submartingales.
\end{remark}

\begin{example}
Let $X\in L^2(\mathcal{F}_T)$ and $A\in D_{\mathcal{F}}^2(0,T)$ be given
such that $A$ is an increasing process. By the monotonicity of $\mathcal{E}^g
$, i.e., Theorem \ref{th2.2}, we have, for $t\in [0,T]$,
\[
\begin{array}{l}
Y_t:=\mathcal{E}_{t,T}^g[X]=\mathcal{E}_{t,T}^g[X;0]\hbox{\ is a }\mathcal{E}%
^g\hbox{--martingale,} \\
Y_t^{+}:=\mathcal{E}_{t,T}^g[X;A]\hbox{\ is a }\mathcal{E}^g\hbox{%
--supermartingale,} \\
Y_t^{-}:=\mathcal{E}_{t,T}^g[X;-A]\hbox{ is a }\mathcal{E}^g\hbox{%
--submartingale.}
\end{array}
\]
\end{example}

As in classical situation, an interesting and hard problem is the inverse
one: if $Y$ is an $\mathcal{E}^g$--supermartingale, can we find an
increasing and predictable process $A$ such that $Y_t\equiv \mathcal{E}%
_{t,T}^g[X;A]$? This nonlinear version of Doob--Meyer's decomposition
theorem will be stated as follows. It plays a crucially important role in
this paper.

The following result is a nonlinear version of optional sampling theorem for
$g$--supermartingale. See \cite{Peng2003b} (also Theorem \ref{m8.14} for a
more general situation).

\begin{proposition}
\label{p2.2} Let $g$ satisfy (\ref{h2.1})--(i) and (ii) and let $Y\in D_{%
\mathcal{F}}^2(0,T)$ be an $\mathcal{E}^g$--martingale (resp. $\mathcal{E}^g$%
--supermartingale, $\mathcal{E}^g$--submartingale). Then for each stopping
times $0\leq \sigma \leq \tau \leq T$, we have
\begin{equation}
\mathcal{E}_{\sigma ,\tau }^g[Y_\tau ]=Y_\sigma \hbox{, (resp.
}\leq Y_\sigma \hbox{, }\geq Y_\sigma \hbox{)}.  \label{e2.14}
\end{equation}
\end{proposition}

We have the following $\mathcal{E}^g$--supermartingale decomposition theorem
of Doob--Meyer's type. This nonlinear decomposition theorem was obtained in
\cite{Peng1999}. But the formulation using the new notation $\mathcal{E}%
_{t,T}^g[\cdot ;A]$ is new. In fact we think this is the intrinsic
formulation since it becomes necessary in the more abstract situation of the
$\mathcal{E}$--supermartingale decomposition theorem, i.e., Theorem \ref
{m6.1} which can considered as a generalization of the following result.

\begin{proposition}
\label{p2.3} We assume (\ref{h2.1})--(i) and (ii). Let $Y\in D_{\mathcal{F}%
}^2(0,T)$ be an $\mathcal{E}^g$--supermartingale. Then there exists a unique
increasing process $A\in D_{\mathcal{F}}^2(0,T)$ (thus predictable) with $%
A_0=0$, such that
\begin{equation}
Y_t=\mathcal{E}_{t,T}^g[Y_T;A],\;\forall 0\leq t\leq T.  \label{e2.15}
\end{equation}
\end{proposition}

\begin{corollary}
\label{c2.2}Let $K\in D_{\mathcal{F}}^2(0,T)$ be given and let $Y\in D_{%
\mathcal{F}}^2(0,T)$ be an $\mathcal{E}^g[\cdot ;K]$--supermartingale in the
following sense
\begin{equation}
\mathcal{E}_{s,t}^g[Y_t;K]\leq Y_s,\;\forall 0\leq s\leq t\leq T.
\label{e2.16}
\end{equation}
Then there exists a unique increasing process $A\in D_{\mathcal{F}}^2(0,T)$
with $A_0=0$, such that
\begin{equation}
Y_t=\mathcal{E}_{t,T}^g[Y_T;K+A],\;\forall 0\leq t\leq T.  \label{e2.17}
\end{equation}
\end{corollary}

\smallskip\noindent\textbf{Proof. }By the notations of (\ref{e-gbar}) with $\tau =T$,
we have
\begin{equation}
\mathcal{E}_{s,t}^g[Y_t;K]+K_s=\mathcal{E}_{s,t}^{\bar g}[Y_t+K_t].
\label{e2.18}
\end{equation}
It follows that (\ref{e2.16}) is equivalent to
\begin{equation}
\mathcal{E}_{s,t}^{\bar g}[Y_t+K_t]\leq Y_s+K_s,\;\forall 0\leq s\leq t\leq
T.  \label{e2.19}
\end{equation}
In other words, $Y+K$ is an $\mathcal{E}^{\bar g}$--supermartingale in the
sense of (\ref{e2.13}). By the above supermartingale decomposition theorem,
Proposition \ref{p2.3}, there exists an increasing process $A\in D_{\mathcal{%
F}}^2(0,T)$ with $A_0=0$, such that
\begin{equation}
Y_t+K_t=\mathcal{E}_{t,T}^{\bar g}[Y_T+K_T;A],\;\forall 0\leq t\leq T,
\label{e2.20}
\end{equation}
or, equivalently (\ref{e2.17}). $\Box $\medskip\

\section{$\mathcal{E}_{s,t}[\cdot ;K]$ and related properties\label{ss4}}

\subsection{$\mathcal{E}_{s,t}[\cdot ;K]$ and it's main properties}

In this section, $\left\{ \mathcal{E}_{s,t}[\cdot ]\right\} _{0\leq s\leq
t\leq T}$ is a fixed $\mathcal{F}_t$--consistent evaluation satisfying
(A1)--(A5) (without (A4$_0$)) as well as (\ref{e3.1a}). To prove our main
theorem, we shall use $\mathcal{E}_{s,t}[\cdot ]$ to generate an operator $%
\mathcal{E}_{s,t}[\cdot ;K]$, $K\in D_{\mathcal{F}}^2(0,T)$ which plays the
same role as $\mathcal{E}_{s,t}^g[\cdot ;K]$ defined in (\ref{e2.3}).

To this end we first define such operator on the space of step processes
defined by

\begin{equation}
D_{\mathcal{F}}^{2,0}(0,T):=\{K\in D_{\mathcal{F}}^2(0,T),K_t=\sum_{i=0}^N%
\xi _i1_{[t_i,t_{i+1})}(t)\hbox{, for some }0=t_0<t_1\ldots
<t_N=T\}. \label{e3.2}
\end{equation}

We are given $K\in D_{\mathcal{F}}^{2,0}(0,T)$ in the form $%
K_t=\sum_{i=0}^N\xi _i1_{[t_i,t_{i+1})}(t),\;t\in [0,T]$. For each $%
i=0,1,2,\cdots ,N-1$, we define, for $T_i\leq s\leq t\leq T_{i+1}$ and $X\in
L^2(\mathcal{F}_t)$
\[
\mathcal{E}_{s,t}^i[X;K]:=\mathcal{E}_{s,t}[X+K_t-K_s].
\]
We have

\begin{lemma}
\label{m3.2}For each $i=0,1,2,\cdots ,N-1$, $\mathcal{E}_{s,t}[\cdot ;K]$, $%
t_i\leq s\leq t\leq t_{i+1}$ is an $\mathcal{F}_t$--consistent evaluation.
\end{lemma}

\smallskip\noindent\textbf{Proof. }It is easy to check that (A1), (A2) and (A3) holds.
We now prove (A4), i.e., for each $t_i\leq s\leq t\leq t_{i+1}$ and $X\in
L^2(\mathcal{F}_t)$,
\begin{equation}
1_A\mathcal{E}_{s,t}^i[X;K]=1_A\mathcal{E}_{s,t}^i[1_AX;K],\;\forall A\in
\mathcal{F}_s.  \label{e3.4}
\end{equation}
We have
\begin{eqnarray*}
1_A\mathcal{E}_{s,t}^i[X;K_{\cdot }] &=&1_A\mathcal{E}_{s,t}[X+K_t-K_s] \\
\ &=&1_A\mathcal{E}_{s,t}[1_A(1_AX+K_t-K_s)] \\
\ &=&1_A\mathcal{E}_{s,t}[1_AX+K_t-K_s] \\
\ &=&1_A\mathcal{E}_{s,t}^i[1_AX;K_{\cdot }].
\end{eqnarray*}
Thus (A4) holds. $\Box $\medskip\

By Proposition \ref{m2.3}, there exists a unique $\mathcal{F}_t$--consistent
evaluation that coincides with $\mathcal{E}^i[\cdot ;K]$ for each interval $%
[T_i,T_{i+1}]$.

\begin{definition}
\label{m3.1}We denote this unique $\mathcal{F}_t$--consistent evaluation
that coincides with $\mathcal{E}^i[\cdot ;K]$ by $\mathcal{E}_{s,t}[\cdot ;K]
$:
\[
\mathcal{E}_{s,t}[X;K]:X\in L^2(\mathcal{F}_t)\rightarrow L^2(\mathcal{F}_s).
\]
\end{definition}

\begin{lemma}
\label{m3.3}If there is some function $g$ satisfying (\ref{h2.1}) such that $%
\mathcal{E}$ coincides with $\mathcal{E}^g$, i.e., for each $0\leq s\leq
t\leq T$ and $X\in L^2(\mathcal{F}_T)$ we have $\mathcal{E}_{s,t}[X]=%
\mathcal{E}_{s,t}^g[X]$, then, for each $K\in D_{\mathcal{F}}^{2,0}(0,T)$, $%
\mathcal{E}_{s,t}[\cdot ;K]$ also coincides with $\mathcal{E}_{s,t}^g[\cdot
;K]$.
\end{lemma}

\smallskip\noindent\textbf{Proof. } It is easy to check that $\mathcal{E}_{s,t}[X;K]=%
\mathcal{E}_{s,t}^i[X;K]=\mathcal{E}_{s,t}^g[X;K]$, $t_i\leq s\leq t\leq
t_{i+1}$. Thus we can apply Proposition \ref{m2.3} to prove this lemma for $%
0\leq s\leq t\leq T$. $\Box $\medskip\

\begin{lemma}
\label{m3.5}$\mathcal{E}$ is dominated by $\mathcal{E}^{g_\mu }$ in the
following sense: for each $K$, $K^{\prime }\in D_{\mathcal{F}}^{2,0}(0,T)$
and for each $0\leq s\leq t\leq T$, $X$, $X^{\prime }\in L^2(\mathcal{F}_t)$%
, we have
\begin{eqnarray}
\mathcal{E}_{s,t}^{-g_\mu }[X-X^{\prime };(K-K^{\prime })_{\cdot }]\leq
\mathcal{E}_{s,t}[X;K_{\cdot }]-\mathcal{E}_{s,t}[X^{\prime };K_{\cdot
}^{\prime }]  \label{e3.6} \\
\ \leq \mathcal{E}_{s,t}^{g_\mu }[X-X^{\prime };(K-K^{\prime })_{\cdot }]%
\hbox{, a.s.}  \nonumber
\end{eqnarray}
\begin{equation}
\mathcal{E}_{s,t}^{-g_\mu }[0;K+K^0]\leq \mathcal{E}_{s,t}[0;K]\leq \mathcal{%
E}_{s,t}^{g_\mu }[0;K+K^0]  \label{e3.5}
\end{equation}
\end{lemma}

\smallskip\noindent\textbf{Proof. }We only prove (\ref{e3.6}). The proof of (\ref{e3.5}%
) is similar. Without loss of generality, we can set $K_t=\sum_{i=0}^N\xi
_i1_{[t_i,t_{i+1})}(t)$ and $K_t^{\prime }=\sum_{i=0}^N\xi _i^{\prime
}1_{[t_i,t_{i+1})}(t)$ for some $0=t_0<t_1\ldots <t_N=T$. When $t_i\leq
s\leq t\leq t_{i+1}$, we have, since $\mathcal{E}[\cdot ]$ satisfies (A5)
\begin{eqnarray*}
\mathcal{E}_{s,t}[X;K_{\cdot }]-\mathcal{E}_{s,t}[X^{\prime };K_{\cdot
}^{\prime }] &=&\mathcal{E}_{s,t}[X+K_t-K_s]-\mathcal{E}_{s,t}[X^{\prime
}+K_t^{\prime }-K_s^{\prime }] \\
&\leq &\mathcal{E}_{s,t}^{g_\mu }[(X-X^{\prime })+(K_t-K_t^{\prime
})-(K_s-K_s^{\prime })] \\
&=&\mathcal{E}_{s,t}^{g_\mu }[X-X^{\prime };(K-K^{\prime })_{\cdot }].
\end{eqnarray*}
Now let $t_{i-1}\leq s\leq t_i\leq t\leq t_{i+1}$, for some $1\leq i\leq N-1$%
, we have
\begin{eqnarray*}
\mathcal{E}_{s,t}[X;K_{\cdot }]-\mathcal{E}_{s,t}[X^{\prime };K_{\cdot
}^{\prime }] &=&\mathcal{E}_{s,t_i}[\mathcal{E}_{t_i,t}[X;K_{\cdot
}];K_{\cdot }]-\mathcal{E}_{s,t_i}[\mathcal{E}_{t_i,t}[X^{\prime };K_{\cdot
}^{\prime }];K_{\cdot }^{\prime }] \\
&\leq &\mathcal{E}_{s,t_i}^{g_\mu }[\mathcal{E}_{t_i,t}[X;K_{\cdot }]-%
\mathcal{E}_{t_i,t}[X^{\prime };K_{\cdot }^{\prime }];(K-K^{\prime })_{\cdot
}] \\
&\leq &\mathcal{E}_{s,t_i}^{g_\mu }[\mathcal{E}_{t_i,t}^{g_\mu }[X-X^{\prime
};(K-K^{\prime })_{\cdot }];(K-K^{\prime })_{\cdot }] \\
&=&\mathcal{E}_{s,t}^{g_\mu }[X-X^{\prime };(K-K^{\prime })_{\cdot }].
\end{eqnarray*}
We can repeat this procedure to prove that
\[
\mathcal{E}_{s,t}[X;K_{\cdot }]-\mathcal{E}_{s,t}[X^{\prime };K_{\cdot
}^{\prime }]\leq \mathcal{E}_{s,t}^{g_\mu }[X-X^{\prime };(K-K^{\prime
})_{\cdot }],\;\forall 0\leq s\leq t\leq T.
\]
We then have obtained the second inequality of (\ref{e3.6}). The first
inequality is obtained by changing the positions of $(X,K)$ and $(X^{\prime
},K^{\prime })$ in the second inequality of (\ref{e3.6}) and by observing
that
\begin{equation}
-\mathcal{E}_{s,t}^{g_\mu }[X^{\prime }-X;(K^{\prime }-K)_{\cdot }]=\mathcal{%
E}_{s,t}^{-g_\mu }[X-X^{\prime };(K-K^{\prime })_{\cdot }].  \label{e3.7}
\end{equation}
The proof is complete. $\Box $\medskip\

\begin{corollary}
\label{m3.6}We have the following estimate
\begin{eqnarray}
\ \ \ \ \ \ \ \ E[\sup_{s\in [0,t]}|\mathcal{E}_{s,t}[X;K_{\cdot }]+K_s-(%
\mathcal{E}_{s,t}[X^{\prime };K_{\cdot }^{\prime }]+K_s^{\prime })|^2]
\label{e3.8} \\
\ \leq CE[(X-X^{\prime }+K_t-K_t^{\prime })^2]+CE\int_0^t(K_s-K_s^{\prime
})^2ds.  \nonumber
\end{eqnarray}
where $C$ only depends on $T$ and $\mu $.
\end{corollary}

\smallskip\noindent\textbf{Proof. }Since both $g_\mu $ and $-g_\mu $ satisfies
conditions (\ref{h2.1}) for $g$. We set $y_s^1:=\mathcal{E}_{s,t}^{g_\mu
}[X^{\prime }-X;(K^{\prime }-K)_{\cdot }]$ and $y_s^2=\mathcal{E}%
_{s,t}^{-g_\mu }[X^{\prime }-X;(K^{\prime }-K)_{\cdot }]$. We observe that
\[
\mathcal{E}_{s,t}^{g_\mu }[0;0]=\mathcal{E}_{s,t}^{-g_\mu }[0;0]\equiv 0.
\]
We then can apply estimate (\ref{e2.2}) with $\tau =t$, to get, for $i=1,2$,

\begin{eqnarray*}
&&\ \ \ \ E[\sup_{s\in [0,t]}|y_s^i+K_s-K_s^{\prime }|^2] \\
\ &\leq &CE[(X-X^{\prime }+K_t-K_t^{\prime })^2]+CE\int_0^t(K_s-K_s^{\prime
})^2ds.
\end{eqnarray*}
This with (\ref{e3.6}) derives (\ref{e3.8}). The proof is complete. $\Box $%
\medskip\

For each $K\in D_{\mathcal{F}}^2(0,T)$ and for each $0\leq s\leq t\leq T$, $%
X\in L^2(\mathcal{F}_t)$, we take a sequence $\{K^i\}_{i=1}^\infty $ in $D_{%
\mathcal{F}}^{2,0}(0,T)$ such that $\{K^i\}_{i=1}^\infty $ converges in $L_{%
\mathcal{F}}^2(0,T)$ to $K$ and such that $K_s^i=K_s$, $K_t^i=K_t$. It
follows that $\{\mathcal{E}_{s,t}[X;K_{\cdot }^i]+K_s^i\}_{i=1}^\infty $ is
a Cauchy sequence in $L^2(\mathcal{F}_s)$. Consequently, $\{\mathcal{E}%
_{s,t}[X;K_{\cdot }^i]\}_{i=1}^\infty $ is a Cauchy sequence in $L^2(%
\mathcal{F}_s)$.

\begin{remark}
A sequence $\{K^i\}_{i=1}^\infty $ satisfying the above condition can be
realized by, for example, taking $0=t_0^i<t_1^i<\cdots t_i^i=T$, $%
\max_j(t_{j+1}^i-t_j^i)\rightarrow 0$, with $s=t_{j_1}^i$ and $t=t_{j_2}^i$
for some $j_1\leq j_2\leq i$, and then define
\[
K_t^i:=\sum_{j=0}^iK_{t_j^i}1_{[t_j^i,t_{j+1}^i)}(t),\;t\in [0,T].
\]
\end{remark}

\begin{definition}
\label{m3.7}We denote the limit of the Cauchy sequence $\{\mathcal{E}%
_{s,t}[X;K_{\cdot }^i]\}_{i=1}^\infty $ in $L^2(\mathcal{F}_s)$ by $\mathcal{%
E}_{s,t}[X;K_{\cdot }]$.
\end{definition}

The following property still holds true for $K\in D_{\mathcal{F}}^2(0,T)$.

\begin{proposition}
\label{m3.8}We assume (A1)--(A5) as well as (\ref{e3.1a}). Then $\mathcal{E}%
[\cdot ;K]$ is dominated by $\mathcal{E}^{g_\mu }$ in the following sense,
for each $K$, $K^{\prime }\in D_{\mathcal{F}}^2(0,T)$ and for each $0\leq
s\leq t\leq T$, $X$, $X^{\prime }\in L^2(\mathcal{F}_t)$, we have
\begin{eqnarray}
\mathcal{E}_{s,t}^{-g_\mu }[X-X^{\prime };(K-K^{\prime })_{\cdot }]\leq
\mathcal{E}_{s,t}[X;K_{\cdot }]-\mathcal{E}_{s,t}[X^{\prime };K_{\cdot
}^{\prime }]  \label{e3.10} \\
\ \leq \mathcal{E}_{s,t}^{g_\mu }[X-X^{\prime };(K-K^{\prime })_{\cdot }]%
\hbox{, a.s.}  \nonumber
\end{eqnarray}
\begin{equation}
\mathcal{E}_{s,t}^{-g_\mu }[0;K_{\cdot }+K^0]\leq \mathcal{E}%
_{s,t}[0;K_{\cdot }]\leq \mathcal{E}_{s,t}^{g_\mu }[0;K_{\cdot }+K^0]
\label{e3.10a}
\end{equation}
\end{proposition}

\smallskip\noindent\textbf{Proof. }Let $\{K^i\}_{i=1}^\infty $ and $\{K^{^{\prime
}i}\}_{i=1}^\infty $ be sequences in $D_{\mathcal{F}}^{2,0}(0,T)$ satisfying
the conditions of Definition \ref{m3.7} for $K$ and $K^{\prime }$,
respectively. By lemma \ref{m3.5}, we have
\begin{eqnarray}
\mathcal{E}_{s,t}^{-g_\mu }[X-X^{\prime };(K^i-K^{\prime i})_{\cdot }] &\leq
&\mathcal{E}_{s,t}[X;K_{\cdot }^i]-\mathcal{E}_{s,t}[X^{\prime };K_{\cdot
}^{\prime i}]  \label{e3.11} \\
\ &\leq &\mathcal{E}_{s,t}^{g_\mu }[X-X^{\prime };(K^i-K^{\prime
j})_{\cdot }]\hbox{, a.s.}  \nonumber
\end{eqnarray}
We then pass to the limit to get (\ref{e3.10}). The proof of (\ref{e3.10a})
is similar. $\Box $\medskip\

Exactly as Corollary \ref{m3.6}, we have immediately the following result.

\begin{corollary}
\label{m3.9}For $K$, $K^{\prime }\in D_{\mathcal{F}}^2(0,T)$, we have
\begin{eqnarray}
\ \ \ \ \ \ E[\sup_{s\in [0,t]}|\mathcal{E}_{s,t}[X;K_{\cdot }]+K_s-(%
\mathcal{E}_{s,t}[X^{\prime };K_{\cdot }^{\prime }]+K_s^{\prime })|^2]
\label{e3.12} \\
\ \leq CE[(X-X^{\prime }+K_t-K_t^{\prime })^2]+CE\int_0^t(K_s-K_s^{\prime
})^2ds.  \nonumber
\end{eqnarray}
\end{corollary}

The following properties comes immediately from Lemma \ref{m3.2} and Lemma
\ref{m3.5}.

\begin{proposition}
\label{m3.10}For a given $K\in D_{\mathcal{F}}^2(0,T)$, $\mathcal{E}%
_{s,t}[\cdot ;K]$ is an $\mathcal{F}_t$--consistent evaluation: (A1)--(A5),
i.e., \\ \textbf{(A1)} $\mathcal{E}_{s,t}[X;K]\geq \mathcal{E}%
_{s,t}[X^{\prime };K]$, a.s., if $X\geq X^{\prime }$, a.s. ;\\ \textbf{(A2)}
$\mathcal{E}_{t,t}[X;K]=X$; \\ \textbf{(A3)} $\mathcal{E}_{r,s}[\mathcal{E}%
_{s,t}[X;K];K]=\mathcal{E}_{r,t}[X;K]$, $\forall 0\leq r\leq s\leq t$;\\
\textbf{(A4) }$1_A\mathcal{E}_{s,t}[X;K]=1_A\mathcal{E}_{s,t}[1_AX;K]$, $%
\forall A\in \mathcal{F}_s$; \\ \textbf{(A5) }for each $K$, $K^{\prime }\in
D_{\mathcal{F}}^2(0,T)$, inequalities (\ref{e3.10a}) and (\ref{e3.10}) hold
true.
\end{proposition}

\begin{proposition}
\label{m3.11}If there is some function $g$ satisfying (\ref{h2.1}) such that
$\mathcal{E}$ coincides with $\mathcal{E}^g$, i.e., for each $0\leq s\leq
t\leq T$ and $X\in L^2(\mathcal{F}_T)$ we have $\mathcal{E}_{s,t}[X]=%
\mathcal{E}_{s,t}^g[X]$, then, for each $K\in D_{\mathcal{F}}^2(0,T)$, $%
\mathcal{E}_{s,t}[\cdot ;K]$ also coincides with $\mathcal{E}_{s,t}^g[\cdot
;K]$.
\end{proposition}

\subsection{Two corollaries from Theorem \ref{m7.1}}

The situation without assumption (A4$_0$) can be derived by Theorem \ref
{m7.1}:

\begin{corollary}
\label{m7.1a} Let $\mathcal{E}_{s,t}[\cdot ]:L^2(\mathcal{F}_t)\rightarrow
L^2(\mathcal{F}_s)$, $0\leq s\leq t\leq T$, satisfy (A1)--(A5) and
\begin{equation}
\mathcal{E}_{s,t}^{-g_\mu }[0;K^0]\leq \mathcal{E}_{s,t}[0]\leq \mathcal{E}%
_{s,t}^{g_\mu }[0;K^0],  \label{e3.1a}
\end{equation}
with a given $K^0\in D_{\mathcal{F}}^2(0,T)$. Then there exists a function $%
g(\omega ,t,y,z)$ satisfying (\ref{h2.1}) with $g(s,0,0)\equiv 0$, such that
\begin{equation}
\mathcal{E}_{s,t}[X]=\mathcal{E}_{s,t}^g[X;K^0],  \label{e7.1a}
\end{equation}
\end{corollary}

\smallskip\noindent\textbf{Proof. }\textbf{\ }We define $\mathcal{E}_{s,t}^0[X]:=%
\mathcal{E}_{s,t}[X;-K_{\cdot }^0]$. The above proposition ensures that $%
\mathcal{E}^0[\cdot ]$ satisfies all conditions (A1)--(A5). Moreover, by (%
\ref{e3.10a}),
\[
0=\mathcal{E}_{s,t}^{-g_\mu }[0;K^0-K^0]\leq
\mathcal{E}_{s,t}^0[0]\leq \mathcal{E}_{s,t}^{g_\mu
}[0;K^0-K^0]=0\hbox{, a.s..}
\]
Thus $\mathcal{E}_{s,t}^0[\cdot ]$ also satisfies (A4$_0$). It follows from
Theorem \ref{m7.1} that there exists $g$ with $g(s,0,0)\equiv 0$, such that
\[
\mathcal{E}_{s,t}^0[X]=\mathcal{E}_{s,t}^g[X]
\]
or
\[
\mathcal{E}_{s,t}[X;-K_{\cdot }^0]=\mathcal{E}_{s,t}^g[X].
\]
Consequently
\[
\mathcal{E}_{s,t}[X]=\mathcal{E}_{s,t}^g[X;K^0].
\]
\hfill$\Box $\medskip

\begin{corollary}
\label{m7.1b}Let $\mathcal{E}_{s,t}[\cdot ]:L^2(\mathcal{F}_t)\rightarrow
L^2(\mathcal{F}_s)$, $0\leq s\leq t\leq T$, satisfy (A1)--(A5) and
\begin{equation}
\mathcal{E}_{s,t}^{-g_\mu +g^0}[0]\leq \mathcal{E}_{s,t}[0]\leq \mathcal{E}%
_{s,t}^{g_\mu +g^0}[0],  \label{e3.1b}
\end{equation}
with a given $g^0\in L_{\mathcal{F}}^2(0,T)$. Then there exists a function $%
g(\omega ,t,y,z)$ satisfying (\ref{h2.1}) with $g(s,0,0)\equiv g_s^0$, such
that, for each $0\leq s\leq t\leq T$,
\begin{equation}
\mathcal{E}_{s,t}[X]=\mathcal{E}_{s,t}^g[X],\;\;\forall X\in L^2(\mathcal{F}%
_t).  \label{e7.1b}
\end{equation}
\end{corollary}

\smallskip\noindent\textbf{Proof. }We set $K_t^0:=\int_0^tg_s^0ds$, $t\in [0,T]$. By
the definition of $\mathcal{E}^g[\cdot ;K]$, condition (\ref{e3.1b}) reads
as
\[
\mathcal{E}_{s,t}^{-g_\mu }[0;K^0]\leq \mathcal{E}_{s,t}[0]\leq \mathcal{E}%
_{s,t}^{g_\mu }[0;K^0].
\]
It follows from Corollary \ref{m7.1a} that there exists a function $\bar
g(\omega ,t,y,z)$ satisfying (\ref{h2.1}) with $\bar g(s,0,0)\equiv 0$, such
that, for each $0\leq s\leq t\leq T$ and $X\in L^2(\mathcal{F}_t)$ we have $%
\mathcal{E}_{s,t}[X]=\mathcal{E}_{s,t}^{\bar g}[X;K^0]$, or equivalently, $%
\mathcal{E}_{s,t}[X]=\mathcal{E}_{s,t}^g[X]$, where we set $g(s,y,z):=\bar
g(s,y,z)+g_s^0$. The proof is complete. \hfill$\Box $\medskip

\section{$\mathcal{E}[\cdot ;K]$--martingales \label{ss5}}

Hereinafter, $\mathcal{E}[\cdot ]$ will be a fixed $\mathcal{F}_t$%
--consistent evaluation satisfying (A1)--(A5) and (A4$_0$). We introduce the
notion of $\mathcal{E}[\cdot ;K]$--martingale:

\begin{definition}
\label{m4.1}Let $K\in D_{\mathcal{F}}^2(0,T)$ be given. A process $Y\in L_{%
\mathcal{F}}^2(t_0,t_1)$ satisfying $E[\mathrm{ess}\sup_{s\in
[t_1,t_1]}|Y_s|^2]<\infty $, is said to be an $\mathcal{E}[\cdot ;K]$%
--martingales (resp. $\mathcal{E}[\cdot ;K]$--supermartingale, $\mathcal{E}%
[\cdot ;K]$--submartingale) on $[t_0,t_1]$ if for each $t_0\leq s\leq t\leq
t_1$, we have
\begin{equation}
\mathcal{E}_{s,t}[Y_t;K]=Y_s\hbox{, (resp. }\leq Y_s\hbox{, }\geq Y_s\hbox{%
), a.s.}  \label{e4.1}
\end{equation}
\end{definition}

\begin{proposition}
\label{m4.2}We assume (A1)--(A5) and (A4$_0$). Let $K$ and $K^{\prime }\in
D_{\mathcal{F}}^2(0,T)$ be given. Then for each fixed $t_1\in [0,T]$ and $X$%
, $X^{\prime }\in L^2(\mathcal{F}_{t_1})$, the process defined by
\begin{equation}
Y_s^{t_1,X,K}:=\mathcal{E}_{s,t_1}[X;K],\;s\in [0,t_1]  \label{e4.2}
\end{equation}
is an $\mathcal{E}[\cdot ;K]$--martingale, an $\mathcal{E}^{g_\mu }[\cdot
;K] $--submartingale as well as an $\mathcal{E}^{-g_\mu }[\cdot ;K]$%
--supermartingale. The difference of the processes
\begin{equation}
\ Y_s=\mathcal{E}_{s,t_1}[X;K]-\mathcal{E}_{s,t_1}[X^{\prime };K^{\prime
}],\;s\in [0,t_1]  \label{e4.3}
\end{equation}
is also an $\mathcal{E}^{g_\mu }[\cdot ;K-K^{\prime }]$--submartingale and
an $\mathcal{E}^{-g_\mu }[\cdot ;K-K^{\prime }]$--supermartingale.
\end{proposition}

\smallskip\noindent\textbf{Proof. }The first assertion comes directly from (A3) of
Proposition \ref{m3.10}. Now, for each $0\leq s\leq t\leq t_1$,
\begin{eqnarray*}
Y_s &=&\mathcal{E}_{s,t_1}[X;K]-\mathcal{E}_{s,t_1}[X^{\prime };K^{\prime }]
\\
\ &=&\mathcal{E}_{s,t}[\mathcal{E}_{t,t_1}[X;K];K]-\mathcal{E}_{s,t}[%
\mathcal{E}_{t,t_1}[X^{\prime };K^{\prime }];K^{\prime }] \\
\ &\leq &\mathcal{E}_{s,t}^{g_\mu }[\mathcal{E}_{t,t_1}[X;K]-\mathcal{E}%
_{t,t_1}[X^{\prime };K^{\prime }];K-K^{\prime }] \\
\ &\leq &\mathcal{E}_{s,t}^{g_\mu }[Y_t;K-K^{\prime }].
\end{eqnarray*}
Thus $Y$ is an $\mathcal{E}^{g_\mu }[\cdot ;K-K^{\prime }]$ submartingale on
$[0,t_1]$. Similarly we can prove that it is an $\mathcal{E}^{-g_\mu }[\cdot
;K-K^{\prime }]$--supermartingale. $\Box $\medskip\

We will prove that $\mathcal{E}_{s,t_1}[X;K]$, $s\in [0,t_1]$ have an RCLL
modification. The following upcrossing inequality can be found in \cite
{Peng2003b}. We denote $\mathcal{E}^\mu [\cdot ]:=\mathcal{E}_{0,T}^{g_\mu
^{*}}[\cdot ]$, $\mathcal{E}^{-\mu }[\cdot ]:=\mathcal{E}_{0,T}^{-g_\mu
^{*}}[\cdot ]$ with $g_\mu ^{*}(z)=\mu |z|$, $z\in R^d$.

\begin{theorem}
We assume that $g$ satisfies (i) and (ii) of (\ref{h2.1}). Let $%
Y=(Y_t)_{t\in [0,T]}$ be a $\mathcal{E}^g$--supermartingale, $D$ be a
denumerable dense subset of $[0,T]$. Then for each $a,b\in R$, $r,s\in [0,T]$
such that $a<b$ and $r<s$, we have
\begin{equation}
\mathcal{E}^{-\mu }[U_a^b(Y,D\cap [r,s])\leq \frac{e^{2\mu (s-r)}}{b-a}\{%
\mathcal{E}^\mu [(Y_s-a)^{-}]+\mathcal{E}^\mu [\int_r^se^{\mu
t}|g_t^0|dt]+a\mu (s-r)\},\   \label{eA10}
\end{equation}
where $\mu $ is the Lipschitz constant of $g$ and $g_s^0=g(s,0,0)$. In
particular
\begin{equation}
\mathcal{E}^{-\mu }[U_a^b(Y,D)]\leq \frac{e^{2\mu T}}{b-a}\{\mathcal{E}^\mu
[(Y_T-a)^{-}]+\mathcal{E}^\mu [\int_0^Te^{\mu t}|g_t^0|dt]\ +a\mu T\}.
\label{upcross3}
\end{equation}
Moreover, $U_a^b(Y,D)<\infty $, a.s.
\end{theorem}

\begin{proposition}
For each $X\in L^2(\mathcal{F}_{t_1})$, the process $(\mathcal{E}%
_{t,t_1}[X])_{0\leq t\leq t_1}$ has an RCLL modification.
\end{proposition}

\smallskip\noindent\textbf{Proof. }Without loss of generality, we set $t_1=T$. Since $%
\mathcal{E}_{\cdot ,T}[X]$ is an $\mathcal{E}^{-g_\mu }$--supermartingale,
by upcrossing inequality, it is classical that the $\mathcal{F}_t$--adapted
process $\bar Y$ defined by
\[
\bar Y_t:=\lim_{s\in \mathbf{Q}\cap (0,T],s\searrow t}\mathcal{E}%
_{s,T}[X],\;t\in [0,T).
\]
is RCLL. Thus it suffices to prove that, for each $t\in [0,T)$, $\bar Y_t=$ $%
\mathcal{E}_{t,T}[X]$, a.s. Indeed, let $\{s_n\}_{n=1}^\infty \subset \cap
(t,T]$ be such that $s_n\searrow t$. By $\mathbf{E}[\sup_{0\leq s\leq T}|%
\mathcal{E}_{s,T}[X]|^2]<\infty $ and $\mathcal{E}_{s_n,T}[X;K]\rightarrow
\bar Y_t$ it follows that
\[
\mathcal{E}_{s_n,T}[X]\rightarrow \bar Y_t\hbox{, in
}L^2(\mathcal{F}_T).
\]
Thus for each $A\in \mathcal{F}_t$, $\mathcal{E}_{s_n,T}[X]1_A\rightarrow
\bar Y_t1_A$, in $L^2(\mathcal{F}_T).$ Or
\[
\mathcal{E}_{s_n,T}[X1_A]\rightarrow \bar Y_t1_A,\hbox{ in }L^2(\mathcal{F}%
_T).
\]
Thus, by (\ref{e2.2}) with $K=K^{\prime }=0$, $\tau =T$, $X=\mathcal{E}%
_{s_n,T}[X1_A]$ and $X^{\prime }=\bar Y_t1_A$, we have
\begin{eqnarray*}
&&\ |\mathcal{E}_{0,s_n}[\mathcal{E}_{s_n,T}[X1_A]]-\mathcal{E}_{0,s_n}[\bar
Y_t1_A]|^2 \\
\ &\leq &\mathbf{E}[|\mathcal{E}_{s_n,T}[X1_A]-\bar Y_t1_A|^2]\rightarrow 0.
\end{eqnarray*}
We also have
\begin{eqnarray*}
|\mathcal{E}_{0,s_n}[\bar Y_t1_A]-\mathcal{E}_{0,t}[\bar Y_t1_A]|^2 &=&|%
\mathcal{E}_{0,t}[\mathcal{E}_{t,s_n}[\bar Y_t1_A]]-\mathcal{E}_{0,t}[\bar
Y_t1_A]|^2 \\
\ &\leq &c\mathbf{E}[|\mathcal{E}_{t,s_n}[\bar Y_t1_A]-\bar Y_t1_A]|^2]
\end{eqnarray*}
It follows from (\ref{e8.29}) that $\mathcal{E}_{0,s_n}[\bar
Y_t1_A]\rightarrow \mathcal{E}_{0,t}[\bar Y_t1_A]$. The above two
convergences imply
\[
\mathcal{E}_{0,s_n}[\mathcal{E}_{s_n,T}[X1_A]]\rightarrow \mathcal{E}%
_{0,T}[\bar Y_t1_A].
\]
But on the other hand, $\mathcal{E}_{0,s_n}[\mathcal{E}_{s_n,T}[X1_A]]=%
\mathcal{E}_{0,t}[\mathcal{E}_{t,T}[X]1_A]$. We thus have
\[
\mathcal{E}_{0,t}[\mathcal{E}_{t,T}[X]1_A]=\mathcal{E}_{0,T}[\bar
Y_t1_A],\;\forall A\in \mathcal{F}_t.
\]
From which it follows that $\mathcal{E}_{t,T}[X]=\bar Y$, a.s. \hfill$\Box $%
\medskip

We will always take an RCLL version of $\mathcal{E}_{\cdot ,t}[\cdot ]$.

\begin{proposition}
We assume (A1)--(A5) and (A4$_0$). Then, for each $X\in L^2(\mathcal{F}%
_{t_1})$ and $K\in D_{\mathcal{F}}^2(0,T)$ the process $(\mathcal{E}%
_{t,t_1}[X;K])_{0\leq t\leq t_1}$ belongs to $D_{\mathcal{F}}^2(0,t_1)$. %
\hfill$\Box $
\end{proposition}

\smallskip\noindent\textbf{Proof. }In the case where $K\in D_{\mathcal{F}}^{2,0}(0,T)$%
, from the definition and the above lemma it follows that $\mathcal{E}%
_{\cdot ,t_1}[X;K]$ is also RCLL. This with (\ref{e3.12}) (by setting $%
K^{\prime }\equiv 0$, $X^{\prime }=0$) we deduce that $\mathcal{E}_{\cdot
,t_1}[X;K]\in D_{\mathcal{F}}^2(0,t_1)$. Now let $K\in D_{\mathcal{F}%
}^2(0,T) $ and let $\{K^i\}_{i=1}^\infty $ be a sequence in $D_{\mathcal{F}%
}^{2,0}(0,T)$ such that $K^i\rightarrow K$ in $L_{\mathcal{F}}^2(0,T)$ and $%
K_{t_1}^i\rightarrow K_{t_1}$ in $L^2(\mathcal{F}_T)$. From (\ref{e3.12}) it
follows that
\begin{eqnarray*}
&&\ \ \ \ \ \ \ E[\sup_{s\in [0,t_1]}|\mathcal{E}_{s,t_1}[X;K_{\cdot
}^i]+K_s^i-(\mathcal{E}_{s,t_1}[X;K_{\cdot }]+K_s)|^2] \\
\ &\leq &CE[(K_{t_1}^i-K_{t_1})^2]+CE\int_0^{t_1}(K_s-K_s^{\prime
})^2ds\rightarrow 0.
\end{eqnarray*}
Since $(\mathcal{E}_{t,t_1}[X;K^i]-K_t^i)_{0\leq t\leq t_1}$ $i=1,2,\cdots $
are in $D_{\mathcal{F}}^2(0,t_1)$, $(\mathcal{E}_{t,t_1}[X;K]-K_t)_{0\leq
t\leq t_1}$ and thus $(\mathcal{E}_{t,t_1}[X;K])_{0\leq t\leq t_1}$are also
in $D_{\mathcal{F}}^2(0,t_1)$. \hfill$\Box $

We then can apply $\mathcal{E}^g$--supermartingale decomposition theorem,
i.e., Proposition \ref{p2.3}, to get the following result.

\begin{proposition}
\label{m4.3}We assume (A1)--(A5) and (A4$_0$). Let $K$ $\in D_{\mathcal{F}%
}^2(0,T)$ be given. For fixed $t\in [0,T]$ and $X\in L^2(\mathcal{F}_t)$,
the process $Y_s^{t,X,K}:=\mathcal{E}_{s,t}[X;K],\;s\in [0,t]$, has the
following expression: there exist processes $(g_{\cdot }^{t,X,K},z_{\cdot
}^{t,X,K})\in L_{\mathcal{F}}^2(0,t;R\times R^d)$ such that
\begin{equation}
Y_s^{t,X,K}=X+K_t-K_s+\int_s^tg_r^{t,X,K}dr-\int_s^tz_r^{t,X,K}dB_r,\;\;s\in
[0,t],  \label{e4.4}
\end{equation}
such that
\begin{equation}
|g_s^{t,X,K}|\leq \mu (|Y_s^{t,X,K}|+|z_s^{t,X,K}|),\;\forall s\in [0,t].
\label{e4.5}
\end{equation}
Moreover let $Y_s^{t,X^{\prime },K^{\prime }}:=\mathcal{E}_{s,t_1}[X^{\prime
};K^{\prime }]$, $s\in [0,t]$, for some other $K^{\prime }\in D_{\mathcal{F}%
}^2(0,T)$, $X^{\prime }\in L^2(\mathcal{F}_t)$ and let $(g_{\cdot
}^{t,X^{\prime },K^{\prime }},z_{\cdot }^{t,X^{\prime },K^{\prime }})$ be
the corresponding expression in (\ref{e4.4}), then we have
\begin{equation}
|g_s^{t,X,K}-g_s^{t,X^{\prime },K^{\prime }}|\leq \mu
(|Y_s^{t,X,K}-Y_s^{t,X^{\prime },K^{\prime }}|+|z_s^{t,X,K}-z_s^{t,X^{\prime
},K^{\prime }}|),\;\forall s\in [0,t].  \label{e4.6}
\end{equation}
\end{proposition}

\smallskip\noindent\textbf{Proof. }Since $(Y_s^{t,X,K})_{s\in [0,t]}$, is an $\mathcal{%
E}^{g_\mu }[\cdot ;K]$--submartingale and $\mathcal{E}^{-g_\mu }[\cdot ;K]$%
--super--martingale, by Proposition \ref{p2.3} and Corollary \ref{c2.2},
there exists an increasing process $A_{\cdot }^{+}\in D_{\mathcal{F}}^2(0,t)$
and $A_{\cdot }^{-}\in D_{\mathcal{F}}^2(0,t)$ with $A_0^{+}=A_0^{-}=0$,
such that
\begin{equation}
Y_s^{t,X,K}=\mathcal{E}_{s,t}^{g_\mu }[X;(K-A^{+})_{\cdot }]=\mathcal{E}%
_{s,t}^{g_\mu }[X;(K+A^{-})_{\cdot }],\;s\in [0,t].  \label{e4.7}
\end{equation}
According to the notion of $\mathcal{E}^g$ defined in (\ref{e2.3}), $%
Y_s^{t,X,K}$ is the solution of the following BSDE on $[0,t]$:
\begin{eqnarray}
Y_s^{t,X,K} &=&X+(K-A^{+})_t-(K-A^{+})_s  \label{e4.8} \\
&&\ +\int_s^t\mu (|Y_r^{t,X,K}|+|Z_r^{+}|)dr-\int_s^tZ_r^{+}dB_r  \nonumber
\end{eqnarray}
and
\begin{eqnarray}
Y_s^{t,X,K} &=&X+(K+A^{-})_t-(K+A^{-})_s  \label{e4.9} \\
&&-\int_s^t\mu (|Y_r^{t,X,K}|+|Z_r^{-}|)dr-\int_s^tZ_r^{-}dB_r.  \nonumber
\end{eqnarray}
It then follows that $Z_s^{t,X,K}:=Z_s^{+}\equiv Z_s^{-}$, $s\in [0,t]$ and
thus
\[
-dA_s^{+}+\mu (|Y_s^{t,X,K}|+|Z_s^{t,X,K}|)ds\equiv dA_s^{-}-\mu
(|Y_s^{t,X,K}|+|Z_s^{t,X,K}|)ds,
\]
or
\begin{equation}
dA_s^{-}+dA_s^{+}\equiv 2\mu (|Y_s^{t,X,K}|+|Z_s^{t,X,K}|)ds,\;s\in [0,t]
\label{e4.10}
\end{equation}
Thus $dA^{+}$ and $dA^{-}$ are absolutely continuous with respect to $ds$.
We denote $a_s^{+}ds=dA_s^{+}$ and $a_s^{-}ds=dA_s^{-}$. It is clear that
\begin{eqnarray*}
0 &\leq &a_s^{+}\leq 2\mu (|Y_s^{t,X,K}|+|Z_s^{t,X,K}|), \\
0 &\leq &a_s^{-}\leq 2\mu (|Y_s^{t,X,K}|+|Z_s^{t,X,K}|),\;dP\times dt\hbox{%
--a.e.}
\end{eqnarray*}
We then can rewrite (\ref{e4.8}) as
\begin{equation}
Y_s^{t,X,K}=X+K_t-K_s+\int_s^t[-a_r^{+}+\mu
(|Y_r^{t,X,K}|+|Z_r^{+}|)]dr-\int_s^tZ_r^{+}dB_r.  \label{e4.11}
\end{equation}
Thus, by setting $g_r^{t,X,K}:=-a_r^{+}+\mu (|Y_r^{t,X,K}|+|Z_r^{+}|)$, we
have the expression (\ref{e4.4}) as well as the estimate (\ref{e4.5}).

It remains to prove (\ref{e4.6}). By (A5) of Proposition \ref{p2.3} $\hat
Y_s=Y_s^{t,X,K}-Y_s^{t,X^{\prime },K^{\prime }}$ is an $\mathcal{E}^{g_\mu
}[\cdot ;K-K^{\prime }]$--submartingale and an $\mathcal{E}^{-g_\mu }[\cdot
;K-K^{\prime }]$--supermartingale on $[0,t]$. Thus we can repeat the above
procedure to prove that there exist processes $(\hat g_{\cdot },\hat
Z_{\cdot })\in L_{\mathcal{F}}^2(0,t;R\times R^d)$ such that
\begin{equation}
\hat Y_s=X-X^{\prime }+(K-K^{\prime })_t-(K-K^{\prime })_s+\int_s^t\hat
g_rdr-\int_s^t\hat Z_rdB_r,\;\;s\in [0,t],  \label{e4.12}
\end{equation}
such that
\begin{equation}
|\hat g_s|\leq \mu (|\hat Y_s|+|\hat Z_s|),\;\forall s\in [0,t].
\label{e4.13}
\end{equation}
But by (\ref{e4.4}) and $\hat Y_s\equiv Y_s^{t,X,K}-Y_s^{t,X^{\prime
},K^{\prime }}$, we immediately have
\begin{equation}
\hat g_s\equiv g_s^{t,X,K}-g_s^{t,X^{\prime },K^{\prime }},\;\hat Z_s\equiv
z_s^{t,X,K}-z_s^{t,X^{\prime },K^{\prime }}.  \label{e4.14}
\end{equation}
This with (\ref{e4.13}) yields (\ref{e4.6}). The proof is complete. $\Box $%
\medskip\

\begin{corollary}
\label{m4.4}Let $K^1$ and $K^2\in D_{\mathcal{F}}^2(0,T)$ and $X^1\in L^2(%
\mathcal{F}_{t_1})$, $X^2\in L^2(\mathcal{F}_{t_2})$ be given for some fixed
$0\leq t_1\leq t_2\leq T$ and let $(g_s^{t_i,X^i,K^i},Z_s^{t_i,X^i,K^i})_{s%
\in [0,t_i]}$, $i=1,2$, be the pair in (\ref{e4.4}) for $Y_s^{t_i,X^i,K^i}=%
\mathcal{E}_{s,t_i}[X^i;(K^i)_{\cdot }]$, $i=1,2$, respectively. Then we
have
\begin{equation}
|g_s^{t,X^1,K^1}-g_s^{t,X^2,K^2}|\leq \mu
(|Y_s^{t_1,X^1,K^2}-Y_s^{t_2,X^2,K^2}|+|z_s^{t_1,X^1,K^1}-z_s^{t_2,X^2,K^2}|),\;\forall s\in [0,t_1].
\label{e4.15}
\end{equation}
\end{corollary}

\smallskip\noindent\textbf{Proof. }With the observation
\[
Y_s^{t_2,X^2,K^2}=\mathcal{E}_{s,t_1}[Y_{t_1}^{t_2,X^2,K^2};(K^2)_{\cdot
}],\;s\in [0,t_1],
\]
it is an immediate consequence of Proposition \ref{m4.3}. $\Box $\medskip\

\begin{corollary}
\label{m4.5}For each $t\in [0,T]$ and $X\in L^2(\mathcal{F}_t)$, $K\in D_{%
\mathcal{F}}^2(0,T)$, the process $(\mathcal{E}_{s,t}[X;K_{\cdot }])_{s\in
[0,t]}$ is also in $D_{\mathcal{F}}^2(0,t)$. If moreover, $K\in S_{\mathcal{F%
}}^2(0,T)$ (resp. It\^o's process), then $(\mathcal{E}_{s,t}[X;K_{\cdot
}])_{s\in [0,t]}$ is also in $S_{\mathcal{F}}^2(0,t)$ (resp. It\^o's
process).
\end{corollary}

\section{BSDE under $\mathcal{E}[\cdot ]$\label{ss6}}

We now consider the following kind of backward stochastic differential
equations: Let $X\in L^2(\mathcal{F}_T)$ be given and let
\begin{equation}
f:(\omega ,t,y)\in \Omega \times [0,T]\times R\rightarrow R  \label{e5.1}
\end{equation}
be a given function. We assume that $f$ satisfies
\begin{equation}
\left\{
\begin{array}{rrl}
\hbox{(i)} & f(\cdot ,y) & \in L_{\mathcal{F}}^2(0,T)\hbox{, for
each }y\in
R, \\
\hbox{(ii)} & |f(t,y) & -f(t,y^{\prime })|\leq c(|y-y^{\prime
}|),\;\forall y,y^{\prime }\in R.
\end{array}
\right.  \label{h5.1}
\end{equation}
We consider the following kind of BSDE
\begin{equation}
Y_t=\mathcal{E}_{t,T}[X;\int_0^{\cdot }f(s,Y_s)ds],\;t\in [0,T].
\label{e5.2}
\end{equation}
We have the following existence and uniqueness result

\begin{theorem}
\label{m5.1}Let $\mathcal{E}[\cdot ]$ satisfy (A1)--(A5) and (A4$_0$) and
let $f$ satisfies (\ref{h5.1}). Then for each $X\in L^2(\mathcal{F}_T)$,
there exists a unique solution $Y\in \mathcal{S}_{\mathcal{F}}^2(0,T)$ of
BSDE (\ref{e5.2}).
\end{theorem}

\smallskip\noindent\textbf{Proof. }With Corollary \ref{m4.5}, we only need to prove
that BSDE (\ref{e5.2}) has a unique solution $Y\in L_{\mathcal{F}}^2(0,T)$.
To this end we define a mapping $\Lambda [\cdot ]:L_{\mathcal{F}%
}^2(0,T)\rightarrow L_{\mathcal{F}}^2(0,T)$ by
\begin{equation}
\Lambda _t[y]:=Y_t=\mathcal{E}_{t,T}[X;\int_0^{\cdot }f(s,y_s)ds],\;y\in L_{%
\mathcal{F}}^2(0,T).  \label{e5.3}
\end{equation}
By Corollary \ref{m4.5}, $\Lambda _{\cdot }[y]$ belongs to $S_{\mathcal{F}%
}^2(0,T)$. From Proposition \ref{m3.8}, we have
\begin{equation}
\Lambda _t[y^1]-\Lambda _t[y^2]\leq \mathcal{E}_{t,T}^{g_\mu
}[0;\int_0^{\cdot }(f(s,y_s^1)-f(s,y_s^2))ds],\;\forall y^1,y^2\in L_{%
\mathcal{F}}^2(0,T).  \label{e5.4}
\end{equation}
By Proposition \ref{p2.0}, we have, with $\beta =\mu ^2+2\mu +1$,
\begin{eqnarray*}
E[(\Lambda _t[y^1]-\Lambda _t[y^2])^2] &\leq &E[(\mathcal{E}_{t,T}^{g_\mu
}[0;\int_0^{\cdot }(f(s,y_s^1)-f(s,y_s^2))ds)^2] \\
&\leq &E[\int_t^Te^{\beta (s-t)}(f(s,y_s^1)-f(s,y_s^2))^2ds] \\
&\leq &CE[\int_t^T|y_s^1-y_s^2|^2ds],
\end{eqnarray*}
where $C=c^2e^{\beta T}$ and $c$ is the Lipschitz constant of $f$. We
multiple $e^{2Ct}$ on both sides and integrate on $[0,T]$,
\begin{eqnarray*}
E[\int_0^Te^{2Ct}(\Lambda _t[y^1]-\Lambda _t[y^2])^2dt] &\leq
&CE[\int_0^T\int_t^Te^{2Ct}|y_s^1-y_s^2|^2dsdt] \\
&=&CE[\int_0^T(\int_0^se^{2Ct}dt)|y_s^1-y_s^2|^2ds] \\
&=&C(2C)^{-1}E[\int_0^T(e^{2Ct}-1)|y_s^1-y_s^2|^2ds].
\end{eqnarray*}
We thus have
\begin{equation}
E[\int_0^Te^{2Ct}(\Lambda _t[y^1]-\Lambda _t[y^2])^2dt]\leq \frac
12E[\int_0^Te^{2Ct}|y_t^1-y_t^2|^2dt].  \label{e5.5}
\end{equation}
We observe that the following two norms are equivalent to each others in $L_{%
\mathcal{F}}^2(0,T)$:
\[
\lbrack E\int_0^T|\phi _t|^2dt]^{1/2}\sim [E\int_0^Te^{2Ct}|\phi
_t|^2dt]^{1/2}.
\]
It follows that $\Lambda _t[\cdot ]$ is a contraction mapping. Thus there
exists a unique fixed point $Y\in L_{\mathcal{F}}^2(0,T)$, such that
\[
Y_t=\Lambda _t[Y]=\mathcal{E}_{t,T}[X;\int_0^{\cdot }f(s,Y_s)ds].
\]
The prove is complete. $\Box $\medskip\

\begin{proposition}
\label{m5.2}Let $Y\in S_{\mathcal{F}}^2(0,T)$ be the solution of BSDE (\ref
{e5.2}).\ Then, for each $t\in [0,T)$ we have
\end{proposition}

\begin{equation}
Y_s=\mathcal{E}_{s,t}[Y_t;\int_0^{\cdot }f(s,Y_s)ds],\;s\in [0,t].
\label{e5.6}
\end{equation}

\smallskip\noindent\textbf{Proof. }We set $K_t:=\int_0^tf(s,Y_s)ds$. By proposition
\ref{m4.2}, $Y_s=\mathcal{E}_{s,T}[X;K_{\cdot }]$ is an $\mathcal{E}%
_{s,T}[\cdot ;K_{\cdot }]$--martingale. Thus (\ref{e5.6}) holds. $\Box $%
\medskip\

\begin{proposition}
\label{m5.3}Let $X$, $X^{\prime }\in L^2(\mathcal{F}_T)$ and $\phi \in L_{%
\mathcal{F}}^2(0,T)$ be given. Let $Y\in \mathcal{S}_{\mathcal{F}}^2(0,T)$
be the solution of BSDE (\ref{e5.2}), $Y^{\prime }\in \mathcal{S}_{\mathcal{F%
}}^2(0,T)$ be the solution of
\begin{equation}
Y_t^{\prime }=\mathcal{E}_{t,T}[X^{\prime };\int_0^{\cdot }(f(s,Y_s^{\prime
})+\phi _s)ds].  \label{e5.7}
\end{equation}
Then $Y^{\prime }-Y$ is an $\mathcal{E}^{g_{\mu +c,\mu }+\phi }[\cdot ]$%
--submartingale and an $\mathcal{E}^{-g_{\mu +c,\mu }+\phi }[\cdot ]$%
--supermartingale on $[0,T]$, where $c\geq 0$ is the Lipschitz constant of $%
f $ with respect to $y$ and
\[
g_{\mu +c,\mu }(y,z):=(c+\mu )|y|+\mu |z|.
\]
\end{proposition}

\smallskip\noindent\textbf{Proof. }By Proposition \ref{m3.8} and the above
proposition, we have, for each $0\leq s\leq t\leq T$,
\begin{eqnarray*}
Y_s^{\prime }-Y_s &=&\mathcal{E}_{s,t}[Y_t^{\prime };\int_0^{\cdot
}(f(s,Y_s^{\prime })+\phi _s)ds]-\mathcal{E}_{s,t}[Y_t;\int_0^{\cdot
}f(s,Y_s)ds] \\
\ &\geq &\mathcal{E}_{s,t}^{-g_\mu }[Y_t^{\prime }-Y_t;\int_0^{\cdot
}((f(s,Y_s^{\prime })-f(s,Y_s)+\phi _s)ds].
\end{eqnarray*}
Thus $Y^{\prime }-Y$ is an $\mathcal{E}^{-g_\mu }[\cdot ;K]$
supermartingale, where $K_t:=\int_0^t((f(s,Y_s^{\prime })-f(s,Y_s)+\phi
_s)ds $. By $\mathcal{E}^g$--supermartingale decomposition theorem
(Corollary \ref{c2.2}), there exists an increasing process $A\in D_{\mathcal{%
F}}^2(0,T)$ with $A_0=0$, such that
\[
Y_t-Y_t^{\prime }=\mathcal{E}_{t,T}^{-g_\mu }[X^{\prime }-X;(K+A)_{\cdot }].
\]
$\hat Y_t:=Y_t-Y_t^{\prime }$ solves the following BSDE
\begin{eqnarray*}
\hat Y_t &=&X^{\prime }-X+A_T-A_t+\int_t^T((f(s,Y_s^{\prime })-f(s,Y_s)+\phi
_s)ds \\
&&-\int_t^T\mu (|\hat Y_s|+|\hat Z_s|)ds-\int_t^T\hat Z_sdB_s
\end{eqnarray*}
or, equivalently
\begin{eqnarray*}
\hat Y_t &=&X^{\prime }-X+(A+\hat A)_T-(A+\hat A)_t\  \\
&&\ -\int_t^T[(\mu +c)(|\hat Y_s|+\mu |\hat Z_s|+\phi _s)ds-\int_t^T\hat
Z_sdB_s.
\end{eqnarray*}
Here we set
\[
\hat A_t=\int_0^t[c|\hat Y_s|+f(s,Y_s^{\prime })-f(s,Y_s)]ds.
\]
It is clear that $\hat A$ and thus $A+\hat A$ is an increasing process. Thus
$\hat Y$ is an $\mathcal{E}^{-g_{c+\mu ,\mu }+\phi }$--supermartingle.
Analogously, we can prove that it is an $\mathcal{E}^{g_{c+\mu ,\mu }+\phi }$%
--submartingle. $\Box $\medskip\

\begin{corollary}
\label{m5.5}If $X^{\prime }\geq X$ and $\phi _s\geq 0$, $dt\times dP$--a.e.,
then we have $Y_t^{\prime }-Y_t\geq 0$, $dt\times dP$--a.e..
\end{corollary}

\section{$\mathcal{E}$--supermartingale decomposition theorem: intrinsic
formulation\label{ss7}}

Our objective of this section is to prove the following $\mathcal{E}$%
--supermartingale decomposition theorem of Doob--Meyer's type. Since $(%
\mathcal{E}_{s,t}[\cdot ])_{s\leq t}$ is abstract and nonlinear, it is
necessary to introduce the intrinsic form (\ref{e6.1}). This theorem plays
an important role in the proof of the main theorem of this paper. It can be
also considered as a generalization of Proposition \ref{p2.3}.

\begin{theorem}
\label{m6.1}We assume (A1)--(A5) as well as (A4$_0$). Let $Y\in S_{\mathcal{F%
}}^2(0,T)$ be an $\mathcal{E}[\cdot ]$--supermartingale. Then there exists
an increasing process $A\in S_{\mathcal{F}}^2(0,T)$ with $A_0=0$, such that $%
Y$ is an $\mathcal{E}[\cdot ;A]$--martingale, i.e.,
\begin{equation}
Y_t=\mathcal{E}_{t,T}[Y_T;A_{\cdot }],\;t\in [0,T].  \label{e6.1}
\end{equation}
\end{theorem}

\begin{remark}
\label{m6.1Rem1}This theorem has an interesting interpretation in finance:
the fact that $Y\in S_{\mathcal{F}}^2(0,T)$ is an $\mathcal{E}[\cdot ]$%
--supermartingale is equivalent to that there exists an increasing process $%
A $ such that $Y$ is the flow of the dynamical evaluation of the sum of the
final payoff $Y_T$ at $T$ plus the flow of the ``dividend'' $A$ during the
whole period $[0,T]$.
\end{remark}

\begin{remark}
\label{m6.1Rem2}In the case where $(\mathcal{E}_{s,t}[\cdot ])_{0\leq s\leq
t\leq T}$ is a system of linear mappings, (\ref{e6.1}) becomes
\[
Y_t+A_t=\mathcal{E}_{t,T}[Y_T+A_T],\;t\in [0,T],
\]
i.e., as in classical situation, $Y+A$ is an $\mathcal{E}[\cdot ]$%
--martingale. But, the intrinsic formulation that can be applied to
nonlinear situation is that $Y$ is an $\mathcal{E}[\cdot ;A]$--martingale.
\end{remark}

In order to prove this theorem, we need to extent the definition of $%
\mathcal{E}_{s,t}[\cdot ]$ from deterministic times $s$, $t$ $\in [0,T]$ to $%
\mathcal{E}_{\sigma ,\tau }[\cdot ]$ of stopping times $\sigma $ and $\tau $:

\begin{theorem}
\label{m8.13}There exists a unique extension of $(\mathcal{E}_{s,t}[\cdot
])_{0\leq s\leq t\leq T}$ to the $\mathcal{F}_\tau $--consistent nonlinear
evaluation:
\[
\mathcal{E}_{\sigma ,\tau }[\cdot ]:L^2(\mathcal{F}_\tau )\rightarrow L^2(%
\mathcal{F}_\sigma ),\;\sigma ,\tau \in \mathcal{S}_T,\;\sigma \leq \tau ,
\]
for each $X$, $X^{\prime }\in L^2(\mathcal{F}_\tau )$, \\\textbf{(A1)} $%
\mathcal{E}_{\sigma ,\tau }[X]\geq \mathcal{E}_{\sigma ,\tau }[X^{\prime }]$%
, a.s., if $X\geq X^{\prime }$, a.s.; \\\textbf{(A2)} $\mathcal{E}_{\tau
,\tau }[X]=X$, a.s.;\\\textbf{(A3)} $\mathcal{E}_{\sigma ,\rho }[\mathcal{E}%
_{\rho ,\tau }[X]]=\mathcal{E}_{\sigma ,\tau }[X]$, $\forall 0\leq \sigma
\leq \rho \leq \tau $, $\rho \in \mathcal{S}_T$;\\\textbf{(A4')} $1_A%
\mathcal{E}_{\sigma ,\tau }[X]=\mathcal{E}_{\sigma ,\tau }[1_AX]$,\ $\forall
A\in \mathcal{F}_\sigma $; \\\textbf{(A5)} $\mathcal{E}^{g_\mu }$%
--domination: $\forall X,X^{\prime }\in L^2(\mathcal{F}_\tau ),\;\;\forall
K,K^{\prime }\in \mathcal{S}_{\mathcal{F}}^2(0,T).$%
\begin{equation}
\mathcal{E}_{\sigma ,\tau }[X;K]-\mathcal{E}_{\sigma ,\tau }[X^{\prime
};K^{\prime }]\leq \mathcal{E}_{\sigma ,\tau }^{g_\mu }[X-X^{\prime
};K-K^{\prime }].\;\;  \label{e8.34}
\end{equation}
\end{theorem}

\begin{remark}
The ``unique extension'' is in the following sense: if the system $\mathcal{E%
}_{\sigma ,\tau }^{\prime }[\cdot ]:L^2(\mathcal{F}_\tau )\rightarrow L^2(%
\mathcal{F}_\sigma )$, $\sigma ,\tau \in \mathcal{S}_T,\;\sigma \leq \tau $,
satisfies also the above (A1)--(A3), (A4') and (A5) such that $\mathcal{E}%
_{s,t}^{\prime }[X]=\mathcal{E}_{s,t}[X]$, a.s., for each deterministic
times $0\leq s\leq t\leq T$ and for each $X\in L^2(\mathcal{F}_t)$, then we
have $\mathcal{E}_{\sigma ,\tau }^{\prime }[X]=\mathcal{E}_{\sigma ,\tau
}[X] $, a.s., for each $\sigma $, $\tau \in \mathcal{S}_T$, $\sigma \leq
\tau $, and for each $X\in L^2(\mathcal{F}_\tau )$.
\end{remark}

We will give the proof of this Theorem in the last section.

To prove Theorem \ref{m6.1}, we need to introduce a sequence of BSDEs of the
following form: for $n=1,2,\cdots $,
\begin{equation}
y_t^n=\mathcal{E}_{t,T}[Y_T;n\int_0^{\cdot }(Y_s-y_s^n)ds].  \label{e6.2}
\end{equation}
The solution $y^n\in S_{\mathcal{F}}^2(0,T)$ has the following interesting
property.

\begin{lemma}
\label{m6.2}For each $n=1,2,\cdots $, we have $y_t^n\leq Y_t$, $dt\times dP$%
--a.e..
\end{lemma}

\smallskip\noindent\textbf{Proof. }For each fixed $n$, and any $\delta >0$, we define
\[
\sigma ^\delta :=\inf \{t\geq 0;y_t^n\geq Y_t+\delta \}\wedge T.
\]
If for all $\delta >0$ we always have $P\{\sigma ^\delta =T\}=1$, then we
have our conclusion. Otherwise, there exists a $\delta >0$ such that
\[
P(\sigma ^\delta >T)>0.
\]
We then define
\[
\tau :=\inf \{t\geq \sigma ^\delta ;y_t^n\leq Y_t\}.
\]
Since $y_T^n=Y_T$, we have $P(\tau \leq T)=1$. Since $\mathcal{E}%
_{s,t}[\cdot ;n\int_0^{\cdot }(Y_s-y_s^n)ds]$, $0\leq s\leq t\leq T$, is an $%
\mathcal{F}_t$--consistent evaluation satisfying (A1)--(A5), by Theorem \ref
{m8.13}, it can be uniquely extended to $\mathcal{E}_{\rho ,\sigma }[\cdot
;n\int_0^{\cdot }(Y_s-y_s^n)ds]$, $\rho \leq \sigma $, $\rho $, $\sigma \in
\mathcal{S}_T$. This with the fact that $(y_t^n)_{t\in [0,T]}$ is an
martingale under $\mathcal{E}_{s,t}[\cdot ;n\int_0^{\cdot }(Y_s-y_s^n)ds]$,
it follows from the Optional Stopping Theorem \ref{m8.14} that
\[
y_\rho ^n=\mathcal{E}_{\rho ,\sigma }[y_\sigma ^n;n\int_0^{\cdot
}(Y_s-y_s^n)ds],\;\hbox{and }\mathcal{E}_{\rho ,\sigma }[Y_\sigma
]\leq Y_\rho \hbox{, a.s.}
\]
We then have
\begin{eqnarray*}
y_{\sigma ^\delta }^n &=&\mathcal{E}_{\sigma ^\delta ,\tau }[y_\tau
^n;\int_0^{\cdot }n(Y_s-y_s^n)ds] \\
\ &\leq &\mathcal{E}_{\sigma ^\delta ,\tau }[y_\tau ^n]+\mathcal{E}_{\sigma
^\delta ,\tau }^{g_\mu }[0;\int_0^{\cdot }n(Y_s-y_s^n)ds].
\end{eqnarray*}
The inequality is from (\ref{e8.34}). But since $n(Y_s-y_s^n)\leq 0$, on $%
[\sigma ^\delta ,\tau ]$, by the definition of $\mathcal{E}^g[\cdot ;K]$
(see (\ref{e2.3})) and the Comparison Theorem of BSDE we have $\mathcal{E}%
_{\sigma ^\delta ,\tau }^{g_\mu }[0;\int_0^{\cdot }n(Y_s-y_s^n)ds]\leq 0$.
Thus
\[
y_{\sigma ^\delta }^n\leq \mathcal{E}_{\sigma ^\delta ,\tau }[y_\tau ^n]=%
\mathcal{E}_{\sigma ^\delta ,\tau }[Y_\tau ]\leq Y_{\sigma ^\delta }.
\]
But this contradicts with $y_t^n\geq Y_t+\delta $ on $\{\sigma ^\delta <T\}$
and $P(\sigma ^\delta <T)>0$. The proof is complete. $\Box $\medskip\

By Lemma \ref{m6.2}, we can rewrite BSDE (\ref{e6.2}) as:
\begin{equation}
y_t^n=\mathcal{E}_{t,T}[Y_T;n\int_0^{\cdot }(Y_s-y_s^n)^{+}ds].  \label{e6.3}
\end{equation}
By comparison theorem (Corollary \ref{m5.5}), we have
\begin{equation}
y_t^1\leq y_t^2\leq \cdots \leq Y_t.  \label{e6.4}
\end{equation}
Since $y^1$ and $Y$ are both in $S_{\mathcal{F}}^2(0,T)$ it follows from $%
|y_t^n|\leq |y_t^1|+|Y_t|$ that, there exists a constant $C>0$ which is
independent of $n$ such that
\begin{equation}
E[\sup_{t\in [0,T]}|y_t^n|^2]\leq C.  \label{e6.5}
\end{equation}
We define
\begin{equation}
A_t^n:=n\int_0^t(Y_s-y_s^n)^{+}ds=n\int_0^t(Y_s-y_s^n)ds,\;t\in [0,T].
\label{e6.6a}
\end{equation}
It follows from Proposition \ref{m4.3} that $y^n$ has the expression
\begin{equation}
y_t^n=Y_T+A_T^n-A_t^n+\int_t^Tg_s^nds-\int_t^Tz_s^ndB_s,\;\;t\in [0,T],
\label{e6.6}
\end{equation}
where $(g^n,z^n)\in L_{\mathcal{F}}^2(0,T;R\times R^d)$ satisfies, for each $%
m,n=1,2,\cdots ,$%
\begin{eqnarray}
\;|g_t^n| &\leq &\mu |y_t^n|+\mu |z_t^n|,  \label{e6.7} \\
|g_t^n-g_t^m| &\leq &\mu |y_t^n-y_t^m|+\mu |z_t^n-z_t^m|,\;  \label{e6.8} \\
\forall t &\in &[0,T],\;dt\times dP,\hbox{ a.e.}  \nonumber
\end{eqnarray}
We have the following estimates.

\begin{lemma}
\label{m6.3}There exists a constant $C>0$ which is independent of $n$ such
that
\begin{equation}
E[\int_0^T|z_t^n|^2dt\leq C,\;\;E[|A_T^n|^2]\leq C.  \label{e6.9}
\end{equation}
\end{lemma}

\smallskip\noindent\textbf{Proof. } From (\ref{e6.6}) we have
\begin{eqnarray*}
A_T^n &=&y^n(0)-y_T^n-\int_0^Tg_s^nds+\int_0^Tz_s^ndB_s \\
\ &\le &|y^n(0)|+|y_T^n|+\int_0^T\mu (|y_s^n|+|z_s^n|)ds+|\int_0^Tz_s^ndB_s|.
\end{eqnarray*}
With (\ref{e6.5}) it follows that there are two constants $c_1$ and $c_2$,
independent of $n$, such that
\begin{equation}
E|A_T^n|^2\le c_1+c_2E\int_0^T|z_s^n|^2ds.  \label{e6.10}
\end{equation}
On the other hand, It\^o's formula applied to $|y^n(\cdot )|^2$ gives:

\begin{eqnarray*}
E[|y^n(0)|^2] &=&E|Y_T|^2+E\int_0^T[2y_s^n\cdot g_s^n-|z_s^n|^2]ds \\
&&\ +2E\int_0^Ty_s^ndA_s^n \\
\ &\le &E|Y_T|^2+E\int_0^T[2\mu |y_s^n|(|z_s^n|+|y_s^n|)-|z_s^n|^2]ds\  \\
&&\ +E[2A_T^n\sup_{0\le s\le T}|y_s^n|].
\end{eqnarray*}
Thus, by $2\mu |y^n||z^n|\leq 2\mu ^2|y^n|^2+1/2|z^n|^2$ and $2A_T^n\sup
|y_s^n|\leq \frac 1{4c_2}|A_T^n|^2+4c_2\sup |y_s^n|^2$,
\begin{eqnarray*}
E\int_0^T|z_s^n|^2ds\ &\le &2E|Y_T|^2+E\int_0^T(4\mu ^2+4\mu +1)|y_s^n|^2ds
\\
&&\ +8c_2E[\sup_{0\le s\le T}|y_s^n|^2]+{\frac 1{2c_2}}[E|A_T^n|^2] \\
\ &\leq &2E|Y_T|^2+E\int_0^T(4\mu ^2+4\mu +1)|y_s^n|^2ds \\
&&\ +8c_2[E\sup_{0\le s\le T}|y_s^n|^2]+{\frac{c_1}{2c_2}}+{\frac 12}%
E\int_0^T|z_s^n|^2ds.
\end{eqnarray*}
The last inequality is due to (\ref{e6.10}). Then the first estimate of (\ref
{e6.9}) yields immediately from (\ref{e6.5}). From which and (\ref{e6.10})
we obtain the second one. The proof is complete. $\Box $\medskip\

We rewrite (\ref{e6.6}) in the following forward version:
\begin{equation}
y_t^n=y_0^n-A_t^n-\int_0^tg_s^nds+\int_0^tz_s^ndB_s,\;\;t\in [0,T],
\label{e6.11}
\end{equation}
in order to apply the following monotonicity limit theorem (see
\cite[Peng 1999]{Peng1999}, Theorem 2.1).

\begin{proposition}
\label{m6.4}Let $\{A^n\}_{n=1}^\infty $ be a sequence of increasing
processes in $S_{\mathcal{F}}^2(0,T)$ with $A_0^n=0$, and let $%
\{(g^n,z^n)\}_{n=1}^\infty $ be uniformly bounded in $L_{\mathcal{F}}^2(0,T)$%
:
\begin{equation}
E\int_0^T[|g_s^n|^2+|z_s^n|^2]\leq C.  \label{e6.12}
\end{equation}
If $\{(y_t^n)_{t\in [0,T]}\}_{n=1}^\infty $ increasingly converges to $%
(y_t)_{t\in [0,T]}$ with $E[\sup_{t\in [0,T]}|y_t|^2]<\infty $, then this
limit process has the following form
\begin{equation}
y_t=y_0-A_t-\int_0^tg_sds+\int_0^tz_sdB_s,\;t\in [0,T],  \label{e6.13}
\end{equation}
where $A\in D_{\mathcal{F}}^2(0,T)$ is an increasing process such that $%
A_0=0 $, $(g,z)\in L_{\mathcal{F}}^2(0,T;R\times R^d)$, and
\begin{equation}
\begin{array}{rl}
A_t^n\rightharpoonup & A_t\;\hbox{weakly in }L^2(\mathcal{F}_T),\hbox{ }%
\forall t\in [0,T], \\
(g^n,z^n)\rightharpoonup & (g,z)\;\;\hbox{weakly in }L_{\mathcal{F}%
}^2(0,T;R\times R^n), \\
\lim_{n\rightarrow \infty } &
E\int_0^T|z_t^n-z_t|^pdt=0,\;\hbox{for each fiexd }p\in
[1,2)\hbox{.}
\end{array}
\label{e6.14}
\end{equation}
Moreover, if $y$ is continuous, i.e., $y\in S_{\mathcal{F}}^2(0,T)$, then we
have
\[
\lim_{n\rightarrow \infty }E\int_0^T|z_t^n-z_t|^2dt=0.
\]
\end{proposition}

\medskip\ We now can proceed to give

\noindent\textbf{Proof of Theorem \ref{m6.1}. }Since $A^n$ defined by (\ref
{e6.6a}) is bounded by the second estimate of (\ref{e6.9}), it follows that $%
y_t^n\nearrow Y_t$ $dt\times dP$--a.e.. On the other hand, by (\ref{e6.5}), (%
\ref{e6.7}) and the first estimate of (\ref{e6.9}), $\{(g^n,z^n)\}_{n=1}^%
\infty $ is also uniformly bounded in $L_{\mathcal{F}}^2(0,T)$. We then can
apply Proposition \ref{m6.4} to derive that
\begin{equation}
Y_t=y_0-A_t-\int_0^tg_sds+\int_0^tz_sdB_s,\;t\in [0,T].  \label{e6.15}
\end{equation}
Since $Y\in S_{\mathcal{F}}^2(0,T)$, then $\lim_{n\rightarrow \infty
}E\int_0^T|z_t^n-z_t|^2dt=0$. We also have $A\in S_{\mathcal{F}}^2(0,T)$.
But observe that $y^n$ also converges to $Y$ in $L_{\mathcal{F}}^2(0,T)$.
Thus, in considering (\ref{e6.8}), $\{g^n\}_{n=1}^\infty $ also converges
strongly to $g$ in $L_{\mathcal{F}}^2(0,T)$. It follows that $%
A_t^n\rightarrow A_t$ in $L^2(\mathcal{F}_T)$ for each $t$, and thus $A^n$
also converges strongly in $L_{\mathcal{F}}^2(0,T)$. By Corollary~\ref{m3.9}
\begin{eqnarray}
&&E[|\mathcal{E}_{t,T}[Y_T;A_{\cdot }]-\mathcal{E}_{t,T}[Y_T;A_{\cdot
}^n]|^2]  \label{e6.16} \\
&\leq &CE[(A_T-A_T^n)^2]+CE\int_0^T(A_t-A_t^n)^2ds\rightarrow 0.  \nonumber
\end{eqnarray}
We then have
\begin{equation}
Y_t=\lim_{n\rightarrow \infty }y_t^n=\lim_{n\rightarrow \infty }\mathcal{E}%
_{t,T}[Y_T;A_{\cdot }^n]=\mathcal{E}_{t,T}[Y_T;A_{\cdot }].  \label{e6.17}
\end{equation}
The proof is complete. $\Box $\medskip\

\section{Proof of Theorem \ref{m7.1}\label{ss8}}

For each fixed $(t,y,z)\in [0,T]\times R\times R^d$, we consider the
solution $Y^{t,y,z}\in S_{\mathcal{F}}^2(0,T)$ of a It\^o's equation on $%
[t,T]$, and a BSDE on $[0,t]$:
\begin{eqnarray}
dY_s^{t,y,z} &=&-\mu (|Y_s^{t,y,z}|+|z|)ds+zdB_s,\;\;s\in (t,T],
\label{e7.2} \\
Y_t^{t,y,z} &=&y.  \label{e7.4}
\end{eqnarray}
It is easy to check that $Y^{t,y,z}$ is an $\mathcal{E}^{g_\mu }[\cdot ]$%
--martingale. Thus, by (\ref{e3.1}), it is also an $\mathcal{E}[\cdot ]$%
--supermartingale. By Decomposition Theorem \ref{m6.1}, there exists an
increasing process $A^{t,y,z}\in S_{\mathcal{F}}^2(0,T)$ with $A_0^{t,y,z}=0$%
, such that
\[
Y_s^{t,y,z}=\mathcal{E}_{s,T}[Y_T^{t,y,z};A_{\cdot }^{t,y,z}]. \label{e7.5}
\]
By Proposition \ref{m4.3} and Corollary \ref{m4.4}, there exists $($ $%
g^{t,y,z},Z^{t,y,z})\in L_{\mathcal{F}}^2(0,T)$ such that
\begin{equation}
-dY_s^{t,y,z}=dA_s^{t,y,z}+g_s^{t,y,z}ds-Z_s^{t,y,z}dB_s,\;s\in [t,T],
\label{e7.6}
\end{equation}
and such that, for each different $(t,y,z)$, $(t^{\prime },y^{\prime
},z^{\prime })$ $\in [0,T]\times R\times R^d$%
\begin{equation}
|g_s^{t,y,z}-g_s^{t^{\prime },y^{\prime },z^{\prime }}|\leq \mu
|Y_s^{t,y,z}-Y_s^{t^{\prime },y^{\prime },z^{\prime }}|+\mu
|Z_s^{t,y,z}-Z_s^{t^{\prime },y^{\prime },z^{\prime }}|,\;s\in [t\vee
t^{\prime },T],  \label{e7.7}
\end{equation}
and
\begin{equation}
|g_s^{t,y,z}|\leq \mu |Y_s^{t,y,z}|+\mu |Z_s^{t,y,z}|,\;s\in
[t,T],\;ds\times dP\hbox{--a.e.}  \label{e7.7a}
\end{equation}

Now for each $X\in L^2(\mathcal{F}_{t^{\prime }})$, we set
\begin{equation}
\bar Y_s^{t^{\prime },X}:=\mathcal{E}_{s,t^{\prime }}[X]=\mathcal{E}%
_{s,t^{\prime }}[X;0].  \label{e7.8a}
\end{equation}
We use once more Proposition \ref{m4.3} and Corollary \ref{m4.4}: there
exists $($ $\bar g^{t^{\prime },X},\bar Z^{t^{\prime },X})\in L_{\mathcal{F}%
}^2(0,t^{\prime })$ such that, for $s\in [0,t^{\prime }]$,
\begin{equation}
-d\bar Y_s^{t^{\prime },X}=\bar g_s^{t^{\prime },X}ds-\bar Z_s^{t^{\prime
},X}dB_s,\;\bar Y_{t^{\prime }}=X,  \label{e7.8}
\end{equation}
such that
\begin{equation}
|g_s^{t,y,z}-\bar g_s^{t^{\prime },X}|\leq \mu |Y_s^{t,y,z}-\bar
Y_s^{t^{\prime },X}|+\mu |Z_s^{t,y,z}-\bar Z_s^{t^{\prime
},X}|,\;s\in [t,t^{\prime }],\;ds\times dP\hbox{--a.e.}
\label{e7.9}
\end{equation}
and, for $X$, $X\in L^2(\mathcal{F}_{t^{\prime }})$,
\[
|\bar g_s^{t^{\prime },X}-\bar g_s^{t^{\prime },X^{\prime }}|\leq
\mu |\bar Y_s^{t^{\prime },X}-\bar Y_s^{t^{\prime },X^{\prime
}}|+\mu |\bar Z_s^{t^{\prime },X}-\bar Z_s^{t^{\prime },X^{\prime
}}|,\;s\in [0,t^{\prime }],\;ds\times dP\hbox{--a.e..}
\label{e7.10}
\]
On the other hand, comparing to (\ref{e7.2}) and (\ref{e7.6}), we have
\[
Z_s^{t,y,z}\equiv 1_{[t,T]}(s)z. \label{e7.11}
\]
Thus (\ref{e7.7}), (\ref{e7.7a}) and (\ref{e7.9}) become, respectively,
\begin{equation}
|g_s^{t,y,z}-g_s^{t^{\prime },y^{\prime },z^{\prime }}|\leq \mu
|Y_s^{t,y,z}-Y_s^{t^{\prime },y^{\prime },z^{\prime }}|+\mu
|z-z^{\prime }|,\;s\in [t\vee t^{\prime },T],\;ds\times
dP\hbox{--a.e.,}  \label{e7.12}
\end{equation}
\begin{equation}
|g_s^{t,y,z}|\leq \mu |Y_s^{t,y,z}|+\mu |z|,\;  \label{e7.12a}
\end{equation}
and
\begin{equation}
|g_s^{t,y,z}-\bar g_s^{t^{\prime },X}|\leq \mu |Y_s^{t,y,z}-\bar
Y_s^{t^{\prime },X}|+\mu |z-\bar Z_s^{t^{\prime },X}|,\;s\in
[t,t^{\prime }],\;ds\times dP\hbox{--a.e.}  \label{e7.13}
\end{equation}

Now, for each $n=1,2,3,\cdots $, we set $t_i^n=i2^{-n}T$, $i=0,1,2,\cdots
,2^n$, and define
\begin{equation}
g^n(s,y,z):=\sum_{i=0}^{2^n-1}g_s^{t_i^n,y,z}1_{[t_i^n,t_{i+1}^n)}(s),\;s\in
[0,T].  \label{e7.14}
\end{equation}
It is clear that $g^n$ is an $\mathcal{F}_t$--adapted process. We also have

\begin{lemma}
\label{m7.2}For each fixed $(y,z)\in R\times R^d$, $\{g^n(\cdot
,y,z)\}_{n=1}^\infty $ is a Cauchy sequence in $L_{\mathcal{F}}^2(0,T)$.
\end{lemma}

To prove this lemma, we need the following classical result of It\^o's SDE.
The proof is classic.

\begin{lemma}
\label{m7.2a}We have the following estimate: there exist a constant
depending only on $\mu $ and $T$ such that, for each $(t,y,z)\in [0,T]\times
R\times R^d$,
\begin{equation}
E[|Y_s^{t,y,z}-y|^2]\leq C(|y|^2+|z|^2+1)(s-t),\;\forall s\in [t,T].
\label{e7.15}
\end{equation}
\end{lemma}

We can give

\medskip\

\noindent \textbf{Proof of Lemma \ref{m7.2}.} Let $0<m<n$ be two integers.
For each $s\in [0,T)$, there are some integers $i\leq 2^m-1$ and $j\leq
2^n-1 $ with $t_i^m\leq t_j^m$, such that $s\in [t_i^m,t_{i+1}^m)$, $s\in
[t_j^n,t_{j+1}^n)$. We have, by (\ref{e7.12})
\begin{eqnarray*}
|g^m(s,y,z)-g^n(s,y,z)| &=&|g_s^{t_i^m,y,z}-g_s^{t_j^n,y,z}| \\
&\leq &\mu |Y_s^{t_i^m,y,z}-Y_s^{t_j^n,y,z}| \\
&\leq &\mu |Y_s^{t_i^m,y,z}-y|+\mu |Y_s^{t_j^n,y,z}-y|.
\end{eqnarray*}
By (\ref{e7.15}),
\begin{eqnarray*}
E[|g^m(s,y,z)-g^n(s,y,z)|^2] &\leq &2\mu
^2E[|Y_s^{t_i^m,y,z}-y|^2+|Y_s^{t_j^n,y,z}-y|^2] \\
&\leq &2\mu ^2E[|Y_s^{t_i^m,y,z}-y|^2+|Y_s^{t_j^n,y,z}-y|^2] \\
&\leq &2\mu ^2C(|y|^2+|z|^2+1)(2^{-m}+2^{-n})T.
\end{eqnarray*}
Thus
\[
\sup_{s\in [0,T)}E[|g^m(s,y,z)-g^n(s,y,z)|^2]\leq 2\mu
^2C(|y|^2+|z|^2+1)(2^{-m}+2^{-n})T. \label{e7.16}
\]
Thus $\{g^n(\cdot ,y,z)\}_{n=1}^\infty $ is a Cauchy sequence in $L_{%
\mathcal{F}}^2(0,T)$. $\Box $\medskip\

\begin{definition}
\label{m7.3}For each $(y,z)\in R\times R^d$, we denote $g(\cdot ,y,z)\in L_{%
\mathcal{F}}^2(0,T)$, the Cauchy limit of $\{g^n(\cdot ,y,z)\}_{n=1}^\infty $
in $L_{\mathcal{F}}^2(0,T)$.
\end{definition}

We will prove that this function is just what we are looking for in Theorem
\ref{m7.1}. We still need to investigate some important properties of $g$.
We have the following estimates for the function $g$.

\begin{lemma}
\label{m7.4}The limit $g:\Omega \times [0,T]\times R\times R^d\rightarrow
R^d $ satisfies the following properties:
\begin{equation}
\left\{
\begin{array}{rrl}
\hbox{(i)}\, &  & g(\cdot ,y,z)\in L_{\mathcal{F}}^2(0,T)\hbox{, for each }%
(y,z)\in R\times R^d; \\
\hbox{(ii)}\, &  & |g(s,y,z)-g(s,y^{\prime },z^{\prime })|\leq \mu
(|y-y^{\prime }|+|z-z^{\prime }|),\;\forall y,y^{\prime }\in
R,\;z,z^{\prime
}\in R^d; \\
\hbox{(iii)}\, &  & g(s,0,0)\equiv 0; \\
\hbox{(iv)}\, &  & |g(s,y,z)-\bar g^{t,X}|\leq \mu |y-\bar
Y_s^{t,X}|+\mu |z-\bar Z_s^{t,X}|,\;s\in [0,t].
\end{array}
\right.  \label{e7.17}
\end{equation}
where, in (iv), $t\in [0,T]$ and $X\in L^2(\mathcal{F}_t)$ are arbitrarily
given. $(\bar Y^{t,X},\bar Z^{t,X})$ is the process defined in (\ref{e7.8a})
and (\ref{e7.8}).
\end{lemma}

\smallskip\noindent\textbf{Proof. }(i) is clear. To prove (ii), we choose $%
t_i^n=i2^{-n}T$, $i=0,1,2,\cdots ,2^n$ as in (\ref{e7.14}). For each $s\in
[0,T)$. We have, once more by (\ref{e7.12}),
\begin{eqnarray}
|g^n(s,y,z)-g^n(s,y^{\prime },z^{\prime })|
&=&%
\sum_{j=0}^{2^n-1}1_{[t_j^n,t_{j+1}^n)}(s)|g_s^{t_j^n,y,z}-g_s^{t_j^n,y,z}|
\label{e7.17a} \\
&\leq &\mu
\sum_{j=0}^{2^n-1}1_{[t_j^n,t_{j+1}^n)}(s)(|Y_s^{t_j^n,y,z}-Y_s^{t_j^n,y,z}|+|z-z^{\prime }|)
\nonumber \\
&\leq &\mu
\sum_{j=0}^{2^n-1}1_{[t_j^n,t_{j+1}^n)}(s)(|Y_s^{t_j^n,y,z}-y|+|Y_s^{t_j^n,y,z}-y^{\prime }|)
\nonumber \\
&&+\mu (|y-y^{\prime }|+|z-z^{\prime }|)  \nonumber
\end{eqnarray}
For the first term $I^n(s)$ of the left hand, we have, by (\ref{e7.15}),
\begin{eqnarray*}
E[|I^n(s)|^2] &\leq &2\mu
^2%
\sum_{i=0}^{2^n-1}1_{[t_j^n,t_{j+1}^n)}(s)E[|Y_s^{t_j^n,y,z}-y|^2+|Y_s^{t_j^n,y^{\prime },z^{\prime }}-y^{\prime }|^2]
\\
&\leq &2\mu
^2\sum_{i=0}^{2^n-1}1_{[t_j^n,t_{j+1}^n)}(s)C(|y|^2+|z|^2+|y^{\prime
}|^2+|z^{\prime }|^2+2)2^{-n}T.
\end{eqnarray*}
Thus $I^n(\cdot )\rightarrow 0$ in $L_{\mathcal{F}}^2(0,T)$ as $n\rightarrow
\infty $. (ii) is obtained by passing to the limit in both sides of (\ref
{e7.17a}). (iii) is proved similarly by using (\ref{e7.12a}) and (\ref{e7.15}%
).

To prove (iv), We apply (\ref{e7.13}),

\begin{eqnarray*}
|g^n(s,y,z)-\bar g_s^{t,X}|
&=&\sum_{i=0}^{2^n-1}1_{[t_j^n,t_{j+1}^n)}(s)|g_s^{t_j^n,y,z}-\bar g_s^{t,X}|
\\
\ &\leq &\sum_{i=0}^{2^n-1}1_{[t_j^n,t_{j+1}^n)}(s)[\mu
|Y_s^{t_j^n,y,z}-\bar Y_s^{t,X}|+\mu |z-\bar Z_s^{t,X}|] \\
\ &\leq &\mu
\sum_{i=0}^{2^n-1}1_{[t_j^n,t_{j+1}^n)}(s)|Y_s^{t_j^n,y,z}-y|+\mu |y-\bar
Y_s^{t,X}|+\mu |z-\bar Z_s^{t,X}|.
\end{eqnarray*}
Then we pass to the limit on both sides. $\Box $\medskip\

Finally, We give

\medskip\

\noindent \textbf{Proof of Theorem \ref{m7.1}.} For each fixed $t\in [0,T]$
and $X\in L^2(\mathcal{F}_t)$, we consider $\bar Y_s^{t,X}:=\mathcal{E}%
_{s,t}[X]$, $s\in [0,T]$. By Proposition \ref{m4.3} and Corollary \ref{m4.4}%
, we can write
\[
\bar Y_s^{t,X}=X+\int_s^t\bar g_r^{t,X}dr-\int_s^t\bar Z_r^{t,X}dB_r,\;s\in
[0,t].
\]
On the other hand, let $(Y^{t,X},Z^{t,X})$ be the solution of the following
BSDE
\[
Y_s=X+\int_s^tg(r,Y_r,Z_r)dr-\int_s^tZ_rdB_r,\;s\in [0,t].
\]
By Lemma \ref{m7.4}--(i) and (ii), this BSDE is well--posed. We then apply
It\^o's formula to $|\bar Y^{t,X}-Y|^2$ in the interval $[0,t]$, take
expectation and then apply (iv) of Lemma \ref{m7.4}. Exactly as the
classical proof of the uniqueness of BSDE, we have \newpage
\begin{eqnarray*}
E|\bar Y_s^{t,X}-Y_s|^2 &+&E\int_s^t|\bar Z_r^{t,X}-Z_r|^2dEr \\
&=&2E\int_s^t(\bar Y_r^{t,X}-Y_r)(\bar g_r^{t,X}-g(r,Y_r,Z_r))dr \\
&\leq &2E\int_s^t(|\bar Y_r^{t,X}-Y_r|\cdot |\bar g_r^{t,X}-g(r,Y_r,Z_r)|)dr
\\
&\leq &2E\int_s^t|\bar Y_r^{t,X}-Y_r|\cdot \mu (|\bar Y_r^{t,X}-Y_r|+|\bar
Z_r^{t,X}-Z_r|)dr \\
&\leq &E\int_s^t2(\mu +\mu ^2)|\bar Y_r^{t,X}-Y_r|^2+\frac 12|\bar
Z_r^{t,X}-Z_r|^2)dr.
\end{eqnarray*}
It then follows by using Gronwall's inequality that $\bar Y_s^{t,X}\equiv
Y_s $, $s\in [0,t]$. Recall that $Y_s=\mathcal{E}_{s,t}^g[X]$. We thus have
the desired result. The proof is complete. $\Box $\medskip\

\section{Proof of Theorem \ref{m8.13} and optional stopping theorem for $%
\mathcal{E}_{\sigma ,\tau }[\cdot ]$\label{ss9}}

We now consider the situation where the time indices $s$ and $t$ in $%
\mathcal{E}_{s,t}[\cdot ]$ is replaced by stopping times $\sigma $, $\tau
\in \mathcal{S}_T$, $\sigma \leq \tau \leq T$.\ We will extend $\mathcal{E}%
_{s,t}[\cdot ]$ to $\mathcal{E}_{\sigma ,\tau }[\cdot ]$ and prove Theorem
\ref{m8.13}. We will also obtain a generalized version of the optional
stopping theorem for $\mathcal{E}$--super and $\mathcal{E}$%
--sub--martingale. We note that it is not at all a trivial task to define $%
\mathcal{E}_{\sigma ,\tau }[\cdot ]$, especially for the second parameter $%
\tau $. We will first consider the situation of discrete--valued stopping
times, i.e., $\sigma $, $\tau \in \mathcal{S}_T^0$. Then we will pass to the
limit to treat the $\mathcal{S}_T$ case.

\subsection{Simple case: $\mathcal{E}_{\sigma ,\tau }^g[\cdot ]$ with $%
\sigma $, $\tau \in \mathcal{S}_T^0$ \label{ss9.1}}

$\mathcal{E}_{\sigma ,\tau }^g[\cdot ]$, $\sigma $, $\tau \in \mathcal{S}_T$
will provide us a concrete example. In this situation. For a given $X\in
\mathcal{F}_\tau $, we can directly solve the BSDE
\begin{equation}
Y_s=X+\int_s^\tau g(r,Y_r,Z_r)dr-\int_s^\tau Z_rdB_r,\;s\in [0,\tau ],
\label{e8.1}
\end{equation}
or equivalently, on $s\in [0,T],$%
\begin{equation}
Y_s=X+\int_s^T1_{[0,\tau ]}(r)g(r,Y_r,Z_r)dr-\int_s^\tau 1_{[0,\tau
]}(r)Z_rdB_r  \label{e8.3}
\end{equation}
and then define
\begin{equation}
\mathcal{E}_{\sigma ,\tau }^g[X]:=Y_\sigma .  \label{e8.3a}
\end{equation}
It is clear that, when $\sigma =s$ and $\tau =t$ for deterministic time
parameters $s\leq t$, then $\mathcal{E}_{\sigma ,\tau }^g[\cdot ]=\mathcal{E}%
_{s,t}^g[\cdot ]$. We have

\begin{proposition}
\label{m8.0}The system of operators
\[
\mathcal{E}_{\sigma ,\tau }^g[\cdot ]:L^2(\mathcal{F}_\tau )\rightarrow L^2(%
\mathcal{F}_\sigma ),\;\sigma \leq \tau ,\;\sigma ,\tau \in \mathcal{S}_T,
\]
is an $\mathcal{F}_t$--consistent nonlinear evaluation, i.e., it satisfies
(A1)--(A5) in the following sense: for each $X$, $X^{\prime }\in L^2(%
\mathcal{F}_\tau )$, \\\textbf{(a1)} $\mathcal{E}_{\sigma ,\tau }^g[X]\geq
\mathcal{E}_{\sigma ,\tau }^g[X^{\prime }]$, a.s., if $X\geq X^{\prime }$,
a.s. \\\textbf{(a2)} $\mathcal{E}_{\tau ,\tau }^g[X]=X$; \\\textbf{(a3)} $%
\mathcal{E}_{\rho ,\sigma }^g[\mathcal{E}_{\sigma ,\tau }^g[X]]=\mathcal{E}%
_{\rho ,\tau }^g[X]$, $\forall 0\leq \rho \leq \sigma \leq \tau $; \\\textbf{%
(a4')} $1_A\mathcal{E}_{\sigma ,\tau }^g[X]=\mathcal{E}_{\sigma ,\tau
}^g[1_AX]$,\ $\forall A\in \mathcal{F}_\tau $; \\\textbf{(a5)} for each $%
0\leq \sigma \leq \tau \leq T$,
\begin{equation}
\mathcal{E}_{\sigma ,\tau }^g[X]-\mathcal{E}_{\sigma ,\tau }^g[X^{\prime
}]\leq \mathcal{E}_{\sigma ,\tau }^{g_\mu }[X-X^{\prime }],\;\;\forall
X,X^{\prime }\in L^2(\mathcal{F}_\tau ).  \label{e8.a5}
\end{equation}
\end{proposition}

\smallskip\noindent\textbf{Proof. }The proof is analogous to the situation of $%
\mathcal{E}_{s,t}^g[\cdot ]$. We omit it. $\Box $\medskip\

To define $\mathcal{E}_{\sigma ,\tau }[\cdot ]$, we first consider the
situation $\mathcal{E}_{s\wedge \tau ,t\wedge \tau }^g[\cdot ]$, where $%
0\leq s\leq t\leq T$ and $\tau \in \mathcal{S}_T^0$. We often let $\tau $ be
characterized by
\begin{equation}
\;\cup _{i=1}^n\{\tau =t_i\}=\Omega ,\;0=t_0\leq t_1<\cdots <t_n\leq
t_{n+1}=T.  \label{e8.4}
\end{equation}
We consider a more special case where
\begin{equation}
t_i\leq s<t\leq t_{i+1},\;\hbox{for some\ }i\in \{1,2,\cdots
,n\}\hbox{.} \label{e8.5}
\end{equation}

\begin{lemma}
\label{m8.1}In the situation (\ref{e8.4}) with $s$ and $t$ limited to (\ref
{e8.5}), we have, for each $X\in \mathcal{F}_{t\wedge \tau }$%
\begin{equation}
\left\{
\begin{array}{rl}
\hbox{(i)} & \mathcal{E}_{t\wedge \tau ,t\wedge \tau }^g[X]=X; \\
\hbox{(ii)} & \mathcal{E}_{s\wedge \tau ,t\wedge \tau
}^g[X]=1_{\{t\wedge \tau \leq s\}}X+1_{\{t\wedge \tau
=t\}}\mathcal{E}_{s,t}^g[X].
\end{array}
\right.  \label{e8.6}
\end{equation}
\end{lemma}

\smallskip\noindent\textbf{Proof. }(i) is easy. To prove (ii), we first observe that
\begin{equation}
\{t\wedge \tau \leq s\}^C=\{t\wedge \tau =t\}  \label{e8.7}
\end{equation}
and $\{t\wedge \tau \leq s\}=\{t\wedge \tau \leq t_i\}$. Thus $1_{\{t\wedge
\tau \leq s\}}X\in \mathcal{F}_{t_i}$. We also have $1_{\{t\wedge \tau
=t\}}X\in \mathcal{F}_t$. We now solve $Y_{s\wedge \tau }=\mathcal{E}%
_{s\wedge \tau ,t\wedge \tau }^g[X]$ by, as in (\ref{e8.3}),
\begin{equation}
Y_{s\wedge \tau }=\ X+\int_s^T1_{[0,t\wedge \tau
]}(r)g(r,Y_r,Z_r)dr-\int_s^T1_{[0,t\wedge \tau ]}(r)Z_rdB_r.  \label{e8.8}
\end{equation}
Since $1_{[0,t\wedge \tau ]}=1_{\{t\wedge \tau \leq
t_i\}}1_{[0,t_i]}+1_{\{t\wedge \tau =t\}}1_{[0,t]}$. By respectively
multiplying $1_{\{t\wedge \tau \leq t_i\}}$ and $1_{\{t\wedge \tau =t\}}$ on
both sides of (\ref{e8.8}), we have, on $s\in [t_i,t)$,
\begin{equation}
Y_{s\wedge \tau }1_{\{t\wedge \tau \leq t_i\}}=\ X1_{\{t\wedge \tau \leq
t_i\}},  \label{e8.9}
\end{equation}
and
\begin{eqnarray*}
Y_{s\wedge \tau }\ 1_{\{t\wedge \tau =t\}} &=&\ 1_{\{t\wedge \tau
=t\}}X+\int_s^T1_{[0,t]}(r)1_{\{t\wedge \tau =t\}}g(r,Y_r,Z_r)dr \\
&&-\int_s^T1_{[0,t]}1_{\{t\wedge \tau =t\}}(r)Z_rdB_r \\
&=&1_{\{t\wedge \tau =t\}}X+\int_s^tg(r,1_{\{t\wedge \tau
=t\}}Y_r,1_{\{t\wedge \tau =t\}}Z_r)dr-\int_s^t1_{[0,t_{i+1}]}Z_rdB_r.
\end{eqnarray*}
We observe that, the last relation solves a BSDE on $[t_i,t]$. Thus
\[
Y_{s\wedge \tau }\ 1_{\{t\wedge \tau =t\}}=1_{\{t\wedge \tau =t\}}\mathcal{E}%
_{s,t}^g[1_{\{t\wedge \tau =t\}}X]=1_{\{t\wedge \tau =t\}}\mathcal{E}%
_{s,t}^g[X].
\]
This with (\ref{e8.9}) and (\ref{e8.7}), we then have (ii). $\Box $\medskip\

\subsection{$\mathcal{E}_{\sigma ,\tau }[\cdot ]$ with $\sigma $, $\tau \in
S_T^0$\label{ss9.2}}

Let a stopping time $\tau \in \mathcal{S}_T^0$ be characterized by (\ref
{e8.4}). For each $i=0,1,\cdots ,n$, for each $t_i\leq s<t\leq t_{i+1}$, $%
X\in \mathcal{F}_{t\wedge \tau }$, we define
\begin{equation}
\left\{
\begin{array}{rl}
\hbox{(i)} & \mathcal{E}_{t,t}^i[X]=\mathcal{E}_{t\wedge \tau
,t\wedge \tau
}[X]:=X; \\
\hbox{(ii)} & \mathcal{E}_{s,t}^i[X]=\mathcal{E}_{s\wedge \tau
,t\wedge \tau
}[X]:=1_{\{t\wedge \tau \leq s\}}X+1_{\{t\wedge \tau =t\}}\mathcal{E}%
_{s,t}[X].\;
\end{array}
\ \right.  \label{e8.10}
\end{equation}

The reason that we set $\mathcal{E}_{\sigma ,\tau }[\cdot ]$ satisfying (ii)
is as follows

\begin{lemma}
\label{m8.a3}Let $\mathcal{E}_{\sigma ,\tau }^{\prime }[\cdot ]:L^2(\mathcal{%
F}_\tau )\rightarrow L^2(\mathcal{F}_\sigma )$, $\sigma $, $\tau \in
\mathcal{S}_T^0$, $\sigma \leq \tau $, be a system of operators satisfying
the following ($\mathcal{F}_\tau $--consistent) conditions: for each $X\in
\mathcal{F}_\tau $
\[
\begin{array}{cll}
\hbox{(a)} & \mathcal{E}_{\tau ,\tau }^{\prime }[X]=X; &  \\
\hbox{(b)} & \mathcal{E}_{\rho ,\sigma }^{\prime
}[\mathcal{E}_{\sigma ,\tau }^{\prime }[X]]=\mathcal{E}_{\rho
,\tau }^{\prime }[X], & \forall \rho \leq
\sigma \leq \tau ; \\
\hbox{(c)} & 1_A\mathcal{E}_{\sigma ,\tau }^{\prime }[X]=1_A\mathcal{E}%
_{\sigma ,\tau }^{\prime }[1_AX],\; & \forall A\in \mathcal{F}_\sigma ; \\
\hbox{(d)} & \mathcal{E}_{\sigma ,\tau }^{\prime
}[X]=\sum_{i=1}^m1_{\{\sigma =s_i\}}\mathcal{E}_{s_i\wedge \tau
,\tau }^{\prime }[X], & \hbox{ if }\sum_{i=1}^m\{\sigma
=s_i\}=\Omega .
\end{array}
\]
We assume that $\mathcal{E}^{\prime }$ coincides with $\mathcal{E}$ in the
sense that $\mathcal{E}_{s,t}^{\prime }[X]=\mathcal{E}_{s,t}[X]$, for all
(deterministic) $0\leq s\leq t\leq T$, and $X\in \mathcal{F}_T$. Then,
necessarily, $\mathcal{E}^{\prime }$ satisfies (i) and (ii) of (\ref{e8.10}).
\end{lemma}

\smallskip\noindent\textbf{Proof. }(i) comes directly from (a). We now prove (ii).
Since $\{s\wedge \tau =t\wedge \tau \}=\{t\wedge \tau \leq s\}$, we have, by
(d) and (a),
\begin{eqnarray}
1_{\{t\wedge \tau \leq s\}}\mathcal{E}_{s\wedge \tau ,t\wedge \tau }^{\prime
}[X] &=&1_{\{t\wedge \tau \leq s\}}1_{\{s\wedge \tau =t\wedge \tau \}}%
\mathcal{E}_{s\wedge \tau ,t\wedge \tau }^{\prime }[X]  \label{e8.10a} \\
\ &=&1_{\{t\wedge \tau \leq s\}}X.  \nonumber
\end{eqnarray}
By (d) we have
\begin{equation}
1_{\{t\wedge \tau =t\}}\mathcal{E}_{t\wedge \tau ,t}^{\prime
}[X]=1_{\{t\wedge \tau =t\}}\mathcal{E}_{t,t}^{\prime }[X]=1_{\{t\wedge \tau
=t\}}X.  \label{e8.10b}
\end{equation}
Since
\begin{eqnarray*}
\{t\wedge \tau &=&t\}\cap \{s\wedge \tau \leq r\} \\
&=&\left\{
\begin{array}{l}
\{t\wedge \tau =t\}\cap \Omega \in \mathcal{F}_s,\;\hbox{if }r\geq s; \\
\{t\wedge \tau =t\}\cap \{s\wedge \tau =s\}\cap \{s\wedge \tau
\leq r\}=\emptyset \hbox{, if }r<s\hbox{,}
\end{array}
\right.
\end{eqnarray*}
thus $\{t\wedge \tau =t\}\in \mathcal{F}_{s\wedge \tau }$. Thus, by (c) and (%
\ref{e8.10b}),
\begin{eqnarray*}
1_{\{t\wedge \tau =t\}}\mathcal{E}_{s\wedge \tau ,t\wedge \tau }^{\prime
}[X] &=&1_{\{t\wedge \tau =t\}}\mathcal{E}_{s\wedge \tau ,t\wedge \tau
}^{\prime }[1_{\{t\wedge \tau =t\}}X] \\
&=&1_{\{t\wedge \tau =t\}}\mathcal{E}_{s\wedge \tau ,t\wedge \tau }^{\prime
}[1_{\{t\wedge \tau =t\}}\mathcal{E}_{t,t}^{\prime }[X]] \\
&=&1_{\{t\wedge \tau =t\}}\mathcal{E}_{s\wedge \tau ,t\wedge \tau }^{\prime
}[1_{\{t\wedge \tau =t\}}\mathcal{E}_{t\wedge \tau ,t}^{\prime }[X]] \\
&=&1_{\{t\wedge \tau =t\}}\mathcal{E}_{s\wedge \tau ,t\wedge \tau }^{\prime
}[\mathcal{E}_{t\wedge \tau ,t}^{\prime }[X]].
\end{eqnarray*}
It follows from (b) that
\begin{eqnarray*}
1_{\{t\wedge \tau =t\}}\mathcal{E}_{s\wedge \tau ,t\wedge \tau }^{\prime
}[X] &=&1_{\{t\wedge \tau =t\}}1_{\{s\wedge \tau =s\}}\mathcal{E}_{s\wedge
\tau ,t}^{\prime }[X] \\
&=&1_{\{t\wedge \tau =t\}}\mathcal{E}_{s,t}^{\prime }[X].
\end{eqnarray*}
This with (\ref{e8.10a}) and the fact $\{t\wedge \tau =t\}\cup \{t\wedge
\tau \leq s\}=\Omega $ that (ii) holds. $\Box $\medskip\

It is clear that on $[t_i,t_{i+1}]$, $(\mathcal{E}_{s,t}^i[X])_{0\leq s\leq
t\leq T}$ is an $\mathcal{F}_t^\tau =\mathcal{F}_{\tau \wedge t}$%
--consistent evaluation. It immediately follows from Proposition \ref{m2.3}
and Remark \ref{m2.3r} that

\begin{lemma}
\label{m8.1a} $\{\mathcal{E}_{s,t}^i[\cdot ]\}_{i=0}^n$ defined in (\ref
{e8.10}) generates a unique $\mathcal{F}_{t\wedge \tau }$--consistent
evaluation on $[0,T]$. We denote this evaluation by
\begin{equation}
\mathcal{E}_{s\wedge \tau ,t\wedge \tau }[\cdot ]:L^2(\mathcal{F}_{\tau
\wedge t})\rightarrow L^2(\mathcal{F}_{\tau \wedge s}),\;0\leq s\leq t\leq T.
\label{e8.11}
\end{equation}
It satisfies (A1)--(A5) in the following sense: for each $0\leq s\leq t\leq
T $ and $X$, $X^{\prime }\in L^2(\mathcal{F}_{t\wedge \tau })$, \\\textbf{%
(a1)} $\mathcal{E}_{s\wedge \tau ,t\wedge \tau }[X]\geq \mathcal{E}_{s\wedge
\tau ,t\wedge \tau }[X^{\prime }]$, a.s., if $X\geq X^{\prime }$, a.s. \\%
\textbf{(a2)} $\mathcal{E}_{t\wedge \tau ,t\wedge \tau }[X]=X$; \\\textbf{%
(a3)} $\mathcal{E}_{r\wedge \tau ,s\wedge \tau }[\mathcal{E}_{s\wedge \tau
,t\wedge \tau }[X]]=\mathcal{E}_{r\wedge \tau ,t\wedge \tau }[X]$, $\forall
0\leq r\leq s\leq t$; \\\textbf{(a4')} $1_A\mathcal{E}_{s\wedge \tau
,t\wedge \tau }[X]=\mathcal{E}_{s\wedge \tau ,t\wedge \tau }[1_AX]$,\ $%
\forall A\in \mathcal{F}_{s\wedge \tau }$. \\\textbf{(a5)} For each $\tau
\in \mathcal{S}_T^0$, $X,X^{\prime }\in L^2(\mathcal{F}_{t\wedge \tau })$
and $K,\;K^{\prime }\in \mathcal{S}_T^2(0,T)$
\begin{equation}
\mathcal{E}_{s\wedge \tau ,t\wedge \tau }[X;K]-\mathcal{E}_{s\wedge \tau
,t\wedge \tau }[X^{\prime };K^{\prime }]\leq \mathcal{E}_{s\wedge \tau
,t\wedge \tau }^{g_\mu }[X-X^{\prime };K-K^{\prime }].  \label{e8.11a5}
\end{equation}
\end{lemma}

\smallskip\noindent\textbf{Proof. }We first prove that, for each $i$, $\mathcal{E}%
_{s,t}^i[\cdot ]=\mathcal{E}_{s\wedge \tau ,t\wedge \tau }[\cdot ]$, $%
t_i\leq s\leq t\leq t_{i+1}$ satisfies (a1)--(a3), (a4') and (a5).
(a1)--(a3) are easy to check. To prove (a4'), we observe that, for each $%
A\in \mathcal{F}_{s\wedge \tau }$,
\begin{eqnarray*}
A\cap \{t\wedge \tau &=&t\}=A\cap \{s\wedge \tau \leq s\}\cap \{t\wedge \tau
=t\}\in \mathcal{F}_s, \\
A\cap \{t\wedge \tau &\leq &s\}=A\cap \{s\wedge \tau \leq s\}\cap \{t\wedge
\tau \leq s\}\in \mathcal{F}_s.
\end{eqnarray*}
By (\ref{e8.10})--(ii), (a4) follows from
\begin{eqnarray*}
1_A\mathcal{E}_{s\wedge \tau ,t\wedge \tau }[X] &=&1_{\{t\wedge \tau \leq
s\}}1_AX+1_{\{t\wedge \tau =t\}}1_A\mathcal{E}_{s,t}[X] \\
&=&1_{\{t\wedge \tau \leq s\}}1_AX+1_{\{t\wedge \tau =t\}}1_A\mathcal{E}%
_{s,t}[1_{\{t\wedge \tau =t\}}1_AX] \\
&=&1_{\{t\wedge \tau \leq s\}}1_AX+1_{\{t\wedge \tau =t\}}1_A\mathcal{E}%
_{s,t}[1_AX] \\
&=&1_A\mathcal{E}_{s\wedge \tau ,t\wedge \tau }[1_AX].
\end{eqnarray*}
This with $\mathcal{E}_{s\wedge \tau ,t\wedge \tau }[0]=0$ yields (a4'). It
then follows from Proposition \ref{m2.3} and Remark \ref{m2.3r} that, there
exists a unique $\mathcal{F}_{t\wedge \tau }$--consistent evaluation
satisfying (a1)--(a4) that coincides to $\mathcal{E}^i$ on each $%
[t_i,t_{i+1}]$. (a4) plus $\mathcal{E}_{s\wedge \tau ,t\wedge \tau
}[0]\equiv 0$ implies (a4').

We now prove (a5). We only prove the second relation. The proof of the first
one is similar. We still let $\tau $ be characterized by (\ref{e8.4}). We
already have (a5) when $t_i\leq s\leq t\leq t_{i+1}$, for each $%
i=0,1,2,\cdots $. Now if $s\in [t_i,t_{i+1})$, $t\in (t_{i+1},t_{i+2}]$, for
some $i=0,1,2,\cdots $, we have
\begin{eqnarray*}
&&\mathcal{E}_{s\wedge \tau ,t\wedge \tau }[X;K]-\mathcal{E}_{s\wedge \tau
,t\wedge \tau }[X^{\prime };K^{\prime }] \\
&=&\mathcal{E}_{s\wedge \tau ,t_{i+1}\wedge \tau }[\mathcal{E}%
_{t_{i+1}\wedge \tau ,t\wedge \tau }[X;K]]-\mathcal{E}_{s\wedge \tau
,t_{i+1}\wedge \tau }[\mathcal{E}_{t_{i+1}\wedge \tau ,t\wedge \tau
}[X^{\prime };K^{\prime }]] \\
&\leq &\mathcal{E}_{s\wedge \tau ,t_{i+1}\wedge \tau }^{g_\mu }[\mathcal{E}%
_{t_{i+1}\wedge \tau ,t\wedge \tau }[X]-\mathcal{E}_{t_{i+1}\wedge \tau
,t\wedge \tau }[X^{\prime }];K-K^{\prime }] \\
&\leq &\mathcal{E}_{s\wedge \tau ,t_{i+1}\wedge \tau }^{g_\mu }[\mathcal{E}%
_{t_{i+1}\wedge \tau ,t\wedge \tau }^{g_\mu }[X-X^{\prime };K-K^{\prime
}];K-K^{\prime }] \\
&=&\mathcal{E}_{s\wedge \tau ,t\wedge \tau }^{g_\mu }[X-X^{\prime
};K-K^{\prime }].
\end{eqnarray*}
Thus the inequality is still true. We can repeat this procedure to conclude
that, for all $0\leq s\leq t\leq T$, the inequality holds. The proof of the
first inequality of (\ref{e8.11a5}) is analogous. $\Box $\medskip\

\begin{remark}
\label{r-rcll}Since $(\mathcal{E}_{s\wedge \tau ,\tau }[X;K])_{t\in [0,T]}$
is generated by finite steps from $\mathcal{E}_{\cdot ,t}[X;K]\in D_{%
\mathcal{F}}^2(0,t)$, thus $(\mathcal{E}_{s\wedge \tau ,\tau }[X;K])_{t\in
[0,T]}$ is still in $D_{\mathcal{F}}^2(0,T)$.
\end{remark}

\begin{lemma}
\label{m8.3a}If $\mathcal{E}[\cdot ]=\mathcal{E}^g[\cdot ]$, then for each $%
0\leq s\leq t\leq T$,
\begin{equation}
\mathcal{E}_{s\wedge \tau ,t\wedge \tau }^g[X]=\mathcal{E}_{s\wedge \tau
,t\wedge \tau }[X],\;\forall X\in \mathcal{F}_{t\wedge \tau }.  \label{e8.12}
\end{equation}
\end{lemma}

\smallskip\noindent\textbf{Proof. }By Lemma \ref{m8.1}, the evaluation $\mathcal{E}%
_{s\wedge \tau ,t\wedge \tau }[\cdot ]$ defined by (\ref{e8.10}) and (\ref
{e8.11}) coincides with $\mathcal{E}_{s\wedge \tau ,t\wedge \tau }^g[\cdot ]$
on $[t_i,t_{i+1}]$, for all $i=0,\cdots ,n$. It follows again from
Proposition \ref{m2.3} that (\ref{e8.12}) holds. $\Box $\medskip\

\begin{lemma}
\label{m8.a4}Let $\tau \in \mathcal{S}_T^0$ be characterized by (\ref{e8.4}%
). Then, for each $X\in \mathcal{F}_\tau $, $s$, $t\in [0,T]$ and $%
i=1,2,\cdots ,n$, we have
\begin{equation}
1_{\{t\wedge \tau =s\}}\mathcal{E}_{t\wedge \tau ,\tau }[X]=1_{\{t\wedge
\tau =s\}}\mathcal{E}_{s\wedge \tau ,\tau }[X].  \label{e8.10aa}
\end{equation}
\end{lemma}

\smallskip\noindent\textbf{Proof. }This problem can be divided into three cases: case
1: $t<s$. In this case (\ref{e8.10aa}) holds since $1_{\{t\wedge \tau
=s\}}\equiv 0$. Case 2: $t=s$. (\ref{e8.10aa}) is clearly true. We now
consider the last case: $t>s$. We assume that $t\in (t_k,t_{k+1}]$, with $%
i\leq k\leq n$. In this case we have, by Lemma \ref{m8.1a} and (\ref{e8.10}%
)--(ii),
\begin{eqnarray*}
1_{\{t\wedge \tau =s\}}\mathcal{E}_{s\wedge \tau ,\tau }[X] &=&1_{\{t\wedge
\tau =s\}}\mathcal{E}_{s\wedge \tau ,\tau }[\cdots \mathcal{E}_{t_k\wedge
\tau ,t\wedge \tau }[\mathcal{E}_{t\wedge \tau ,\tau }[X]]\cdots ] \\
&=&1_{\{t\wedge \tau =s\}}1_{\{t\wedge \tau \leq s\}}\mathcal{E}_{s\wedge
\tau ,t_{i+1}\wedge \tau }[\cdots \mathcal{E}_{t_k\wedge \tau ,t\wedge \tau
}[\mathcal{E}_{t\wedge \tau ,\tau }[X]]\cdots ] \\
&=&1_{\{t\wedge \tau =s\}}1_{\{t\wedge \tau \leq t_{i+1}\}}\mathcal{E}%
_{t_{i+1}\wedge \tau ,t_{i+2}\wedge \tau }[\cdots \mathcal{E}_{t_k\wedge
\tau ,t\wedge \tau }[\mathcal{E}_{t\wedge \tau ,\tau }[X]]\cdots ] \\
&&\vdots \\
&=&1_{\{t\wedge \tau =s\}}1_{\{t\wedge \tau \leq t_k\}}\mathcal{E}%
_{t_k\wedge \tau ,t\wedge \tau }[\mathcal{E}_{t\wedge \tau ,\tau }[X]] \\
&=&1_{\{t\wedge \tau =s\}}\mathcal{E}_{t\wedge \tau ,\tau }[X].
\end{eqnarray*}

$\Box $\medskip\

We now consider the general case of $\mathcal{E}_{\sigma ,\tau }[\cdot ]$
for $\sigma $, $\tau \in \mathcal{S}_T^0$ with $\sigma \leq \tau $. Let $%
\sigma $ be characterized by
\begin{equation}
\;\cup _{i=1}^m\{\sigma =s_i\}=\Omega ,\;0\leq s_1<\cdots <s_m\leq T.
\label{e8.15}
\end{equation}

\begin{definition}
\label{m8.4a} Let $\sigma $ be characterized by (\ref{e8.15}). $\mathcal{E}%
_{\sigma ,\tau }[\cdot ]$ is defined by, as in classical situations,
\begin{equation}
\mathcal{E}_{\sigma ,\tau }[X]:=\mathcal{E}_{t\wedge \tau ,\tau
}[X]|_{t=\sigma }=\sum_{i=1}^m1_{\{\sigma =s_i\}}\mathcal{E}_{s_i\wedge \tau
,\tau }[X],\;X\in \mathcal{F}_\tau .  \label{e8.16}
\end{equation}
\end{definition}

\begin{remark}
By Lemma \ref{m8.a4}, it is clear that (\ref{e8.16}) is satisfied in the
case $\sigma =t\wedge \tau $.
\end{remark}

\begin{lemma}
\label{m8.4b}Let $\sigma $, $\tau \in \mathcal{S}_T^0$ be such that $\sigma
\leq \tau $. Then for each $0\leq s\leq t\leq T$ and $X\in L^2(\mathcal{F}%
_\tau )$ we have
\begin{equation}
\mathcal{E}_{s\wedge \sigma ,t\wedge \sigma }[\mathcal{E}_{t\wedge \sigma
,\tau }[X]]=\mathcal{E}_{s\wedge \sigma ,\tau }[X].  \label{e8.17}
\end{equation}
\end{lemma}

\smallskip\noindent\textbf{Proof. }Without loss of generality we let both $\sigma $
and $\tau $ be valued in $\{t_0,t_1,\cdots ,t_n\}$ with
\begin{equation}
0=t_0\leq t_1<\cdots <t_n\leq t_{n+1}=T.  \label{e8.18}
\end{equation}
($\{\sigma =t_i\}$ or $\{\tau =t_i\}$ may be an empty set for some $i$). For
a fixed $i\in \{0,1,\cdots ,n-1\}$, we consider, the case $t_i\leq s<t\leq
t_{i+1}$ (the case $s=t$ is clearly true). By (\ref{e8.10})--(ii), we have
\begin{eqnarray}
\mathcal{E}_{s\wedge \sigma ,t\wedge \sigma }[\mathcal{E}_{t\wedge \sigma
,\tau }[X]] &=&1_{\{t\wedge \sigma \leq s\}}\mathcal{E}_{t\wedge \sigma
,\tau }[X]+1_{\{t\wedge \sigma =t\}}\mathcal{E}_{s,t}[\mathcal{E}_{t\wedge
\sigma ,\tau }[X]]  \label{e8.22a} \\
\ &=&1_{\{t\wedge \sigma \leq s\}}\mathcal{E}_{s\wedge \sigma ,\tau
}[X]+1_{\{t\wedge \sigma =t\}}\mathcal{E}_{s,t}[1_{\{t\wedge \sigma =t\}}%
\mathcal{E}_{t\wedge \sigma ,\tau }[X]].\mathcal{\ }  \nonumber
\label{e8.22b}
\end{eqnarray}
The second term is due to the assumption (A4) of $\mathcal{E}$, with the
observation that $1_{\{t\wedge \sigma =t\}}$ is $\mathcal{F}_{t_i}$ (and
thus $\mathcal{F}_s$) measurable. We repeatedly use (\ref{e8.16}) and (A4)
to the second term
\begin{eqnarray}
1_{\{t\wedge \sigma =t\}}\mathcal{E}_{s,t}[1_{\{t\wedge \sigma =t\}}\mathcal{%
E}_{t\wedge \sigma ,\tau }[X]] &=&1_{\{t\wedge \sigma =t\}}\mathcal{E}%
_{s,t}[1_{\{t\wedge \sigma =t\}}\mathcal{E}_{t\wedge \tau ,\tau }[X]]
\label{e8.22c} \\
\ &=&1_{\{t\wedge \sigma =t\}}\mathcal{E}_{s,t}[\mathcal{E}_{t\wedge \tau
,\tau }[X]]  \nonumber \\
\ &=&1_{\{t\wedge \sigma =t\}}1_{\{t\wedge \tau =t\}}\mathcal{E}_{s,t}[%
\mathcal{E}_{t\wedge \tau ,\tau }[X]]  \nonumber \\
\ &=&1_{\{t\wedge \sigma =t\}}1_{\{t\wedge \tau =t\}}\mathcal{E}_{s\wedge
\tau ,t\wedge \tau }[\mathcal{E}_{t\wedge \tau ,\tau }[X]]  \nonumber \\
\ &=&1_{\{t\wedge \sigma =t\}}1_{\{s\wedge \sigma =s\}}\mathcal{E}_{s\wedge
\tau ,\tau }[X].  \nonumber
\end{eqnarray}
Here we use the fact that $\{t\wedge \sigma =t\}\subset \{t\wedge \tau =t\}$%
, since $\tau \geq \sigma $, and $\{t\wedge \sigma =t\}\subset \{s\wedge
\sigma =s\}$. By (\ref{e8.16}), $1_{\{s\wedge \sigma =s\}}\mathcal{E}%
_{s\wedge \tau ,\tau }[X]=1_{\{s\wedge \sigma =s\}}\mathcal{E}_{s\wedge
\sigma ,\tau }[X]$. Thus (\ref{e8.22a}) finally becomes
\begin{eqnarray*}
\mathcal{E}_{s\wedge \sigma ,t\wedge \sigma }[\mathcal{E}_{t\wedge \sigma
,\tau }[X]] &=&1_{\{t\wedge \sigma \leq s\}}\mathcal{E}_{s\wedge \sigma
,\tau }[X]+1_{\{t\wedge \sigma =t\}}\mathcal{E}_{s\wedge \sigma ,\tau }[X] \\
\ &=&\mathcal{E}_{s\wedge \sigma ,\tau }[X],
\end{eqnarray*}
since $\{t\wedge \sigma \leq s\}+\{t\wedge \sigma =t\}=\Omega $.

We have proved (\ref{e8.17}) for the case $t_i\leq s\leq t\leq t_{i+1}$, $%
i\in \{0,1,\cdots ,n\}$. For the general situation, say $s\in [t_i,t_{i+1})$%
, $t\in [t_j,t_{j+1}]$, $0\leq i<j\leq n$, (\ref{e8.17}) can be deduced by
\begin{eqnarray*}
&&\mathcal{E}_{s\wedge \sigma ,t\wedge \sigma }[\mathcal{E}_{t\wedge \sigma
,\tau }[X]] \\
&=&\mathcal{E}_{s\wedge \sigma ,t_i\wedge \sigma }[\mathcal{E}_{t_i\wedge
\sigma ,t_{i+1}\wedge \sigma }[\cdots \mathcal{E}_{t_j\wedge \sigma ,t\wedge
\sigma }[\mathcal{E}_{t\wedge \sigma ,\tau }[X]]]] \\
&=&\mathcal{E}_{s\wedge \sigma ,t_i\wedge \sigma }[\mathcal{E}_{t_i\wedge
\sigma ,t_{i+1}\wedge \sigma }[\cdots \mathcal{E}_{t_j\wedge \sigma ,\tau
}[X]]] \\
&&\vdots \\
&=&\mathcal{E}_{s\wedge \sigma ,t_i\wedge \sigma }[\mathcal{E}_{t_i\wedge
\sigma ,\tau }[X]] \\
&=&\mathcal{E}_{s\wedge \sigma ,\tau }[X].
\end{eqnarray*}
The proof is complete. $\Box $\medskip\

By this result and the definition of $\mathcal{E}_{\sigma ,\tau }[\cdot ]$
in (\ref{e8.16}), we immediately have

\begin{lemma}
\label{m8.4c}Let $\rho $, $\sigma $ and $\tau \in \mathcal{S}_T^0$ be such
that $\rho \leq \sigma \leq \tau $. Then for each $X\in L^2(\mathcal{F}_\tau
)$ we have
\begin{equation}
\mathcal{E}_{\rho ,\sigma }[\mathcal{E}_{\sigma ,\tau }[X]]=\mathcal{E}%
_{\sigma ,\tau }[X].  \label{e8.17a}
\end{equation}
\end{lemma}

We now consider an $\mathcal{E}$--supermartingale. We will prove the
following optional stopping theorem

\begin{lemma}
\label{m8.5}Let $Y\in D_{\mathcal{F}}^2(0,T)$ be an $\mathcal{E}$%
--martingale (respectively $\mathcal{E}$--supermartingale, $\mathcal{E}$%
--submartingale). Then for each $\sigma $, $\tau \in \mathcal{S}_T^0$ such
that $\sigma \leq \tau $, we have
\begin{equation}
\mathcal{E}_{\sigma ,\tau }[Y_\tau ]=Y_\sigma \hbox{, (resp. }\leq
Y_\sigma \hbox{, }\geq Y_\sigma \hbox{) a.s.}  \label{e8.24}
\end{equation}
\end{lemma}

\smallskip\noindent\textbf{Proof. }We only prove the case for $\mathcal{E}$%
--supermartingale. It is clear that, once we have
\begin{equation}
\mathcal{E}_{t\wedge \tau ,\tau }[Y_\tau ]\leq Y_{t\wedge \tau }\hbox{, \ }%
\forall t\in [0,T],  \label{e8.25}
\end{equation}
then, by (\ref{e8.16}), we can also prove (\ref{e8.24}). We still let $\tau $
be characterized by (\ref{e8.4}). We will prove this inequality by
induction. Firstly, when $t\geq t_n$, (\ref{e8.25}) holds since $\mathcal{E}%
_{t\wedge \tau ,\tau }[Y_\tau ]=\mathcal{E}_{\tau ,\tau }[Y_\tau ]=Y_\tau $.
Now suppose that for a fixed $i\in \{1,\cdots ,n\}$, (\ref{e8.25}) holds for
$t\geq t_i$, we shall prove that it also holds for $t\geq t_{i-1}$. We need
to check the case $t\in [t_{i-1},t_i)$. By (\ref{e8.10})--(ii) and applying
(A4) (since $1_{\{t_i\wedge \tau =t_i\}}$ is $\mathcal{F}_t$--measurable),
we have
\begin{eqnarray*}
\mathcal{E}_{t\wedge \tau ,t_i\wedge \tau }[Y_{t_i\wedge \tau }]
&=&1_{\{t_i\wedge \tau \leq t\}}Y_{t_i\wedge \tau }+1_{\{t_i\wedge \tau
=t_i\}}\mathcal{E}_{t,t_i}[Y_{t_i\wedge \tau }] \\
\ &=&1_{\{t_i\wedge \tau \leq t\}}Y_{t_i\wedge \tau }+1_{\{t_i\wedge \tau
=t_i\}}\mathcal{E}_{t,t_i}[1_{\{t_i\wedge \tau =t_i\}}Y_{t_i}] \\
\ &=&1_{\{t_i\wedge \tau \leq t\}}Y_{t_i\wedge \tau }+1_{\{t_i\wedge \tau
=t_i\}}\mathcal{E}_{t,t_i}[Y_{t_i}] \\
\ &\leq &1_{\{t_i\wedge \tau \leq t\}}Y_{t_i\wedge \tau }+1_{\{t_i\wedge
\tau =t_i\}}Y_t \\
\ &=&Y_{t\wedge \tau }.
\end{eqnarray*}
The last step is from $\{t_i\wedge \tau \leq t\}+\{t_i\wedge \tau
=t_i\}=\Omega $ and then $t\wedge \tau =t_i\wedge \tau 1_{\{t_i\wedge \tau
\leq t\}}+t1_{\{t_i\wedge \tau =t_i\}}$. From this result we derive
\begin{eqnarray*}
\mathcal{E}_{t\wedge \tau ,\tau }[Y_\tau ] &=&\mathcal{E}_{t\wedge \tau
,t_i\wedge \tau }[\mathcal{E}_{t_i\wedge \tau ,\wedge \tau }[Y_\tau ]] \\
\ &\leq &\mathcal{E}_{t\wedge \tau ,t_i\wedge \tau }[Y_{t_i\wedge \tau }] \\
\ &\leq &Y_{t\wedge \tau }.
\end{eqnarray*}
Thus (\ref{e8.25}) holds for $t\geq t_{i-1}$. It follows by induction that (%
\ref{e8.25}) holds for $t\in [0,T]$. The proof is complete. $\Box $\medskip\

We conclude

\begin{lemma}
\label{m8.5a}The system of operators
\[
\mathcal{E}_{\sigma ,\tau }[\cdot ]:L^2(\mathcal{F}_\tau )\rightarrow L^2(%
\mathcal{F}_\sigma ),\;\sigma \leq \tau ,\;\sigma ,\tau \in \mathcal{S}_T^0
\]
is an $\mathcal{F}_t$--consistent nonlinear evaluation in the following
sense: for each $X$, $X^{\prime }\in L^2(\mathcal{F}_\tau )$, \\\textbf{(a1)}
$\mathcal{E}_{\sigma ,\tau }[X]\geq \mathcal{E}_{\sigma ,\tau }[X^{\prime }]$%
, a.s., if $X\geq X^{\prime }$, a.s. \\\textbf{(a2)} $\mathcal{E}_{\tau
,\tau }[X]=X$; \\\textbf{(a3)} $\mathcal{E}_{\rho ,\sigma }[\mathcal{E}%
_{\sigma ,\tau }[X]]=\mathcal{E}_{\rho ,\tau }[X]$, $\forall 0\leq \rho \leq
\sigma \leq \tau $; \\\textbf{(a4')} $1_A\mathcal{E}_{\sigma ,\tau }[X]=%
\mathcal{E}_{\sigma ,\tau }[1_AX]$,\ $\forall A\in \mathcal{F}_\sigma $. \\%
We also have\\\textbf{(a5)} For each $0\leq \sigma \leq \tau \leq T$,
\begin{equation}
\mathcal{E}_{\sigma ,\tau }[X;K]-\mathcal{E}_{\sigma ,\tau }[X^{\prime
};K^{\prime }]\leq \mathcal{E}_{\sigma ,\tau }^{g_\mu }[X-X^{\prime
};K-K^{\prime }],\;\;\forall X,X^{\prime }\in L^2(\mathcal{F}_\tau ).
\label{e8.26}
\end{equation}
Moreover, $\mathcal{E}_{\sigma ,\tau }[\cdot ]$ is the unique extension of $%
\mathcal{E}_{s,t}[\cdot ]$ in the following sense: for each system of
operator
\[
\mathcal{E}_{\sigma ,\tau }^{\prime }[\cdot ]:L^2(\mathcal{F}_\tau
)\rightarrow L^2(\mathcal{F}_\sigma ),\;\sigma \leq \tau ,\;\sigma ,\tau \in
\mathcal{S}_T^0
\]
satisfying satisfies (a1)--(a4') such $\mathcal{E}_{s,t}^{\prime }[X]=%
\mathcal{E}_{s,t}[X]$, for each (deterministic) $0\leq s\leq t\leq T$ and
for each $X\in \mathcal{F}_t$, then, for each $\sigma \leq \tau ,\;\sigma
,\tau \in \mathcal{S}_T^0$ and for each $X\in \mathcal{F}_\tau $, we have $%
\mathcal{E}_{\sigma ,\tau }^{\prime }[X]=\mathcal{E}_{\sigma ,\tau }^{\prime
}[X]$.\
\end{lemma}

\smallskip\noindent\textbf{Proof. }(a1) and (a2) are easily checked from Definition (%
\ref{e8.16}) of $\mathcal{E}_{\sigma ,\tau }[\cdot ]$ and Lemma \ref{m8.1a}.
(a3) is given in Lemma \ref{m8.4c}. (a5) can be proved by using (\ref{e8.16}%
) and the (a5) part of Lemma \ref{m8.1a}. It remains to prove (a4'). Let $%
A\in \mathcal{F}_\sigma $ and let $\sigma $ be characterized by (\ref{e8.15}%
). From (\ref{e8.16}) and the fact that $A\cap \{\sigma =s_i\}\in \mathcal{F}%
_{s_i\wedge \sigma }$, we derive, using the (a4') part of Lemma \ref{m8.1a},
\begin{eqnarray*}
1_A\mathcal{E}_{\sigma ,\tau }[X] &=&\sum_{i=1}^m1_A1_{\{\sigma =s_i\}}%
\mathcal{E}_{s_i\wedge \tau ,\tau }[X] \\
\ &=&\sum_{i=1}^m1_A1_{\{\sigma =s_i\}}\mathcal{E}_{s_i\wedge \tau ,\tau
}[1_A1_{\{\sigma =s_i\}}X] \\
\ &=&\sum_{i=1}^m1_A1_{\{\sigma =s_i\}}\mathcal{E}_{s_i\wedge \tau ,\tau
}[1_AX] \\
\ &=&1_A\mathcal{E}_{\sigma ,\tau }[1_AX].
\end{eqnarray*}
This with $\mathcal{E}_{\sigma ,\tau }[0]\equiv 0$ implies (a4'). The
uniqueness is a direct consequence of Lemma \ref{m8.a3} and the uniqueness
part of Lemma \ref{m8.1a}. $\Box $\medskip\

\subsection{$\mathcal{S}_T$ case: Proof of Theorem \ref{m8.13} and optional
stopping theorem\label{ss9.3}}

We now extend $\mathcal{E}_{\sigma ,\tau }[.]$ from $\sigma $, $\tau \in
\mathcal{S}_T^0$ to the general case, i.e., $\sigma $, $\tau \in \mathcal{S}%
_T$. Firstly, by Remark \ref{r-rcll}, for each $\tau \in \mathcal{S}_T^0$,
we have $\mathcal{E}_{\cdot \wedge \tau ,\tau }[X]\in D_{\mathcal{F}}^2(0,T)$%
. This with the definition of of $\mathcal{E}_{\sigma ,\tau }[X]$ it follows
that

\begin{lemma}
\label{m8.9}Let $\sigma \in \mathcal{S}_T$, $\tau \in \mathcal{S}_T^0$ be
such that $\sigma \leq \tau $ and let $X\in L^2(\mathcal{F}_\tau )$. Then
for each sequence $\{\sigma _n\}_{n=1}^\infty $ of $\mathcal{S}_T^0$ such
that $\sigma \leq \sigma _n\leq \tau $ and $\lim_{n\rightarrow \infty
}\sigma _n=\sigma $, a.s., $\{\mathcal{E}_{\sigma _n,\tau
}[X]\}_{n=1}^\infty $ is a Cauchy sequence in $L^2(\mathcal{F}_T)$. Moreover
the limit of this sequence is in $L^2(\mathcal{F}_\sigma )$. We denote it by
$\mathcal{E}_{\sigma ,\tau }[X]_{0\leq \sigma \leq \tau \leq T}$:
\begin{equation}
\mathcal{E}_{\sigma ,\tau }[\cdot ]:L^2(\mathcal{F}_\tau )\rightarrow L^2(%
\mathcal{F}_\sigma )\hbox{,\ }\sigma \in \mathcal{S}_T,\;\tau \in \mathcal{S}%
_T^0.  \label{e8.31}
\end{equation}
\end{lemma}

With the above convergence result and Lemma \ref{m8.5a}, we can easily have

\begin{lemma}
\label{m8.10}The system of operators $\mathcal{E}_{\sigma ,\tau }[\cdot ]$,
defined in (\ref{e8.31}) satisfies the following properties: for each $%
\sigma \in \mathcal{S}_T$, $\tau \in \mathcal{S}_T^0$, $\sigma \leq \tau $,
and $X$, $X^{\prime }\in L^2(\mathcal{F}_\tau )$, we have\\\textbf{(a1)} $%
\mathcal{E}_{\sigma ,\tau }[X]\geq \mathcal{E}_{\sigma ,\tau }[X^{\prime }]$%
, a.s., if $X\geq X^{\prime }$, a.s. \\\textbf{(a3)} $\mathcal{E}_{\sigma
,\rho }[\mathcal{E}_{\rho ,\tau }[X]]=\mathcal{E}_{\sigma ,\tau }[X]$, $%
\forall 0\leq \sigma \leq \rho \leq \tau $, $\rho \in \mathcal{S}_T^0$;\\%
\textbf{(a4')} $1_A\mathcal{E}_{\sigma ,\tau }[X]=\mathcal{E}_{\sigma ,\tau
}[1_AX]$,\ $\forall A\in \mathcal{F}_\sigma $. \\\textbf{(a5)} For each $%
0\leq \sigma \leq \tau \leq T$ and $X,X^{\prime }\in L^2(\mathcal{F}_\tau ),$
\begin{equation}
\mathcal{E}_{\sigma ,\tau }[X;K]-\mathcal{E}_{\sigma ,\tau }[X^{\prime
};K]\leq \mathcal{E}_{\sigma ,\tau }^{g_\mu }[X-X^{\prime };K-K^{\prime }].
\label{e8.32}
\end{equation}
Consequently, the estimates in Lemma \ref{m8.8} still hold for $\sigma \in
\mathcal{S}_T$ and $\tau \in \mathcal{S}_T^0$.
\end{lemma}

To proceed, we need the following estimates

\begin{lemma}
\label{m8.7}For each $\sigma $, $\tau \in \mathcal{S}_T$, $\sigma \leq \tau $
and $X\in L^2(\mathcal{F}_\sigma )$, we have the following estimate
\begin{equation}
E[|\mathcal{E}_{\sigma ,\tau }^g[X]-X|^2]\leq CE\int_\sigma ^\tau
|g(s,0,0)|^2ds+CE[(\tau -\sigma )|X|^2]  \label{e8.28}
\end{equation}
as well as
\begin{equation}
E[\sup_{s\in [0,T]}|\mathcal{E}_{t\wedge \tau ,\tau }^g[X]|^2]\leq
CE[|X|^2]+CE\int_0^\tau |g(s,0,0)|^2)ds,  \label{e8.28a}
\end{equation}
where the constant $C$ depends only on $\mu $ and $T$.
\end{lemma}

\smallskip\noindent\textbf{Proof. }(\ref{e8.28a}) is a special case of (\ref{e2.1})
with $K_t\equiv 0$. In order to prove (\ref{e8.28}), we set $\bar y_t\equiv
X $ on $[\sigma ,\tau ]$. Observe that $\mathcal{E}_{\sigma ,\tau
}^g[X]=y_\sigma $, where $(y_t)_{t\in [0,T]}$ is the solution of the BSDE
\[
y_t=X+\int_t^T1_{[\sigma ,\tau ]}(s)g(s,y_s,z_s)ds-\int_t^Tz_sdB_s.
\]
We apply It\^o's formula to $|y_t-\bar y_t|^2e^{\beta t}$ on $[\sigma ,\tau
] $, with $\beta =2\mu +3\mu ^2$. By $|g(s,y_s,z_s)|\leq \mu (|y_s-\bar
y_s|+|\bar y_s|+|z_s|)+|g_s^0|$ with $g_s^0=g(s,0,0)$, we have
\begin{eqnarray*}
&&\ \ E[|y_\sigma -\bar y_\sigma |^2e^{\beta \sigma }]+E\int_\sigma ^\tau
[\beta |y_s-\bar y_s|^2+|z_s|^2]e^{\beta s}ds \\
\ &=&E\int_\sigma ^\tau 2(y_s-\bar y_s)g(s,y_s,z_s)e^{\beta s}ds \\
\ &\leq &E\int_\sigma ^\tau 2|y_s-\bar y_s|\mu (|y_s-\bar y_s|+|\bar
y_s|+|z_s|+|g_s^0|)e^{\beta s}ds \\
\ &\leq &E\int_\sigma ^\tau [(2\mu +3\mu ^2)|y_s-\bar y_s|^2+|\bar
y_s|^2+|z_s|^2+|g_s^0|^2]e^{\beta s}ds.
\end{eqnarray*}
With $\bar y_s\equiv X$, It follows that
\[
E[|y_\sigma -X|^2e^{\beta \sigma }]\leq E\int_\sigma ^\tau
(|X|^2+|g_s^0|^2)e^{\beta s}ds.
\]
Thus we have (\ref{e8.28}) with $C=e^{\beta T}$. $\Box $\medskip\

\begin{corollary}
\label{m8.8}For each $\sigma $, $\tau \in \mathcal{S}_T^0$, $\sigma \leq
\tau $ and $X$, $X^{\prime }\in L^2(\mathcal{F}_\tau )$, we have the
following estimates
\begin{equation}
E[\sup_{s\in [0,T]}|\mathcal{E}_{t\wedge \tau ,\tau }[X]|^2]\leq CE[|X|^2],
\label{e8.30}
\end{equation}
\begin{equation}
|\mathcal{E}_{\sigma ,\tau }[X]-\mathcal{E}_{\sigma ,\tau }[X^{\prime
}]|^2\leq CE[|X-X^{\prime }|^2]  \label{e8.30a}
\end{equation}
and, for each $X\in L^2(\mathcal{F}_\sigma )$,
\begin{equation}
E[|\mathcal{E}_{\sigma ,\tau }[X]-X|^2]\leq CE[(\tau -\sigma )|X|^2],
\label{e8.29}
\end{equation}
where the constant $C$ depends only on $\mu $ and $T$.
\end{corollary}

\smallskip\noindent\textbf{Proof. }By (a5) of Lemma \ref{m8.5a},
\[
\mathcal{E}_{\sigma ,\tau }^{-g_\mu }[X]-X\leq \mathcal{E}_{\sigma ,\tau
}[X]-X\leq \mathcal{E}_{\sigma ,\tau }^{g_\mu }[X]-X.
\]
Thus
\[
|\mathcal{E}_{\sigma ,\tau }[X]-X|\leq |\mathcal{E}_{\sigma ,\tau }^{g_\mu
}[X]-X|+|\mathcal{E}_{\sigma ,\tau }^{-g_\mu }[X]-X|.
\]
It then follows from Lemma \ref{m8.7} that (\ref{e8.29}) holds for $X\in L^2(%
\mathcal{F}_\sigma )$. Using (a5), the proofs of (\ref{e8.30}) and (\ref
{e8.30a}) are similar. Here we need the estimates (\ref{e2.gm1}). $\Box $%
\medskip\

To extend $\tau $ to $\mathcal{S}_T$ we need

\begin{lemma}
\label{m8.a2}Let $\sigma \in \mathcal{S}_T$, $\tau $, $\tau ^{\prime }\in
\mathcal{S}_T^0$ be such that $\sigma \leq \tau \vee \tau ^{\prime }$, and
let $X\in L^2(\mathcal{F}_{\tau \wedge \tau ^{\prime }})$. Then we have
\begin{equation}
E[\ |\mathcal{E}_{\sigma ,\tau }[X]-\mathcal{E}_{\sigma ,\tau ^{\prime
}}[X]|^2]\leq cE[(\tau \vee \tau ^{\prime }-\tau \wedge \tau ^{\prime
})|X|^2],  \label{e8.a2}
\end{equation}
where $c$ depends only on $\mu $ and $T$.
\end{lemma}

\smallskip\noindent\textbf{Proof. }We have
\begin{eqnarray*}
\ |\mathcal{E}_{\sigma ,\tau }[X]-\mathcal{E}_{\sigma ,\tau ^{\prime }}[X]|
&\leq &|\mathcal{E}_{\sigma ,\tau \wedge \tau ^{\prime }}[\mathcal{E}_{\tau
\wedge \tau ^{\prime },\tau }[X]]-\mathcal{E}_{\sigma ,\tau \wedge \tau
^{\prime }}[X]| \\
&&+|\mathcal{E}_{\sigma ,\tau \wedge \tau ^{\prime }}[X]-\mathcal{E}_{\sigma
,\tau \wedge \tau ^{\prime }}[\mathcal{E}_{\tau \wedge \tau ^{\prime },\tau
^{\prime }}[X]]|
\end{eqnarray*}
For the first term, by (\ref{e8.30a}) and then (\ref{e8.29}),
\begin{eqnarray*}
&&\ \ E[|\mathcal{E}_{\sigma ,\tau \wedge \tau ^{\prime }}[\mathcal{E}_{\tau
\wedge \tau ^{\prime },\tau }[X]]-\mathcal{E}_{\sigma ,\tau \wedge \tau
^{\prime }}[X]|^2] \\
\ &\leq &CE[|\mathcal{E}_{\tau \wedge \tau ^{\prime },\tau }[X]-X|^2] \\
\ &\leq &C^2E[(\tau -\tau \wedge \tau ^{\prime })|X|^2].
\end{eqnarray*}
Similarly
\begin{eqnarray*}
&&\ E[|\mathcal{E}_{\sigma ,\tau \wedge \tau ^{\prime }}[X]-\mathcal{E}%
_{\sigma ,\tau \wedge \tau ^{\prime }}[\mathcal{E}_{\tau \wedge \tau
^{\prime },\tau ^{\prime }}[X]]|^2] \\
\ &\leq &CE[|\mathcal{E}_{\tau \wedge \tau ^{\prime },\tau ^{\prime
}}[X]-X|^2] \\
\ &\leq &C^2E[(\tau ^{\prime }-\tau \wedge \tau ^{\prime })|X|^2].
\end{eqnarray*}
From the above three inequalities we have (\ref{e8.a2}). $\Box $\medskip\

By this estimate we have

\begin{lemma}
\label{m8.11}Let $\sigma $, $\tau \in \mathcal{S}_T$, be such that $\sigma
\leq \tau $ and let $X\in L^2(\mathcal{F}_{\tau \wedge \tau ^{\prime }})$.
Then for each sequence $\{\tau _n\}_{n=1}^\infty $ in $\mathcal{S}_T^0$ such
that $\tau \leq \tau _n$ and $\lim_{n\rightarrow \infty }\tau _n=\tau $,
a.s., the sequence $\{\mathcal{E}_{\sigma ,\tau _n}[X]\}_{n=1}^\infty $ is a
Cauchy sequence in $L^2(\mathcal{F}_\sigma )$.
\end{lemma}

\smallskip\noindent\textbf{Proof. }By (\ref{e8.a2})
\begin{equation}
E[\ |\mathcal{E}_{\sigma ,\tau _m}[X]-\mathcal{E}_{\sigma ,\tau
_n}[X]|^2]\leq cE[(\tau _m\vee \tau _n-\tau )|X|^2].  \label{e8.a2a}
\end{equation}
We then have that $\{\mathcal{E}_{\sigma ,\tau _n}[X]\}_{n=1}^\infty $ is a
Cauchy sequence in $L^2(\mathcal{F}_\sigma )$. $\Box $\medskip\

\begin{definition}
\label{m8.12}We denote the limit of the sequence $\{\mathcal{E}_{\sigma
,\tau _n}[X]\}_{n=1}^\infty $ in $L^2(\mathcal{F}_\sigma )$ of the above
Lemma $\mathcal{E}_{\sigma ,\tau }[X]$:
\begin{equation}
\mathcal{E}_{\sigma ,\tau }[\cdot ]:L^2(\mathcal{F}_\tau )\rightarrow L^2(%
\mathcal{F}_\sigma )\hbox{,\ }\sigma ,\tau \in \mathcal{S}_T.
\label{e8.33}
\end{equation}
\end{definition}

\noindent \textbf{Proof of Theorem} \textbf{\ref{m8.13}}. We have already
defined $\mathcal{E}_{\sigma ,\tau }$ in (\ref{e8.33}). With which (A1),
(A4') and (A5) are proved by simply using Lemma \ref{m8.10} and Lemma \ref
{m8.11} and by passing to the limit. (A2) is proved by
\[
0=\mathcal{E}_{\tau ,\tau }^{-g_\mu }[X]-X\leq \mathcal{E}_{\tau
,\tau }[X]-X\leq \mathcal{E}_{\tau ,\tau }^{g_\mu
}[X]-X=0,\;\hbox{a.s.}
\]
Once we have these properties, it is easy to check that the estimates (\ref
{e8.29})--(\ref{e8.30a}) still hold for $\sigma $, $\tau \in \mathcal{S}_T$.

We now prove (A3), i.e.,
\begin{equation}
\mathcal{E}_{\rho ,\sigma }[\mathcal{E}_{\sigma ,\tau }[X]]=\mathcal{E}%
_{\rho ,\tau }[X],\forall 0\leq \rho \leq \sigma \leq \tau .  \label{e8.35}
\end{equation}
We first prove this relation for the case $\rho $, $\sigma \in \mathcal{S}_T$
and $\tau \in \mathcal{S}_T^0$. Let $\{\sigma _n\}_{n=1}^\infty $ be a
sequence in $\mathcal{S}_T^0$ such that $\sigma \leq \sigma _n\leq \tau $, $%
n=1,2,\cdots $ and $\lim_{n\rightarrow \infty }\sigma _n=\sigma $, a.s.. By
Lemma \ref{m8.9} and Lemma \ref{m8.11}, we have
\begin{equation}
\lim_{n\rightarrow \infty }E[|\mathcal{E}_{\rho ,\sigma _n}[\mathcal{E}%
_{\sigma ,\tau }[X]]-\mathcal{E}_{\rho ,\sigma }[\mathcal{E}_{\sigma ,\tau
}[X]]|^2]=0.  \label{e8.36}
\end{equation}
and
\begin{equation}
\lim_{n\rightarrow \infty }E[|\mathcal{E}_{\sigma _n,\tau }[X]-\mathcal{E}%
_{\sigma ,\tau }[X]|^2]=0.  \label{e8.37}
\end{equation}
On the other hand, by (a3) of Lemma \ref{m8.10}, we have
\begin{eqnarray}
&&\mathcal{E}_{\rho ,\tau }[X]-\mathcal{E}_{\rho ,\sigma _n}[\mathcal{E}%
_{\sigma ,\tau }[X]]  \label{e8.38} \\
&=&\mathcal{E}_{\rho ,\sigma _n}[\mathcal{E}_{\sigma _n,\tau }[X]]-\mathcal{E%
}_{\rho ,\sigma _n}[\mathcal{E}_{\sigma ,\tau }[X]].  \nonumber
\end{eqnarray}
It follows from (\ref{e8.30a}) and (\ref{e8.37}) that, as $n$ tends to
infinity,
\[
E[|\mathcal{E}_{\rho ,\tau }[X]-\mathcal{E}_{\rho ,\sigma _n}[\mathcal{E}%
_{\sigma ,\tau }[X]]|^2]\leq CE[|\mathcal{E}_{\sigma _n,\tau }[X]-\mathcal{E}%
_{\sigma ,\tau }[X]|^2]\rightarrow 0.
\]
From this and (\ref{e8.36}) it follows that (\ref{e8.35}) holds for $\rho $,
$\sigma \in \mathcal{S}_T$ and $\tau \in \mathcal{S}_T^0$.

We now prove this relation for the general case: $\rho $, $\sigma $ and $%
\tau \in \mathcal{S}_T$. Let $\{\tau _n\}_{n=1}^\infty $ be a sequence in $%
\mathcal{S}_T^0$ such that $\tau \leq \tau _n$, $n=1,2,\cdots $ and $%
\lim_{n\rightarrow \infty }\tau _n=\tau $, a.s.. We have
\begin{equation}
\mathcal{E}_{\rho ,\sigma }[\mathcal{E}_{\sigma ,\tau _n}[X]]=\mathcal{E}%
_{\rho ,\tau _n}[X].  \label{e8.39}
\end{equation}
From (\ref{e8.30a}), we have
\[
|\mathcal{E}_{\rho ,\sigma }[\mathcal{E}_{\sigma ,\tau _n}[X]]-\mathcal{E}%
_{\rho ,\sigma }[\mathcal{E}_{\sigma ,\tau }[X]]|\leq CE[|\mathcal{E}%
_{\sigma ,\tau _n}[X]-\mathcal{E}_{\sigma ,\tau }[X]|^2].
\]
But by Lemma \ref{m8.11}, both $\{\mathcal{E}_{\sigma ,\tau
_n}[X]\}_{n=1}^\infty $ and $\{\mathcal{E}_{\rho ,\tau _n}[X]\}_{n=1}^\infty
$ are Cauchy sequences in $L^2(\mathcal{F}_T)$. We then can pass to the
limit on both sides of (\ref{e8.39}) to obtain (\ref{e8.35}). It is easy to
check that, once we have (A1)--(A3), (A4$_0$) and (A5), for $\rho $, $\sigma
$, $\tau \in \mathcal{S}_T$, the estimate (\ref{e8.a2}) still holds for $%
\sigma $, $\sigma ^{\prime }\in \mathcal{S}_T$ and $\tau $, $\tau ^{\prime
}\in \mathcal{S}_T$. From these estimates we have the continuity of $%
\mathcal{E}_{\sigma ,\tau }[\cdot ]$ in the following sense: for each $%
\sigma $, $\tau \in \mathcal{S}_T$, $X\in L^2(\mathcal{F}_\tau )$ and
sequences $\{\sigma _i\}_{i=1}^\infty $, $\{\tau _i\}_{i=1}^\infty $ in $%
\mathcal{S}_T$ such that $\sigma \leq \sigma _i\leq \tau $ and $\tau \leq
\tau _i$ with $\lim_{i\rightarrow \infty }\tau _i=\tau $ and $%
\lim_{i\rightarrow \infty }\sigma _i=\sigma $, we have
\[
\mathcal{E}_{\sigma ,\tau }[X]=\lim_{i\rightarrow \infty }\mathcal{E}%
_{\sigma _i,\tau }[X]=\lim_{i\rightarrow \infty
}\mathcal{E}_{\sigma ,\tau _i}[X],\;\hbox{in }L^2(\mathcal{F}_T).
\]
The uniqueness of $\mathcal{E}_{\sigma ,\tau }[\cdot ]$ is due to the
uniqueness part of Lemma \ref{m8.5a} and the continuity of $\mathcal{E}%
_{\sigma ,\tau }[\cdot ]$ in $\sigma $ and $\tau $. The proof is complete. $%
\Box $\medskip\

We also have the following optional stopping theorem.

\begin{theorem}
\label{m8.14}Let $Y\in D_{\mathcal{F}}^2(0,T)$ be an $\mathcal{E}$%
--supermartingale (resp. $\mathcal{E}$--submartingale). Then for each $%
\sigma $, $\tau \in \mathcal{S}_T$ such that $\sigma \leq \tau $, we have
\begin{equation}
\mathcal{E}_{\sigma ,\tau }[Y_\tau ]\leq Y_\sigma \hbox{ (resp.
}\geq Y_\sigma \hbox{), a.s. .}  \label{e8.40}
\end{equation}
\end{theorem}

\smallskip\noindent\textbf{Proof. }We only prove the supermartingale part. We first
consider the case $\sigma \in \mathcal{S}_T$ and $\tau \in \mathcal{S}_T^0$.
Let $\{\sigma _n\}_{n=1}^\infty $ be a sequences in $\mathcal{S}_T^0$ such
that $\sigma \leq \sigma _n\leq \tau $, $n=1,2,\cdots $ and $%
\lim_{n\rightarrow \infty }\sigma _n=\sigma $, a.s.. By Lemma \ref{m8.5},
\[
\mathcal{E}_{\sigma _n,\tau }[Y_\tau ]\leq Y_{\sigma _n}\hbox{,
a.s.}
\]
From the convergence result Lemma \ref{m8.9}, the left hand side converges
to $\mathcal{E}_{\sigma ,\tau }[Y_\tau ]$ in $L^2(\mathcal{F}_T)$. Since $%
Y\in D_{\mathcal{F}}^2(0,T)$, $Y_{\sigma _n}\rightarrow Y_\sigma $, a.s.. We
then have proved (\ref{e8.40}) for the case $\sigma \in \mathcal{S}_T$ and $%
\tau \in \mathcal{S}_T^0$.

Now let $\sigma $, $\tau \in \mathcal{S}_T$ and let $\{\tau
_n\}_{n=1}^\infty $ be a sequences in $\mathcal{S}_T^0$ such that $\tau \leq
\tau _n$, $n=1,2,\cdots $ and $\lim_{n\rightarrow \infty }\tau _n=\tau $,
a.s.. We have proved that
\[
\mathcal{E}_{\sigma ,\tau _n}[Y_{\tau _n}]\leq Y_\sigma \hbox{,
a.s.}
\]

It is clear that $Y_{\tau _n}\rightarrow Y_\tau $ a.s. We also have $%
|Y_{\tau _n}|\leq \sup_{0\leq t\leq T}|Y_t|$. Since the right hand side is
in $L^2(\mathcal{F}_T)$, it follows by Lebesgue's dominated convergence
theorem that $Y_{\tau _n}\rightarrow Y_\tau $ in $L^2(\mathcal{F}_T)$. Since
we have
\begin{eqnarray*}
&&\ |\mathcal{E}_{\sigma ,\tau _n}[Y_{\tau _n}]-\mathcal{E}_{\sigma ,\tau
}[Y_\tau ]| \\
\ &\leq &|\mathcal{E}_{\sigma ,\tau _n}[Y_{\tau _n}]-\mathcal{E}_{\sigma
,\tau _n}[Y_\tau ]|+|\mathcal{E}_{\sigma ,\tau _n}[Y_\tau ]-\mathcal{E}%
_{\sigma ,\tau }[Y_\tau ]|.
\end{eqnarray*}
The second term converges to zero in $L^2(\mathcal{F}_T)$. For the first
term, we apply (\ref{e8.30a}). It follows that
\[
E[|\mathcal{E}_{\sigma ,\tau _n}[Y_{\tau _n}]-\mathcal{E}_{\sigma ,\tau
_n}[Y_\tau ]|^2]\leq CE[|Y_{\tau _n}-Y_\tau |^2]\rightarrow 0.
\]
We finally have (\ref{e8.40}) for the general situation. $\Box $\medskip\

\end{document}